\documentclass[a4paper, 12pt]{amsart}
\usepackage{graphicx} 
\usepackage{float, appendix}
\usepackage{cancel, comment}
\graphicspath{{plots}}
\usepackage{amsmath, amssymb, amsfonts,amsthm, xcolor,comment,enumitem, mathrsfs, amsaddr}

\newtheorem{theorem}{Theorem}
\newtheorem*{theorem*}{Theorem}

\newtheorem{lemma}{Lemma}
\newtheorem{definition}{Definition}
\newtheorem{proposition}{Proposition}
\newtheorem{assumption}{Assumption}

\newcommand{\mg}{\mathcal{G}}
\newcommand{\Prob}{\mathbb{P}}

\usepackage[a4paper, total={6in, 8in}]{geometry}

\newcommand{\poubelle}[1]{}

\newcommand{\R}{{\mathbb{R}}}

\newcommand{\E}{\mathbb{E}}
\renewcommand{\P}{\mathbb{P}}
\newcommand{\Var}{\mathbb{V}{\rm ar}\,}
\newcommand{\Cov}{\mathbb{C}{\rm ov}\,}
\newcommand{\mx}{\mathcal{M}}

\newcommand{\argmax}{\mbox{argmax}}
\newcommand{\argmin}{\mbox{argmin}}
\newcommand{\1}{{{\bf 1}}}

\usepackage{hyperref}

\usepackage{url}

\newcommand{\KS}{\textsf{KS}}
\newcommand{\supp}{\textsf{supp}}

\makeatletter
\renewcommand{\eqref}[1]{(\ref{#1})}
\makeatother


\title[Fundamental limits for empirical approx. of tilted distributions]{Fundamental limits for weighted empirical approximations of tilted distributions}

\author{Sarvesh Ravichandran Iyer$^{1}$, Himadri Mandal$^2$, Dhruman Gupta$^{1}$, Rushil Gupta$^{1}$, Agniv Bandyopadhyay$^3$, Achal Bassamboo$^4$, Varun Gupta$^5$, Sandeep Juneja$^{1}$}
\thanks{${1} : $ Ashoka University. $2 : $ Indian Statistical Institute, Kolkata. $3 : $ Tata Institute of Fundamental Research, $4 : $ Northwestern University, $5 : $ The University of Utah}
\thanks{Note : The three faculty names are listed at the end and arranged in alphabetical order by last name. The first author is the corresponding author.}

\begin{document}

\begin{abstract}
Consider the task of generating samples 
from a tilted distribution of  a random vector whose underlying distribution is unknown, but samples
from it are available. This finds applications in fields such as finance and climate science, and in rare event simulation. In this article, we discuss the asymptotic efficiency of a self-normalized importance sampler of the tilted distribution. We provide a sharp characterization of its accuracy, given the number of samples and the degree of tilt. Our findings reveal a surprising dichotomy: 
while the number of samples needed
to accurately tilt a bounded random vector increases polynomially 
in the tilt amount, it increases at a super polynomial rate for unbounded distributions. 
\end{abstract}

\maketitle

\section{INTRODUCTION}

Exponential tilting of random variables or vectors is a technique in the field of Monte Carlo simulation that relies upon changing the underlying measure of a probability space, making a desired set of outcomes easier to sample from. It has wide applications ranging from finance to climate science (\cite{ragone2018computation}, \cite{McLeishMen2015}). An excellent review of exponential tilting may be found in \cite[Chapter 1]{Alvo2022}. For a historical account of exponential tilting in the context of event estimation, the reader is referred to \cite{JUNEJA2006291} and \cite[Chapter 6]{asmussen2007stochastic}. Beyond exponential tilting, there are other forms of tilting, including hazard rate tilting (see \cite{juneja2002simulating}) and related polynomial tilting (see \eqref{twisted} for a definition). We shall refer to these techniques as \emph{tilting} throughout our article.

Recent work has highlighted the importance of data-driven event estimation. For example, in generative modeling, one may consider the problem of generating tilted samples of an unknown distribution (see \cite{wang2024proteinconformation}). We ask the following question: given iid samples $X_1,X_2,\ldots,X_n$ of a random vector $X$ with unknown distribution, can we produce empirical samples from a tilted version of $X$? Furthermore, for an accurate representation, how would the number of samples depend on the tilt? 

To do this, we first construct a natural estimator using empirical reasoning (see \eqref{empirical}), which is a reweighed version of the empirical estimator for the samples $X_1,\ldots,X_n$. Our main results are matching upper and lower bounds for the asymptotic KS distance of this estimator in the regimes where the tilt is fixed, and where it increases to infinity. In particular, we are able to characterize asymptotic accuracy of the estimator in terms of the \emph{coefficient of variation}, see \eqref{m2tmt2}. 

Such characterizations are novel in sampling literature, since to the best of our knowledge, there are no results characterizing tilting thresholds that demarcate the regions where empirical estimation is accurate and where it is not.

In mathematical finance, Black and Litterman \cite{BlackLitterman1990_AssetAllocation,BlackLitterman1992_GlobalPortfolioOptimization} pioneered a model for financial securities which was motivated by finding a model in proximity to the physical model which satisfies certain expert-suggested views that are expressed as constraints on the expectations of a linear function of the underlying random variables. Technically, this corresponds to finding a probability measure that minimizes an entropic divergence with respect to a given probability measure and satisfies certain moment constraints. The optimal solution for the KL divergence is well known to correspond to a new measure that is an exponentially twisted version of the original measure (see, e.g., \cite{BuchenKelly1996_MaxEntropyOption},
\cite{meucci2010fully}, \cite{stutzer1996simple}).

Tilted distributions are also central to large deviations theory. Conditioned on sums of a large number of independent identically distributed random variables taking a large deviation, the limiting distribution of each individual component (as the number of variables tends to infinity) is a tilted version of the original distribution (See \cite{Dembo2010}).

Further motivation also arises from sampling literature. Indeed, the proposed tilting is an example of self-normalized importance sampling (SNIS) (see \cite[Chapter 9]{Owen2000_SafeEffectiveImportanceSampling}),
also known as ratio sampling (see \cite{Hesterberg1995_WeightedAverageImportanceSampling}) and weighted importance sampling (see \cite[Section 5.6]{Rubinstein1981_SimulationAndTheMonteCarloMethod}). There are multiple applications of SNIS, such as Bayesian data problems\cite{Kong1997_SequentialImputations} and policy selection\cite{KuzborskijVernadeGyorgySzepesvari2021_ConfidentOffPolicy}. The subsphere of exponential tilting is again motivated by a minimization problem involving moment constraints\cite{Clark1966_ExponentialTransform,FuhTengWang2013_EfficientImportanceSampling_arXiv}.

Theoretical results for SNIS are relatively limited (as remarked in \cite[Section 2.3]{CardosoSamsonovThinMoulinesOlsson2022_BR_SNIS}). For instance, upper bounds on the bias and variance of the estimator in the $L^1$ and $L^2$ distance can be found in \cite{AgapiouPapaspiliopoulosSanzAlonsoStuart2017_ImportanceSamplingIntrinsicDimension}. However, these works do not develop lower bounds as we do. Additionally, they do not use the KS distance, which is our focus.

Our answer to the central question lies in the following two heuristics : a large tilt cannot be performed with a small number of samples, and tilting unbounded random vectors is fundamentally harder than tilting bounded ones. We shall develop theoretical versions of the aforementioned heuristics and demonstrate their accuracy through simulations.

The central objects in this article are the distributions of the tilted random vectors, and their empirical approximations. If $X$ is a random vector on $\R^d, d \geq 1$, $g : \mathbb R^d \to \mathbb R^d$ is a suitable site-specific tilt, and $\theta \in \mathbb R^d$ is a vector such that $\mathbb E[e^{\theta^T g(X)}] < \infty$, then we denote by $X_{\theta}$ the tilted random vector, whose distribution is given by  
\begin{equation}\label{twisted}
\Prob[X_{\theta} \in A] = \frac{\E[e^{\theta^Tg(X)}\1_{X \in A}]}{\E[e^{\theta^T g(X)}]}.
\end{equation}

Given independent and identically distributed samples $X_1,X_2,\ldots,X_n$ of $X$, we wish to construct an estimator for $X_{\theta}$. We achieve this as follows : the unknown expectations in \eqref{twisted} can be replaced by their empirical estimators, leading to the ``reweighed'' empirical estimator $R_{n,\theta}$ given by 
\begin{equation}
\Prob[R_{n,\theta} \in A] =\frac{\sum_{i=1}^n e^{\theta^T g(X_i)} \1_{X_i \in A}}{\sum_{i=1}^n e^{\theta^T g(X_i)}} = \sum_{i=1}^n \frac{e^{\theta^T g(X_i)}}{\sum_{j=1}^n e^{\theta^T g(X_j)}}1_{X_i \in A}.\label{empirical}
\end{equation}

The motivation behind the word ``reweighed'' here is the observation that the support of $R_{n,\theta}$ is $\{X_1,\ldots,X_n\}$, which is identical to the support of the empirical estimator $A \mapsto \frac 1n \sum_{i=1}^n \1_{X_i \in A}$. However, the probability assignments are different, leading to the term ``reweighed'' (or self-normalized/weighted as in existing literature).

The quantity central to our investigation is given by \begin{equation}\label{m2tmt2}
M_{\theta} = \frac{\E[e^{2\theta^T g(X)}]}{\E[e^{\theta^T g(X)}]^2}.
\end{equation}
Observe that $M_{\theta} = 1+\textsf{cv}(e^{\theta^T g(X)})^2$, where $$\textsf{cv}(Y) = \frac{\sqrt{\Var(Y)}}{\E[Y]}$$ is the coefficient of variation of a random variable $Y$, which is often used to quantify the accuracy of empirical estimation of $Y$. This should not be surprising, in light of the observation that $R_{n,\theta}$ is obtained by replacing $\E[e^{\theta^T g(X)}]$ by its empirical mean.

To reiterate, in this article our objective is to delineate asymptotic regimes under which the reweighed empirical estimator $R_{n,\theta}$ accurately estimates $X_{\theta}$ as $(n,\theta)$ jointly vary. In particular, if $\KS$ denotes the Kolmogorov-Smirnov distance, then we achieve each of the following aims. 
\begin{enumerate}[label = (\alph*)]
\item For $d=1$, fixing $\theta \in \mathbb R$ we provide a rate of convergence of $\KS(R_{n,\theta},X_{\theta})$ to zero, and demonstrate that the scaled limit's expectation and fluctuations are controlled by $M_{\theta}$.
\item If $\theta_n, n \geq 1$ converges to infinity, we create a regime in terms of $M_{\theta_n}$ within which $\KS(R_{n, \theta_n}, X_{\theta_n})$ converges to zero with an \emph{a priori} rate of convergence.
\item If $X$ is a bounded random variable, we make tail assumptions on $X$ under which the above regime is maximally large, and comment on the asymptotic behavior of $R_{n, \theta_n}$ if $\KS(X_{\theta_n}, R_{n,\theta_n}) \not \to 0$. In particular, we prove that an accurate tilt by $\theta_n$ can be performed with a number of samples which is only \emph{polynomial} in $\theta_n$.
\item We repeat the above exercise for bounded random vectors in multiple dimensions and show that the number of samples required for accurate tilting is polynomial in $\theta_n$.
\item For unbounded random vectors, we prove that accurate tilting requires 
super-polynomially  many samples in $\theta_n$, in contrast to the bounded case. 
\end{enumerate}

We address each of these points in separate sections. Namely, in Section~\ref{fixedtheta} we address tilting by a fixed $\theta$. In Section~\ref{gen1D} we address tilted sampling as $\theta_n \to \infty$. In Section~\ref{1DWeibull} we establish tight regimes on the accuracy of tilting in $1$ dimension, followed by multiple dimensions in Section~\ref{mdweibull}. Finally, unbounded random vectors are discussed in Section~\ref{unbdd}.

In Appendix~\ref{appendix:a} we prove the results from Section~\ref{fixedtheta}. In Appendix~\ref{appendix:b} we prove the results from Section~\ref{gen1D}. In Appendix~\ref{appendix:c} we prove a multitude of results that apply to all the regimes considered in Sections~\ref{1DWeibull}, \ref{mdweibull} and \ref{unbdd}. Following this, in Appendices~\ref{appendix:d}, \ref{appendix:e} and \ref{appendix:f} we prove the results in 
Sections~\ref{1DWeibull}, \ref{mdweibull} and \ref{unbdd} respectively. Finally, Section~\ref{appendix:g} is reserved for experiments.

{\textbf{Notation}} : Whenever we speak of random vectors, we shall also include the one-dimensional case in our purview, and refer to a random vector as a random variable if it is one-dimensional. We always denote by $X$ the random vector which is being subject to tilting. $\theta$ denotes the vector direction/scalar along which the tilt is being affected. 

Throughout this article, $F_Y$ denotes the cumulative distribution function (CDF) of any untilted random vector $Y$. The CDF of $X_{\theta}$ will be denoted by $F_{\theta}$ and that of the empirical distribution $R_{n,\theta}$ will be denoted by $F_{n,\theta}$. Convergence in distribution will be denoted by $\overset{d}{\to}$, while convergence in any other notion will be denoted by $\to$ and mentioned explicitly.

Measures will be indicated by $\mu,\nu$, etc. and integration with respect to them will be indicated by $d\mu, d \nu$ (e.g. $\int_{0}^{\infty} f(x) d\nu(x)$). Sometimes the domain of integration may be omitted if it is the entire support of the measure.

The indicator function of a set $A$ will be denoted by $\1_{A}$, and its complement will be denoted by $A^c$. The topological boundary $\partial A$ of a Borel set $A \subset \mathbb R^d$ is the intersection of $A^c$ with the closure of $A$.

The nomenclature $C_1,C_2,\ldots$ is reserved for constants which are independent of the parameters in the result (typically $X,\theta$), while parameter-dependent quantities will be indicated using subscripts e.g. $M_{\theta}$. 

\section{TILTING FOR FIXED $\theta$}\label{fixedtheta}

In this section, we establish asymptotic consistency of the estimator \eqref{empirical} in the KS distance on the $\sqrt{n}$ scale, in one dimension. Further, we show that the average maximal fluctuations of the scaling limit are $\Theta(M_{\theta})$. 

Recall the true and empirical estimators given by \eqref{twisted} and \eqref{empirical}, respectively, for a random vector $X$ and vector $\theta$. In this section, we assume that $d=1$ and that $X$ is a \emph{continuous} random variable. Furthermore, let $n \geq 1$ be a positive integer, $\theta\in \mathbb R$ be a \emph{fixed} positive real number and $g : \mathbb R \to \mathbb R$ be a strictly continuous increasing function on $\mathbb R$. We will let $n$ tend to infinity in what follows, and assume that $\E[e^{\theta g(X)}] < \infty$.

For the sake of clarity, we redefine the one-dimensional versions of \eqref{twisted} and \eqref{empirical}. That is,
\begin{equation}\label{twisted1D}
\mathbb P(X_\theta \in A) = \frac{\E[e^{\theta g(X)}\1_{X \in A}]}{\mathbb E[e^{\theta g(X)}]},
\end{equation}
and if $X_1,X_2,\ldots,X_n$ are iid with distribution $X$, then \begin{equation}\label{empirical1D}
\mathbb P(R_{n,\theta} \in A) = \frac{\sum_{i=1}^n e^{\theta g(X_i)}\1_{X_i \in A}}{\sum_{i=1}^n e^{\theta g(X_i)}}.
\end{equation}

Note that $F_{n, \theta}(x)$, the CDF of $R_{n, \theta}$ is a random CDF since $R_{n,\theta}$ depends upon the random samples $X_i$. Let $F_{\theta}$ be the CDF of $X_{\theta}$. Define the Kolmogorov-Smirnov ($\KS$) distance between two random variables $U,V$ with CDFs $F_U,F_V$ by 
\begin{equation}\label{KS}
\KS(U,V) = \sup_{x \in \mathbb R} |F_U(x) - F_V(x)|.
\end{equation}

Some comments about the choice of distance follow. Any estimator is typically assessed through the lens of consistency and asymptotic normality, before any further study is conducted. This entails convergence in distribution, which is equivalent to convergence in the KS distance for continuous limits (see \cite[page 3]{Resnick1987}). As a result, the KS distance presents itself naturally in this setup. Furthermore, the KS distance is preserved by scaling limits, as we demonstrate in Lemma~\ref{KStoSL}. 

Thus, the KS distance is fundamental to our analysis. Analysis in any other distance, such as the p-Wasserstein or quantile-based probability metrics, require distinctly different methods which we defer to future work.

Our first main result reads as follows. Recall that a Gaussian field $\mathcal G$ on $\mathbb R$ is an $\mathbb R$-indexed stochastic process $\mathcal G = \{\mathcal G(x) : x \in \mathbb R\}$ such that for every $x_1,x_2,\ldots,x_m \in \mathbb R, m \geq 1$, the random vector $\{\mathcal G(x_1),\ldots,\mathcal G(x_m)\}$ is Gaussian.

\begin{theorem}\label{thm:bd}
There exists a Gaussian random field $\mathcal G_{\theta}$ on $\mathbb R$ such that $$\sqrt{n}\KS(F_{n,\theta}, F_{\theta}) \xrightarrow{d} Z,$$ 
where $Z = \sup_{x}|\mathcal G_{\theta}(x)|$. Furthermore, there exist constants $C_1,C_2>0$ independent of $X$ and $\theta$ such that 
$$
C_1 \sqrt{M_{\theta}} \leq \mathbb E[Z] \leq C_2 \sqrt{M_{\theta}},
$$
where $M_{\theta}$ is as in \eqref{m2tmt2}. Furthermore, 
\begin{equation}\label{dekhna}
\mathbb P(|Z - \E[Z]| > u) \leq e^{-u^2/M_{\theta}}
\end{equation}
for all $u>0$.
\end{theorem}

This theorem has multiple layers to it, which we uncover now. First, note that the appropriate rate at which the $\text{KS}$ distance converges is $\sqrt{n}$, which is independent of $\theta$. However, as we shall observe (see Proposition~\ref{lem:gauss} for the precise result), the covariance functional of $\mathcal G_{\theta}$ depends on $\theta$.

Studying the supremum $Z$ of $|G_{\theta}|$ is essential for demonstrating asymptotic sample complexity bounds, and the theorem states that $Z$ is strongly concentrated about its expectation whenever $M_{\theta}$ is small, demonstrating the cruciality of this quantity to establishing the accuracy of sampling. Furthermore, that the expectation of $Z$ itself is of the order of $M_{\theta}$ implies that the average maximal fluctuation in the limit is $\Theta(M_{\theta})$ with high probability.

Next, We shall provide a brief technical overview of the proof of this theorem. Since $F_{n,\theta}$ is obtained by replacing $\E[e^{\theta^T g(X)}]$ with its empirical estimator, it follows by the law of large numbers and the continuous mapping theorem that $\textsf{KS}(F_{n,\theta},F_{\theta}) \to 0$ almost surely. To understand the fluctuations in this convergence we use the delta-method, which combines the continuous mapping theorem with the Donsker Central-Limit-Theorem. Since the random elements $F_{n,\theta},F_{\theta}$ are infinite-dimensional, we require fairly sophisicated machinery to establish these results, which can be found in \cite[Section 3.9]{vandervaart1996weak}. This helps us establish the existence of a Gaussian process in the scaling limit. Finally, another application of the continuous mapping theorem brings us the limit $Z$ in the theorem above.

In order to establish bounds on $\E[Z]$, we make use of well-established upper and lower bounds on the expected suprema of Gaussian processes, in particular the Sudakov-Fernique theorem (\cite[Theorem A.2.5]{vandervaart1996weak}) for the upper bound and Sudakov's inequality (\cite[Corollary 2.2.8]{vandervaart1996weak}) for the lower bound. However, these results require estimation of the covering numbers of a metric space associated to $G_{\theta}$. In particular, on $\R$ define 
$$
d(s,t) = \sqrt{\E[|G_{\theta}(s) - G_{\theta}(t)|^2]} = \|\mathcal G_{\theta}(s) - \mathcal G_{\theta}(t)\|_{2}
$$
where $\|\cdot\|_2$ denotes the $L^2$ distance between random variables. For any $\epsilon>0$, the covering number $N(\epsilon)$ is the smallest size of an $\epsilon$-cover i.e. a set such that every real number is at most $\epsilon$ away from some point in this set, in the distance $d$. 

We develop estimates for $N(\epsilon)$ in $(\R,d)$ in terms of its diameter, and then show that the diameter is controlled by $\sqrt{M_{\theta}}$. This, along with the Sudakov-Fernique theorem and Sudakov inequality, completes the proof of the expectation bounds.

The concentration bound is a standard consequence of the Borell-TIS inequality (Theorem~\ref{btis}). For the detailed proof of Theorem~\ref{thm:bd}, the reader is referred to Appendix~\ref{appendix:a}.

\section{ACCURATE TILTING IN ONE DIMENSION}\label{gen1D}

In this section, let $X$ be any random variable. Unlike the previous section, we allow $X$ to contain atoms. Let $g: \R \to \R$ be a measurable site-specific tilt function. Assume that $\E[e^{\eta g(X)}] < \infty$ for all $\eta>0$, which is certainly true if $g(X)$ is bounded. As a technical condition, we assume that $g(X)>0$ with non-zero probability, which ensures that $\lim_{\theta_n \to \infty} \E[e^{\theta_n g(X)}] = +\infty$.

Let $n$ be a positive integer and $\theta$ be a positive real number. In this section, $\theta$ and $n$ will be made to tend to infinity jointly. Recall \eqref{twisted1D} and \eqref{empirical1D}, the distributions of $X_{\theta}$ and $R_{n,\theta}$ respectively, and their associated CDFs $F_{\theta}$ and $F_{n,\theta}$. 

Our main result for this section reads as follows.

\begin{theorem}\label{thm:genbd}
Let $\theta_n, n \geq 1$ be a sequence. If $$\lim_{n \to \infty}\frac{M_{\theta_n}}{n} = 0,$$ then
$s_n\KS(F_{n,\theta_n}(x),F_{\theta_n}(x)) \overset{d}{\to} 0$, 
where $s_n$ is any sequence such that \begin{equation}\label{decay}
\lim_{n \to \infty}s_n^2\frac{M_{\theta_n}}{n} =0.
\end{equation}
\end{theorem}

That is, provided that $\theta_n$ doesn't grow too fast, it is possible to asymptotically accurately sample from $X_{\theta}$ with an explicit bound on the $\KS$ accuracy. Note that this includes the previous case where $\theta$ is fixed and $g(x)=x$, which implies that $s_n = o(\sqrt{n})$ works. However, no scaling limit is provided, hence this result is weaker than Theorem~\ref{thm:bd} for this specific case.

An illustrative example of this result is $X = Exp(1)$, with $\theta_n = \frac 12 - \frac 1{\sqrt{n}}$ and site-specific tilt $g(x)=x$. One computes that $$M_{\theta_n} = \frac{(\sqrt{n}+2)^2}{8\sqrt{n}} = \Theta(\sqrt{n}).$$ Therefore, $s_n^2\frac{M_{\theta_n}}{n} \to 0$ provided that $s_n = o\left(n^{\frac 14}\right)$. This is particularly illustrative in light of the fact that $\theta_n \to +\infty$ is not necessary.

Another application of the lemma is when $X$ is a discrete random variable. Let $X$ be uniform on the set $\{1,2,\ldots,6\}$ i.e. $X$ is the result from the roll of a fair dice. Let $\theta_n = C \log n$ for some $C< \frac 16$, $g(x)=x$, and note that $$\frac{M_{\theta_n}}{n} = \frac{\E[e^{\theta_n X}]}{n} = \frac 16\sum_{i=1}^6 e^{\theta_n i - \log n} = \frac 16 \sum_{i=1}^6 n^{Ci-1},$$
which implies that
$$\frac{M_{\theta_n}}{n} \to 0 ,\quad\frac{\frac{M_{\theta_n}}{n}}{n^{6C-1}} = \frac 16.$$
In particular, if $s_n = o(n^{\frac{1-6C}{2}})$ then $s_n^2 \frac{M_{\theta_n}}{n} \to 0$. We remark that this particular application cannot be replicated by any of the other results, whose tail or continuity assumptions would debar $X$.

We shall now provide the following rough idea of the positive bound. Note that the numerator of \eqref{twisted1D} is similar to the numerator of \eqref{empirical1D}, and likewise the denominators are expected to be close for large $n$. Thus, the two CDFs should also be close to each other, once we combine these two observations. More rigorously, we will prove the following result.

\begin{proposition}\label{prop:numden}
For any $t>0$, $n \geq 1$ and $\theta_n \in \mathbb R$, let $T_n = \frac{t}{2}\mathbb E[e^{\theta_n g(X)}]$. Then, we have
\begin{align*}
&\mathbb P\left[\sup_{x} |F_{n,\theta_n}(x) - F_{\theta_n}(x)| \geq t\right] \\ \leq & \mathbb P\left[\sup_{x}\left|\frac 1n \sum_{i=1}^n e^{\theta_n g(X_{i})}\1_{X_{i} \leq x}-\mathbb E[e^{\theta_n g(X)} \1_{X \leq x}]\right| \geq T_n\right]\\ + &\mathbb P \left[\left|\frac 1n \sum_{i=1}^n e^{\theta_n g(X_{i})}-\mathbb E[e^{\theta_n g(X)}]\right| \geq T_n\right].
\end{align*}
\end{proposition}

While the second term above is rather easily bounded, the first term requires a specific concentration result from empirical process theory, which we proceed to apply in our result. We will prove the following result using \cite[Corollary 3.1]{1bc2448d-c21e-3f24-8a81-d9751180f50d}.

\begin{theorem}\label{ledvdg}
Let $Y_1,\ldots,Y_n$ be iid random variables with distribution $Y$ such that $\E[Y^2]<\infty$. Let $$
Z_n = \sup_x \left|\frac 1n\sum_{i=1}^n Y_i 1_{Y_i \leq x} - \mathbb E[Y 1_{Y \leq x}]\right|.
$$
Then, for every $\epsilon, x>0$,
$$
\mathbb P\left(Z_n \geq (1+\epsilon)\E[Z_n]+x\right) \leq \frac{c_{\epsilon}\sqrt{M}}{x \sqrt n}
$$
where $M = \E[Y^2]$.
\end{theorem}

The proof of Theorem~\ref{thm:genbd} follows by treating the second term on the right hand side of Proposition~\ref{prop:numden} using Chebyshev's inequality, and the first term using Theorem~\ref{ledvdg}. We relegate its proof to Appendix~\ref{appendix:b}.

\section{TILTING IN THE 1D WEIBULL REGIME}\label{1DWeibull}

In this section, we discuss the asymptotic accuracy of the estimator \eqref{empirical1D}. In particular, we demonstrate that accuracy of one-dimensional tilting can be completely disseminated, with the results looking strikingly simple in comparison to corresponding results in multiple dimensions. Furthermore, our assumptions in one dimension motivate analogous assumptions in high dimensions.

While the previous section creates a regime in which tilting can be performed accurately, it does not contain any negative results. On the other hand, we saw that Theorem~\ref{thm:bd} provided a tight rate of empirical convergence for the estimator, and this was achieved through a scaling limit. 

Therefore, throughout this section, let $X$ be an \emph{upper bounded} random variable, and $$\mx = \sup\{x : \mathbb P(X \leq x) <1\}$$ be the supremum of the support of $X$. As usual, let $\theta_n, n\geq 1$ be any sequence, and recall $M_{\theta_n}$ from \eqref{m2tmt2}. We shall, under tail assumptions on $X$, establish results on the behavior of $\KS(R_{n, \theta_n}, X_{\theta_n})$ as $n \to \infty$, in the following three regimes:
\begin{enumerate}[label = (\alph*)]
\item  $\frac{M_{\theta_n}}{n} \to 0$.
\item  $\frac{M_{\theta_n}}{n} \to c$ for some $c \in (0,\infty)$.
\item  $\frac{M_{\theta_n}}{n} \to \infty$.
\end{enumerate}

We begin this section with a brief discussion on the importance of scaling limits in our analysis, followed by the results. As remarked earlier, the KS distance is intrinsically connected to scaling limits. In particular, the following result holds. 
\begin{lemma}\label{KStoSL}
Let $X_{1,n},X_{2,n}$ be two sequences of random variables such that for some sequences $a_n$ and $b_n>0$, we have 
\begin{equation}
\frac{X_{i,n} - a_n}{b_n} \overset{d}{\to} Z_i,
\end{equation}
for some random variables $Z_1,Z_2$. 
\begin{enumerate}[label = (\alph*)]
\item If $Z_1\overset{d}{=} Z_2$, then $\KS(X_{1,n},X_{2,n}) \to 0$ if $Z_1$ is continuous. \item If $Z_1(x) \overset{d}{\neq} Z_2(x)$ then $\KS(X_{1,n},X_{2,n}) \not \to 0$ (even if one of $Z_1,Z_2$ has atoms).
\end{enumerate}
\end{lemma}

Scaling limits in the context of tilting were studied in \cite{BKR1} and \cite{BKR2003} (see also the references therein). Chiefly, they proved that scaling limits of exponential tilting could only arise from the {\em extended} Gamma family, which includes the Gamma (and thus exponential) and normal random variables up to parametrization; subsequently the domains of attraction of these limits were also established. Note that these papers do not focus on the empirical estimator \eqref{empirical1D} at all.

It can be seen that the domains of attraction involved the notion of random variables/vectors possessing regularly varying CDFs. Thus, its appearance and importance to our work is fundamental.

In fact, given that exponential twisting by a positive amount biases towards larger samples, it is expected that regulated behavior of exponentially tilted samples is possible only when the underlying random variable behaves well at values near its maximum. This is particularly apt when one wishes to obtain asymptotic results, as is the case in our article. 

Therefore, our assumptions are also inspired by extreme value theory that will aid us in obtaining negative results. Indeed, as demonstrated in \cite[Chapters 1,5]{Resnick1987}, the domains of attraction of such random vectors also involve the notion of multi-variate regular variation. Furthermore,  tail assumptions have been seen in the context of event estimation for random walks\cite{BlanchetGlynn2008_RareEventHeavyTails}, apart from their natural origin in extreme value theory.

However, keeping in mind the applicability of our assumptions, we forego any complicated parameters and focus our attention on an easily defined family of distributions. Finally, we provide standalone proofs of facts already proved in the aforementioned articles, since we plan to reuse asymptotic analysis and insights where necessary.

As a result of Lemma~\ref{KStoSL}, we can focus our attention on scaling limits for $X_{\theta_n}$ and $R_{n,\theta_n}$. Note that as $\theta_n \to \infty$, samples which are closer to $\mx$ receive more weight in the empirical estimator. This observation contributes crucially to our final results. Given that exponential twisting biases towards larger samples, it is expected that regulated behavior of exponentially tilted samples is possible only when the underlying random variable is well-behaved at values near its maximum. This is particularly apt when one wishes to obtain asymptotic results, as is the case in our article. Therefore, we shall now make some assumptions on this behavior which are inspired by extreme value theory, that will aid us in obtaining negative results.

 The famous Fisher-Tippett-Gnedenko theorem states that $\max\{X_1,X_2,...,X_n\}$ possesses a limit only in three well-known scenarios in one dimension. These are known as the Weibull, Gumbel and Frechet limits, and we restrict our attention to the Weibull regime, since it pertains specifically to (upper) bounded random variables.

Before we define the Weibull regime, we require the definition of a regularly varying function (see \cite[Page 18, Section 2]{BinghamGoldieTeugels1987}). A function $f : \R \to \R_+$ is regularly varying at $0$ with index $\alpha \in \mathbb R$ if for all $u>0$, 
\begin{equation}\label{rv}
\lim_{t \to 0} \frac{f(tu)}{f(t)} = u^{\alpha}.
\end{equation}
An analogous definition holds for $f$ being regularly varying at (positive) infinity. As a canonical example, $x^{\alpha}$ is regularly varying of order $\alpha$ at $0$ and infinity, for any $\alpha \in \mathbb R$.

\begin{assumption}\label{weibull}
We assume that the random variable $X$ with CDF $F_X$ falls in the \emph{Weibull regime} with parameter $\alpha$ for some $\alpha>0$. That is, the function $$x \mapsto \mathbb P(X > \mx-x) = 1 - F_X(\mx-x)$$ is regularly varying at $0$ of order $\alpha$.
\end{assumption}

Note that $\alpha=0$ is not allowed : indeed, some bounded random variables which satisfy the above assumption but with $\alpha=0$ do not fall into any of the extreme value regimes. For example, if $\mathbb P(X=\mx) > 0$ then $X$ does not lie in the Weibull regime, since for any $u>0$, 
$$
\lim_{t \to 0} \frac{1-F_X(\mx-tu)}{1-F_X(\mx-u)} = \frac{\mathbb P(X=\mx)}{\mathbb P(X=\mx)} = 1 = u^0.
$$
Since $\alpha=0$ is debarred, $X$ does not satisfy our assumption. Note that we were nevertheless able to obtain results on tilting such a random variable (fair dice roll) in the previous section using Theorem~\ref{thm:genbd}.

Some examples are as follows. If $X=U[0,1]$ then $\mx=1$ and $\mathbb P(X>1-u) = u$ is regularly varying of order $1$ at $0$. Hence, a uniform random variable is in the Weibull regime with parameter $\alpha=1$. We include some more examples below, whose justifications will be provided in Appendix~\ref{appendix:d}.

\begin{lemma}\label{examples1DWeibull}
Each of the following random variables lies in the Weibull regime.
\begin{enumerate}[label = (\alph*)]
\item $Beta(a,b)$, for any $a,b>0$ are in the Weibull regime with parameter $\alpha = b$.
\item Truncated normal/exponential random variables are in the Weibull regime with parameter $\alpha =1$.
\item The sum of two (possibly dependent) random variables $X,Y$, where $X$ and $Y$ themselves lie in the Weibull regime. In fact if $X,Y$ have parameters $\alpha, \beta$ respectively then $X+Y$ has parameter $\alpha+\beta$.
\end{enumerate}
\end{lemma}

Note that it is not necessary that $X$ needs to be a continuous random variable, even in a neighborhood of $\mx$, since regularly varying functions do not necessarily need to be continuous. However, $X$ must be continuous at $\mx$, as we saw earlier.

As remarked before, the Weibull regime allows for the sample maximum to have a scaling limit (cf. \cite[Proposition 0.3]{Resnick1987}).
\begin{lemma}[Fischer-Tippett-Gnedenko]\label{lem:weibull}
Suppose that $X$ is in the Weibull regime, and $X_1,X_2,...,X_n$ are iid with distribution $X$. Then,
$$
\frac{\mx - \max\{X_1,X_2,...,X_n\}}{\mx - F_X^{-1}(1-\frac 1n)} \to -W_{\alpha}
$$
in distribution, where $W_{\alpha}$ is a Weibull random variable having the distribution $$G_{\alpha}(x) = \exp(-(1+\alpha x)^{-1/\alpha}) \text{ for } x<0.$$
\end{lemma}

We will also deal with general tilts in this section. Recall that $g : \R \to \R$ is a site-specific tilt function. Typical examples of $g$ include $g(x) = x$ (exponential tilting), $g(x) = x^{\alpha}$ for some $\alpha>0$, and $g(x) = C \log x$ for some $C>0$ (polynomial tilting). However, the latter form of tilt is rather mild,  and is therefore applied most often while tilting heavy-tailed random variables, which we do not consider.

Thus, we can assume that $g$ has a regularly varying tail at $\mx$ of positive index. However, observe that all the above functions are also strictly increasing and continuous. This motivates the assumption :
\begin{assumption}\label{twist}
    $g : \R \to \R$ is strictly increasing and continuous. Furthermore, $g(\mx) - g(\mx-x)$ is regularly varying at $0$ of index $\beta > 0$.
\end{assumption}

 The following result shows that if $X$ is in the Weibull regime, then $Y = g(X)$ is in the Weibull regime as well.

\begin{lemma}\label{weibtoweib1D}
    If $X$ is in the Weibull regime with index $\alpha>0$ and $g$ satisfies Assumption~\ref{twist} with parameter $\beta>0$, then $g(X)$ is in the Weibull regime with index $\frac{\alpha}{\beta}>0$.
\end{lemma}

We are now ready to state our main theorems of this section. Throughout the rest of this section, let $X$ be in the Weibull regime with parameter $\alpha>0$, maximum value $\mx$ and CDF $F_X$. Keeping Lemma~\ref{KStoSL} in mind, we know that scaling limits for $X_{\theta_n}$ and $R_{n,\theta_n}$ under the three regimes defined at the beginning of this section imply corresponding statements on the asymptotic $\KS$ accuracy of the empirical estimator. Since all the regimes depend upon the asymptotic behavior of $M_{\theta_n}$, our first result establishes \emph{polynomial growth} of $M_{\theta_n}$.

\begin{theorem}\label{thm:m2tmt21D}
We have $$
\left(1-F_X\left(\mx- \frac 1{\theta_n}\right)\right)M_{\theta_n} \to \frac{2^{-\alpha}}{\Gamma(1+\alpha)}.
$$
\end{theorem}

By Assumption~\ref{weibull} and a combination of \cite[Theorem 1.4.1 and Proposition 1.5.1]{BinghamGoldieTeugels1987}, it is clear that $M_{\theta_n}$ grows at most polynomially. This theorem is a simple corollary of the following lemma, which follows from Karamata's Tauberian theorem for regularly varying functions.

\begin{lemma}\label{lem:asy}
The sequence $\E[e^{\theta_n X}]$, as $\theta_n \to \infty$ satisfies $$
\frac{\E[e^{\theta_n X}]}{e^{\theta_n \mx} \left(1-F_X\left(\mx- \frac 1{\theta_n}\right)\right)} \to \Gamma(1+\alpha).
$$
\end{lemma}

The appearance of the $\Gamma$ function above hints at its appearance in the scaling limit for $X_{\theta}$ as well. While the analysis is subtler than the proof of Lemma~\ref{lem:asy}, our next result explicitly establishes the scaling limit for $X_{\theta}$ in this regime as a Gamma random variable.

\begin{theorem}\label{thm:sltrue1D}
Let $X$ be in the Weibull regime with parameter $\alpha>0$ and maximum value $\mx$. As $\theta_n \to \infty$, we have $$
\theta_n(\mx - X_{\theta_n}) \xrightarrow{d} \Gamma(\alpha,1),
$$
where $\Gamma(a,b)$ is the positive random variable with density $$f_{\Gamma(a,b)}(x) = \frac{b^{a}}{\Gamma(a)} x^{a-1}e^{-bx}, \quad x>0.$$
\end{theorem}

Finally, under the three regimes mentioned at the start of this section, we state the scaling limits of $R_{n , \theta_n}$. However, an explanation of the below results follow for clarity.

Observe that as $\theta \to \infty$, $X_{\theta} \to \mathcal{M}$ a.s. Therefore, as $\theta$ increases, more samples near the maximum are required to approximate $X_{\theta}$ accurately. If we are to look at how many samples should be near the maximum for accurate sampling, Theorem~\ref{thm:sltrue1D} asserts that one should consider the number of sample points that are within $\frac 1{\theta_n}$ of $\mx$, whose expectation is $n\mathbb P(X_i \geq \mx - 1/{\theta_n})$.

We shall now understand what happens when $\frac{M_{\theta_n}}{n} \to 0$. Indeed, as we have seen, $M_{\theta_n} \sim \theta_n^{\alpha}$ as $n \to \infty$. So, we have $$
n\mathbb P(X_i \geq \mx - 1/{\theta_n}) \approx n \theta_n^{-\alpha} \to \infty.
$$
In particular, enough points lie near the maximum. Furthermore, their weights $e^{-1} \leq e^{\theta_n(X_i-B)} \leq 1$ are bounded, implying that each contributes equally in the limit. Thus, sampling is asymptotically accurate.

On the other hand, suppose that $\frac{M_{\theta_n}}{n} \to \infty$. Then, $n = \Omega(\theta_n^{-\frac 1{\alpha}})$, and $$
n\mathbb P(X_i \geq \mx - 1/{\theta_n}) \approx n \theta_n^{-\alpha} \to 0.
$$ 
Thus, asymptotically, we expect no points to lie near the maximum, and sampling to be inaccurate. However, the weights $e^{\theta_n(X_i-B)}$ are still likely to be dictated by the largest values of $X_i$. Heuristics from extreme value theory imply that $e^{\theta_n (\max_i X_i - B)}$ is significantly larger than the other values (in fact, our proof technique relies on this observation), dictating that $R_{n,\theta_n}$ should behave like $\max_i X_i$ in this regime.

Finally, consider the intermediate "critical" regime $\frac{M_{\theta_n}}{n} \to K \in (0,\infty)$. In this case, observe that, asymptotically, the number of points close to the maximum is given by 
$$
n\mathbb P(X_i \geq \mx - 1/{\theta_n}) \approx n \theta_n^{-\alpha} \to C_K.
$$
for some $C_K \in (0,\infty)$ depending upon $K$. Thus, the number of points close to the maximum is a random variable; in fact it is Poisson with mean $C_K$. Since the locations of these various limiting maximal points affect the final limit, we see that the final limit is a random measure. It is, indeed, a functional of a limiting Poisson random measure (PRM) associated to the extremal process of the samples $X_i, i \geq 1$. We defer further discussion of these interesting phenomena to Appendix~\ref{appendix:c}, but refer the reader to \cite[Chapter 3]{Resnick1987} for more details on point processes.

\begin{theorem}\label{thm:slemp1D}
For a sequence $\theta_n, n \geq 1$ such that $\theta_n \to \infty$, 
    \begin{enumerate}[label = (\alph*)]
\item If $\frac{M_{\theta_n}}{n} \to 0$, then \begin{equation}
\theta_n(\mx - R_{n,\theta_n}) \to \Gamma(\alpha,1).
\end{equation}
\item If $\frac{M_{\theta_n}}{n} \to c \in (0,\infty)$, then \begin{equation}
\theta_n(\mx - R_{n,\theta_n}) \to Z_{c,PRM}
\end{equation}
for a random variable $Z_{c,PRM} \neq \Gamma(\alpha,1)$ which depends on $c_n$ and a limiting PRM.
\item If $\frac{M_{\theta_n}}{n} \to \infty$,
 then
 \begin{equation}
     \frac{(\mx - R_{n,\theta_n})}{\mx-F_X^{-1}(1-\frac 1n)} \to -W_{\alpha},
 \end{equation}
 where $W_{\alpha}$ is a Weibull random variable of index $\alpha$.
 \end{enumerate}
\end{theorem}

Thus, combining Lemma~\ref{KStoSL} and Theorems ~\ref{thm:m2tmt21D}, ~\ref{thm:sltrue1D} and ~\ref{thm:slemp1D}, we assert that for random variable in the Weibull domain, one can perform tilting accurately with polynomially many samples in the amount of tilt.

We illustrate these results with an example. Let $X= U[0,1]$ and $g(x) = x$, for which we know that $\alpha=1$ and $\mx=1$. Furthermore, $$1-F_X\left(1-\frac 1{\theta_n}\right) = \frac 1{\theta_n}$$ for any $\theta_n>1$. By Theorem~\ref{thm:m2tmt21D}, if $\theta_n \to \infty$ we have $$\frac{M_{\theta_n}}{\theta_n} \to \frac 12.$$
By Theorem~\ref{thm:sltrue1D}, $$
\theta_n(1-X_{\theta_n}) \overset{d}{\to} \Gamma(1,1) = Exp(1),
$$
where $Exp(1)$ is the exponential random variable with parameter $1$ and density $\lambda e^{-x}$ for $x>0$. Note that $\frac{M_{\theta_n}}{n} \to 0$ if and only if $\frac{\theta_n}{n} \to 0$. In this case, by Theorem~\ref{thm:slemp1D}(a) it follows that $$
\theta_n(1-R_{n,\theta_n}) \overset{d}{\to} \Gamma(1,1) = Exp(1).
$$
Thus, by Lemma~\ref{KStoSL}, 
$$
\textsf{KS}(X_{\theta_n}, R_{n,\theta_n}) \to 0
$$
provided that $n = \Omega(\theta_n)$. This shows that a substantial twist can be achieved with very few samples. If, on the other hand, $\frac{\theta_n}{n} \to c$ for some $c \in (0,\infty)$, then $\frac{M_{\theta_n}}{n} \to \frac c2$, and we have by Theorem~\ref{thm:slemp1D}(b) that $$
\theta_n(1-R_{n,\theta_n}) \overset{d}{\to} Z_{\frac c2,PRM}.
$$
This random variable will not be equal to $Exp(1)$, demonstrating that $n = \Theta(\theta_n)$ samples are not enough. Furthermore, if very few samples are considered and $\frac{\theta_n}{n} \to +\infty$, then $\frac{M_{\theta_n}}{n} \to +\infty$ as well. Note that $$1-F_X^{-1}\left(1-\frac 1n\right) = 1-\left(1-\frac 1n\right) = \frac 1n.$$ Therefore, 
$$
n(1-R_{n,\theta_n}) \to -W_{1},
$$
which implies that the scaling limits are totally different for $X = U[0,1]$. 

All the results in this section will be proved in Appendix~\ref{appendix:d}.

\section{TILTING IN HIGH DIMENSIONS}\label{mdweibull}

In this section, let $X$ be a random vector with bounded support on $\mathbb R^d, d \geq 1$. Let $\theta \in \mathbb R^d$ be a non-zero vector, $g : \mathbb R^d \to \mathbb R^d$ be a continuous site-specific tilt weight, and recall the tilted distributions $X_{\theta}$ and $R_{n,\theta}$ from \eqref{twisted} and \eqref{empirical} respectively. Define the support of $X, \textsf{supp}(X)$ as the largest closed set outside which $X$ takes values with probability $0$.

As in the one-dimensional setup, we wish to study accuracy of empirical approximation through scaling limits. However, the $\KS$ distance is fundamentally one-dimensional. The problem with considering an extension as in \cite{jacobs2025efficient}, is that scaling limits are not respected under this extension.

Instead, we consider a weaker notion of convergence, and prove its compatibility with scaling limits.

\begin{lemma}\label{KStoSLhd}
Let $X_{1,n},X_{2,n}$ be two sequences of random vectors such that for some sequences $a_n \in \mathbb R^d$ and $b_n >0$, we have 
\begin{equation}
\frac{X_{i,n} - a_n}{b_n} \overset{d}{\to} Z_i,
\end{equation}
for some random vectors $Z_1,Z_2$. 
\begin{enumerate}[label = \alph*]
\item If $Z_1\overset{d}{=} Z_2$ and $Z_1$ is continuous, then $\sup_{x \in \mathcal{R}} |F_{X_1}(x)-F_{X_2}(x)| \to 0$ for every compact $\mathcal{R} \subset \mathbb R^d$. 
\item If $Z_1(x) \overset{d}{\neq} Z_2(x)$ then for some compact $\mathcal{R} \subset \mathbb R^d$ we have $\sup_{x \in \mathcal{R}} |F_{X_1}(x)-F_{X_2}(x)| \not \to 0$. 
\end{enumerate}
\end{lemma}

To begin the discussion of tilting in multiple dimensions, note that we must define a notion of convergence of $\theta$ to infinity in order to discuss asymptotic accuracy.  A rather natural choice would be to fix a unit vector $\theta \in \mathbb R^d, \|\theta\|_2= 1$, and consider a large parameter $c>0$. We can thus look at scaling limits of $X_{c_n\theta}$ and $R_{n, c_n\theta}$ as $c_n$ converges to infinity. This would correspond to tilting towards samples $x \in \textsf{supp}(X)$ such that $\theta^Tg(x)$ is relatively large i.e. samples which, after a transformation, point in a direction close to $\theta$.

Note that if $c_n \to \infty$, then we expect $X_{c_n \theta}$ to concentrate in regions where $\theta^T g(x)$ is maximized. Suppose there are multiple such points : then the analysis for the limit can get rather complicated and interesting. We impose the following assumption which ensures that a scaling limit is easier to obtain.
\begin{assumption}\label{ass:conv}
The following assumptions hold.
\begin{enumerate}[label = (\alph*)]
\item The functional $x \to \theta^T g(x)$ has a unique maximizer $x_{\theta}$ in $\supp(X)$.
\item For every $\epsilon>0$, $\mathbb P[\theta^T(g(X) - g(x_{\theta})) < \epsilon] > 0$.
\end{enumerate}
\end{assumption}

The assumption ensures that as $c_n \to \infty$, $X_{c_n \theta}$ concentrates around $x_{\theta}$, and not a multitude of different values. Part (b) of the assumption ensures that $x_{\theta}$ is not an isolated point in $\supp(X)$.

This assumption is naturally satisfied in many multidimensional settings with $g(x)=x$. For instance, any random vector whose support is a convex polytope satisfies the above assumption for all directions $\theta$. These include Dirichlet random vectors and uniform random vectors on any such shape, for instance.

Our first result establishes a "law of large numbers" for tilting.

\begin{theorem}
Let $X$ be bounded and satisfy Assumption~\ref{ass:conv} in some direction $\theta$. Then, 
$$
X_{c_n\theta} \xrightarrow{d} x_{\theta}
$$
as $c_n \to \infty$.
\end{theorem}
\begin{proof}
It is enough to prove that if $A$ is a subset of the support of $X$ such that the closure of $A$ doesn't contain $x_\theta$, then $\mathbb P(X_{c_n\theta} \in A) \to 0$. Without loss of generality, we assume that $A$ is closed (or we can replace it by its closure below).

Since $A$ is closed and doesn't contain $x_{\theta}$, by the continuity of the function $\theta^Tg$, there exists $\delta>0$ such that $\theta^Tg(y) < \theta^Tg(x_{\theta})-\delta$ for all $x \in A$. We write
\begin{align}
\mathbb P(X_{c_n\theta} \in A) = &\frac{\mathbb E[e^{c_n\theta^Tg(X)} \1_{X \in A}]}{\mathbb E[e^{c_n\theta^Tg(X)}]}\nonumber \\= &\frac{\mathbb E[e^{c_n\theta^T(-g(x_{\theta})+g(X)) + c_n \delta/2} \1_{X \in A}]}{\mathbb E[e^{c_n\theta^T(-g(x_{\theta}) + g(X)) + c_n \delta/2}]}.\label{One}
\end{align}
Then,
$$
\mathbb E[e^{c_n\theta^T(-g(x_{\theta}) + g(X)) + c_n \delta/2} \1_{X \in A}] \leq e^{-c_n\delta/2}\mathbb E[\1_{A}] \to 0 
$$
as $c_n$ converges to infinity. It now remains to confirm that the denominator of \eqref{One} is non-zero in the limit. Let $$B = \{x \in \supp(X) : \theta^T(g(x) - g(x_{\theta})) < \delta/4\}.$$ We know that $\mathbb P(X \in B)>0$. However,
$$
\mathbb E[e^{c_n\theta^T(-g(x_{\theta}) + g(X)) + c_n \delta/2}] \geq \mathbb E[e^{c \delta/4} \1_{B}]>0
$$
regardless of $c_n$. Thus, it follows that $\mathbb P(X_{c_n\theta} \in A) \to 0$. This completes the proof.
\end{proof}

We shall now proceed to the scaling limit. For this, we require a few additional assumptions on $g$ and the `tail' of $X$ at $x_{\theta}$. This, we address through a multivariate regularly-varying framework (see \cite[Section 2.3]{cones} and \cite[Chapter 5]{Resnick1987}). 

\begin{assumption}\label{ass:mvrv}
We assume that $X$ is regularly varying at $x_{\theta}$ i.e. there exists a non-degenerate Radon measure $\nu$ on $\{y \in \R^d : \theta^T y > 0\}$ such that $$
\frac{1}{U(t)}\int f d \mathbb P\left(\frac{x_{\theta}-X}{t}\right) \to \int f d \nu
$$
for all compactly supported $f$ on $\{\theta^Ty > 0\}$ which are $\nu$-a.s. continuous, and some function $U(t)$ which is regularly varying at $0$ with index $\alpha>0$. 

Furthermore, we assume that \begin{equation}\label{intcond}\nu(\{y : 0<\theta^T y \leq 1\}) \in (0,\infty).\end{equation} This is equivalent to saying that with non-trivial probability, empirical samples will approach $x_{\theta}$ along asymptotic directions which are concentrated around $\theta$, from the support of $X$.

Note that vague convergence is equivalent to $$
\frac{1}{U(t)} \mathbb P\left(\frac{x_{\theta}-X}{t}\in A\right) \to \nu(A)
$$
for all compact $A$ such that $\nu(\partial A) = 0$ where $\partial A$ is the topological boundary of $A$. Using approximation arguments, it is sufficient that the above convergence holds for all $A$ belonging to a determining class, such as all compact hypercubes or balls in $\R^d$.
\end{assumption}

The manner of introduction of $\nu$ necessitates a discussion on applicability and examples. Recall that in one-dimension, the Weibull tail assumption attempted to capture the behaviour of the CDF near the maximum: in particular, that tail probabilities behaved well under scaling. In this case, $\nu$ captures the scaling transformation (see Lemma~\ref{lem:nu}(a)).

However, in one dimension one can only approach the maximum from the left, while in high dimensions, the support of $X$ is only constrained to be a subset of $\{\theta^Ty \leq \theta^T x_\theta\}$, which implies that the geometry of the support of $X$ in the vicinity of $x_\theta$ needs to be carefully studied. As we demonstrate in our examples, the measure $\nu$ roughly answers the following question : around the point $x_{\theta}$, how does the probability mass of $X$ concentrate?

Observe that our method of tilting involves an approach from the direction $c\theta$, towards the point $x_{\theta}$. As explained in the discussion in the previous section, the analysis of the empirical estimator relies upon the number of sample points $X_i$ such that $\theta^T X_i$ is significantly larger than the values of $\theta^T x$ rest of the population. The integrability condition \eqref{intcond} is tantamount to the assertion that the likelihood of such a maximizer lying in a direction near-orthogonal to $\theta$ with respect to $x_{\theta}$ is vanishingly low. In other words, most sample maximizers will indeed concentrate in the direction of approach $\theta$ towards $x_{\theta}$ rather than in an orthogonal direction.

The following facts about the limit measure $\nu$ are useful. While the first and second establish scaling and factorization properties, the third establishes a key integrability property that will be used to determine scaling limits.

\begin{lemma}\label{lem:nu}
Let $X$ satisfy Assumption~\ref{ass:mvrv}, and suppose $\theta,\nu$ and $\alpha$ are as in the assumption. Then, 
\begin{enumerate}[label = (\alph*)]
\item $\nu$ satisfies $\nu(cA) = c^{\alpha}\nu(A)$ for all compact sets $A$ such that $\nu(\partial A) = 0$, and $c>0$.
\item $\nu$ is a product measure in the following sense : there exists a finite measure $\mu$ on $\mathcal{S} = \{\theta^T y = 1\}$ such that for every $a,b>0$ and Borel $B\subset S$,
$$\nu\left(\left\{y : \theta^T y \in [a,b], \frac{y}{\theta^T y} \in B\right\}\right) = (b^{\alpha} - a^{\alpha})\mu(B).$$
\item We have 
$$
\int_{\{y : \theta^T y>0\}} e^{-\theta^T x} d \nu(x) < \infty.
$$
\end{enumerate}
\end{lemma}

Assumption~\ref{ass:mvrv} is rather abstract. This necessitates the need for examples, which are provided by the following lemma.

\begin{lemma}\label{examplesmdweibull}
 The following random vectors all satisfy Assumption~\ref{ass:mvrv}, except (d).
\begin{enumerate}[label = (\alph*)]
    \item If $X$ is a random \emph{variable} in the Weibull regime with parameter $\alpha>0$, CDF $F_X$ and maximum value $\mx$, then it satisfies Assumption~\ref{ass:mvrv} with $\nu([0,y]) = y^{\alpha}$ and $U(t) = 1-F_X(\mx-t)$.
    \item Bounded random vectors having independent components $X_i, 1 \leq i \leq d$ with maximum values $\mx_i, 1 \leq i \leq d$ that lie in the one-dimensional Weibull regime (see Assumption~\ref{weibull}) are regularly varying at $(\mx_1,\mx_2,\ldots,\mx_d)$.
        \item A uniform random variable on any polytope is regularly varying at any extreme point of the polytope.
        \item For any fixed $r>0$, consider either a uniform random variable $X = U(\{x \leq r\})$, or a truncated normal random vector $X = N(0,I)\1_{N(0,I) \leq r}$. Then, $X$ is regularly varying at any point $y \in \R^d$ with $\|y\| = r$, with limit measure $\nu$ being Lebesgue on $\{r^T y >0\}$. However, the resulting limit measure $\nu$ \emph{does not} satisfy the integrability condition \eqref{intcond}.
\end{enumerate}
\end{lemma}

Indeed, in example (c) above, $X_{c_n \theta}$ does not possess a scaling limit as $c_n \to \infty$. Therefore, our analysis cannot be applied to the study of such random vectors, and either an alternate distance, estimator or marginalization will have to be considered. Such a class of examples cannot be dismissed as merely pathological. We will defer the study of such phenomena to future work.

The assumption that we need to make on the site-specific tilt $g$ motivates itself as a high-dimensional version of the corresponding one-dimensional Assumption~\ref{twist}. In order to keep our presentation simple, we restrict ourselves to polynomials.

\begin{assumption}\label{rvhd}
We assume that $g$ is of the form $$
g(x_1,\ldots,x_d) = (x_1^{\alpha_1}, \ldots , x_d^{\alpha_d}),
$$
for some $\alpha_1,\ldots,\alpha_d>0$.
\end{assumption}

Exactly as in the one-dimensional case, general tilting preserves the tail assumption.
\begin{lemma}\label{weibtoweibhd}
If $X$ satisfies Assumption~\ref{ass:mvrv} at $x_{\theta}$ then $g(X)$ satisfies Assumption~\ref{ass:mvrv} at the point $g(x_{\theta})$.
\end{lemma}

The proof of this result is in Appendix~\ref{appendix:e}. From this point we assume that $X$ satisfies Assumptions~\ref{ass:conv} and \ref{ass:mvrv}, and that $U(t)$ is regularly varying of order $\alpha>0$. 

As in the previous section, we begin by considering the quantity $M_{c_n\theta}$ and its asymptotics as $c_n \to \infty$, followed by the scaling limit of $X_{c_n\theta}$. We remark that Assumption~\ref{ass:mvrv} makes the proofs rather easy.

\begin{theorem}\label{m2tmt2hd}
We have
$$
U\left(\frac 1{c_n}\right)M_{c_n \theta} \to \frac{2^{-\alpha}}{\int e^{-\theta^T y} \nu(dy)}.
$$
\end{theorem}

In particular, $M_{c_n \theta}$ grows at most polynomially with $c_n$. Next, we state a scaling limit for $X_{c_n \theta}$ as $c_n \to \infty$.

\begin{theorem}\label{sltruemd}
As $c_n \to \infty$, $$
c_n(x_{\theta}- X_{c_n \theta}) \to Z
$$
where $Z$ is a random variable whose density is proportional to $e^{-\theta^{T}y} \nu(dy)$.
\end{theorem}

We remark that by the scaling \eqref{scaling}, one can prove that $Z$ above has $\Gamma$ marginals, which shows that it extends the one-dimensional Theorem~\ref{thm:m2tmt21D}. However, the marginals may be dependent on each other based on the nature of $\nu$.

Finally, we can obtain asymptotics for $R_{n,c_n\theta}$ in the same three regimes as in the previous section. The first and second regimes intuitively follow as in the previous section, (with the second regime depending upon a PRM).

However, for rapidly growing $\theta_n$, a glance at \eqref{empirical} reveals that $R_{n,c_n\theta}$ should behave asymptotically like the sample maximum of $\theta^T y$ i.e. $\argmax_{y=X_1,\ldots,X_n} \theta^T y$. This is indeed the case, but unlike the one-dimensional case, due to angularity constraints the asymptotic distribution of the latter quantity is fairly complex.

\begin{theorem}\label{thm:emphd}
For a sequence $c_n, n \geq 1$ converging to infinity,
    \begin{enumerate}[label = (\alph*)]
\item If $\frac{M_{c_n\theta}}{n} \to 0$, then \begin{equation}
c_n(x_{\theta} - R_{n,c_n\theta}) \to Z,
\end{equation}
where $Z$ is as in Theorem~\ref{sltruemd}.
\item If $\frac{M_{c_n\theta}}{n} \to C \in (0,\infty)$, then \begin{equation}
c_n(x_{\theta} - R_{n,c_n\theta}) \to Z_{PRM}
\end{equation}
for a random variable $Z_{PRM} \neq Z$.
\item If $\frac{M_{c_n\theta}}{n} \to \infty$,
 then
 \begin{equation}
     \frac{(x_{\theta} - R_{n,c_n\theta})}{U^{-1}(\frac 1n)} \to V,
 \end{equation}
 where $V$ is a random vector depending on $\nu$ and the limiting PRM from the previous result.
 \end{enumerate}
\end{theorem}

We note that the final random variable $V$ has Weibull marginals and therefore extends the one-dimensional result. However, a combination of Lemma~\ref{KStoSLhd}  and Theorems~\ref{m2tmt2hd},~\ref{sltruemd} and ~\ref{thm:emphd} tells us that for multidimensional random variables with regularly varying tails, polynomially many samples are sufficient for tilting.

As an example, let $X = (X_1,X_2)$ where $X_1 \sim U[0,1]$ and $X_2 \sim U[0,1]^2$ are independent. Let $\theta = (1,1)$. In this case, $x_{\theta} = (1,1)$ is the maximizer of $\theta^T x$ for $x \in \supp(X) = [0,1]^2$. 

We claim that $X$ is multivariate regularly varying at $x_{\theta}$. To prove this, observe that if $A = [0,x_1] \times [0,x_2]$ for some $x_1,x_2>0$ then for any $t>0$, \begin{align*}
\mathbb P\left(\frac{x_{\theta} - X}{t} \in A\right) =& \mathbb P(X \in (1,1)-tA) \\
=& \mathbb P(X_1 \in [1-tx_1,1])\mathbb P(X_2 \in [1-tx_2,1]) \\
=& tx_1 (1-\sqrt{1-tx_2})
\end{align*}
since $X_1,X_2$ are independent. Now, we have $$
\frac 1{t^{2}}\mathbb P\left(\frac{x_{\theta} - X}{t} \in A\right) = x_1 \frac{(1-\sqrt{1-tx_2})}{t}.
$$
as $t \to 0$, we get $\frac{x_1x_2}{2}$ in the limit. Therefore, the limit measure $\nu$ is uniquely defined by its values on rectangles, which is given by $\nu([0,x_1] \times [0,x_2]) = \frac{x_1x_2}{2}$, and $U(t) = t^2$, which gives $\alpha=2$.

We can now apply the theorems. In particular, Theorem~\ref{m2tmt2hd} implies that as $c_n \to \infty$, 
$$
\frac{M_{c_n\theta}}{c_n^2} \to \frac{1}{2}.
$$
In particular, $M_{c_n \theta}$ grows at the same rate as $c_n^2$. By Theorem~\ref{sltruemd} we have $$
c_n((1,1)-X_{c_n\theta}) \to Z
$$
where $Z$ has density proportional to $e^{-\theta^T y} \nu(dy)$, which in particular implies that \begin{align*}
&\mathbb P(Z \in [0,x_1] \times [0,x_2]) \\ \propto &\int_{[0,x_1] \times [0,x_2]} e^{-\theta^T y} d\nu(y)\\
\propto & \frac 12\int_{0}^{x_1} e^{-y_1}dy_1 \int_{0}^{x_2} e^{-y_2}dy_2 \\
\propto &\mathbb P(\Gamma(1,1) \in [0,x_1])\mathbb P(\Gamma(1,1) \in [0,x_2])
\end{align*}

In particular, $Z$ is a product of two independent $\Gamma(1,1) = Exp(1)$ random variables in this case. By Theorem~\ref{thm:emphd}(a), if $\frac{M_{c_n\theta}}{n} \to 0$ then $n = \Omega(c_n^2)$, and $$
c_n((1,1)-R_{n,c_n\theta}) \to Z
$$
where $Z$ has independent $Exp(1)$ random variables as its components, as before. On the other hand, if $\frac{c_n^2}{n} \to C \in (0,\infty)$ then $\frac{M_{c_n\theta}}{n} \to \frac{C^2}{2}$, in which case Theorem~\ref{thm:emphd}(b) implies that $$
c_n((1,1)-R_{n,c_n\theta}) \to Z_{\frac{C}{2},PRM}
$$
for some random variable $Z_{\frac C2, PRM}$ which is not equal to $Z$ in distribution. Finally, if $n = o(c_n^2)$ then $\frac{M_{\theta}}{n} \to +\infty$ and $$
\sqrt{n}((1,1)-R_{n,c_n\theta})  \to V
$$
where $V$ is the scaling limit of $\argmax_{i} (X_{i})_{1}+(X_{i})_{2}$. It is possible to prove that $V$ has $-W_1$-marginals, however these marginals will be correlated.

The proofs of all results in this section will be included in Appendix~\ref{appendix:e}.

\section{TILTING IN THE UNBOUNDED SETTING}\label{unbdd}

In this section, we will provide some heuristic arguments as to why tilting unbounded random variables is likely to be hard. Our discussion focuses on the role of the quantity $M_{\theta}$ defined by \eqref{m2tmt2}, and why, unlike the previous two sections, it is either the wrong quantity to focus on, or indicative of requiring too many samples for a small twist. 

\subsection{Tilting Exponential Random Variables}

Consider the exponential random variable $X=Exp(\lambda)$ for some $\lambda>0$, and recall $X_{\theta}$ from \eqref{twisted}. It is easy to see that $X_{\theta}$ exists only till $\theta<\lambda$, and that $X_{\theta } \overset{d}{=} Exp(\lambda-\theta)$. 

If $\frac{\lambda}{2}<\theta < \lambda$, then note that $M_{\theta} = +\infty$. Therefore, unlike theorems in the previous sections, this quantity cannot dictate the rate at which the empirical estimator $R_{n ,\theta}$ given by \eqref{empirical1D} converges to the true distribution $X_{\theta}$, in the regime $\frac{\lambda}{2} < \theta_n< \lambda$. Contrast this with the discussion following Theorem~\ref{thm:genbd}. Thus, unbounded random variables whose moment generating function have finite radii require separate treatment for high tilts : such a question does not arise in the bounded case.

\subsection{Tilting Normal-like Random Variables}

In the bounded setting, the polynomial growth rate of $M_{\theta_n}$ as $\theta \to \infty$  (see Theorems~\ref{thm:m2tmt21D} and ~\ref{m2tmt2hd}) allowed us to tilt bounded distributions accurately with fewer samples (see Theorems~\ref{thm:slemp1D}(a) and ~\ref{thm:emphd}). We will demonstrate a similar trichotomy and prove that for unbounded light tailed random variables, the rate of growth of $M_{\theta}$ is exponential in $\theta$. This verifies our assertion that twisting unbounded random variables is much harder!

To simplify our assumptions as much as possible, let $X$ be a continuous random vector with full support on $\mathbb R^d$ and density $f$. Suppose that there exist $\alpha,K,L>0$ such that
\begin{equation}\label{eq:growth}
\lim_{x \to \infty} \frac{f(x)}{e^{-K\|x\|^{\alpha}}} \to L,
\end{equation}
in that for every $\epsilon>0$ there is an $R>0$ such that if $\|x\|>R$ then $ \frac{f(x)}{e^{-K\|x\|^{\alpha}}} \in (L-\epsilon,L+\epsilon)$. This includes standard Gaussian random vectors, with $\alpha = 2$ and some $K,L>0$, for instance.

As in the previous section, fix a direction $\theta \in \mathbb R^d$ with $\|\theta\| = 1$, and let $c>0$ be arbitrary. In this section, for ease of presentation we only deal with exponential tilting i.e. $g(x)=x$. 

Our results, like the previous sections, are divided into three parts : the asymptotics of $M_{c_n\theta}$, a scaling limit for $X_{c_n\theta}$, and scaling limits for $R_{n,c_n\theta}$ in three regimes depending upon $\frac{M_{c_n\theta}}{n}$. We remark, however, that the results are significantly harder to prove than the bounded case, even under our extremely natural and simple assumptions.

We introduce some notation that will be useful in stating and proving our results. Let $$\Phi_c(x) = c\theta^Tx - K\|x\|^{\alpha}.$$ The following properties of $\Phi_c$ are important.
\begin{lemma}\label{lem:phic}
For each $c>0$,  the function $\Phi_c(x)$\begin{enumerate}[label = (\alph*)]
\item attains its unique maximum at a point $$m_{c} =  \left(\frac{c}{\alpha K}\right)^{\frac 1{\alpha-1}}\theta,$$ and $$
\Phi_c(m_c) =  (\alpha-1)\alpha^{-\alpha/(\alpha-1)}c^{\alpha/(\alpha-1)}K^{-1/(\alpha-1)}.
$$
\item is thrice continuously differentiable in $\mathbb R^d \setminus \{0\}$, and \begin{gather*}
\nabla^2 \Phi_c(m_c) =  -K\alpha \|m_c\|^{\alpha-2} (I+(\alpha-2)\theta \theta^T).
\end{gather*}
(Note : in the second term, the product of vectors is the outer product).
\end{enumerate}
\end{lemma}

We are now ready to state our results.

\begin{theorem}\label{m2tmt2unbdd}
There exists constants $p_{\alpha,K,d},q_{\alpha,K,L,d}$ depending only upon $\alpha,K,L$ and $d$ such that as $c_n \to \infty$, $$
\lim_{c_n \to \infty} \frac{M_{c_n\theta}}{e^{p_{\alpha,K,d}c_n^{\alpha/(\alpha-1)}} c_n^{-d(\alpha-2)/(2\alpha-2)}} = q_{\alpha,K,L,d}.
$$
\end{theorem}

In particular, observe that $M_{c_n \theta}$ grows \emph{exponentially} with $c_n$, since $\alpha/(\alpha-1)>0$. The following scaling limit holds for $X_{c_n \theta}$.
\begin{theorem}\label{slunbdd}
As $c_n \to \infty$, $$(-\nabla^2 \Phi_{c_n}(m_{c_n}))^{-\frac 12}(X_{c_n \theta} - m_{c_n}) \overset{d}{\to} N(0,I).$$
\end{theorem}

The scaling limits for $R_{n,c_n\theta}$ can also be specified. Since the expressions are lengthy, we informally state this result below, and remark that they are motivated by the discussions preceding Theorems~\ref{thm:slemp1D} and ~\ref{thm:emphd}.

\begin{theorem}\label{empunbdd}
If $\frac{M_{c_n}}{n} \to 0$ then $X_{n},R_{n,c_n \theta}$ have the same scaling limit. If $\frac{M_{c_n}}{n} \to C \in (0,\infty)$ then $R_{n, c_n\theta}$ converges to a PRM-based functional. Finally, if $\frac{M_{c_n}}{n} \to \infty$ then $R_{n,c_n \theta}$ has the same scaling limit as the sample maximizer of $\theta^Ty$.
\end{theorem}

These results will be proved in Appendix~\ref{appendix:f}.

\section{Experiments}
\label{appendix:g}

In this section, we present a series of experiments. For simplicity, we get the tilting function to be the identity function, i.e., $g(x) = x$. For the sake of notation, $\theta$ will denote the tilting parameter and $n$ will denote the number of samples, unless otherwise specified.\footnote{The code is available \href{https://github.com/aistats20252404/codebase}{here}}

\subsection{Unbounded Random Variables}

In this section we look at unbounded random variables and how accurately they can be twisted. We start by considering an exponential random variable. 

\begin{figure}[H]
\centering
\fbox{\includegraphics[width=0.95\textwidth]{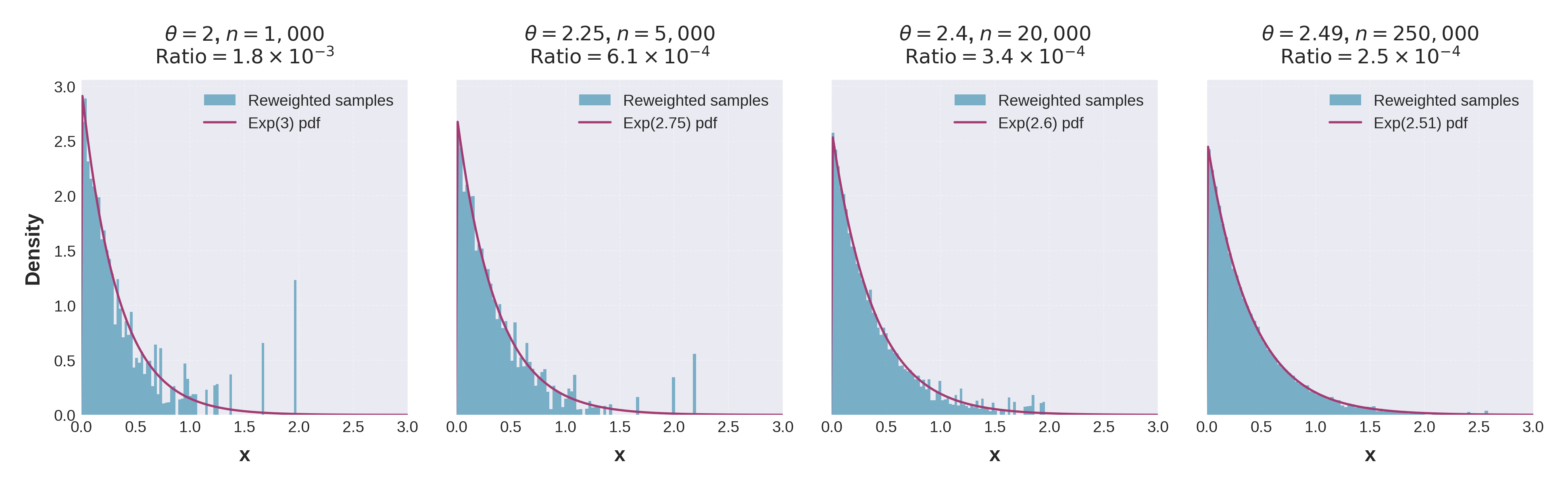}}
\caption{Exponential tilting of Exp(5) distribution with a sequence of $\left(\theta_i, n_i\right)$ s.t. $M_\theta / n \rightarrow 0$.}
\label{exp1}
\end{figure}

As Figure~\ref{exp1} shows, if $\frac{M_{\theta_n}}{n} \to 0$ then exponential tilting is asymptotically accurate, reflecting Theorem~\ref{thm:genbd}.

\subsection{Bounded Random Variables}

The next figures represent twists of bounded random variables in the Weibull regime.

\begin{figure}[H]
\centering
\fbox{\includegraphics[width=0.95\textwidth]{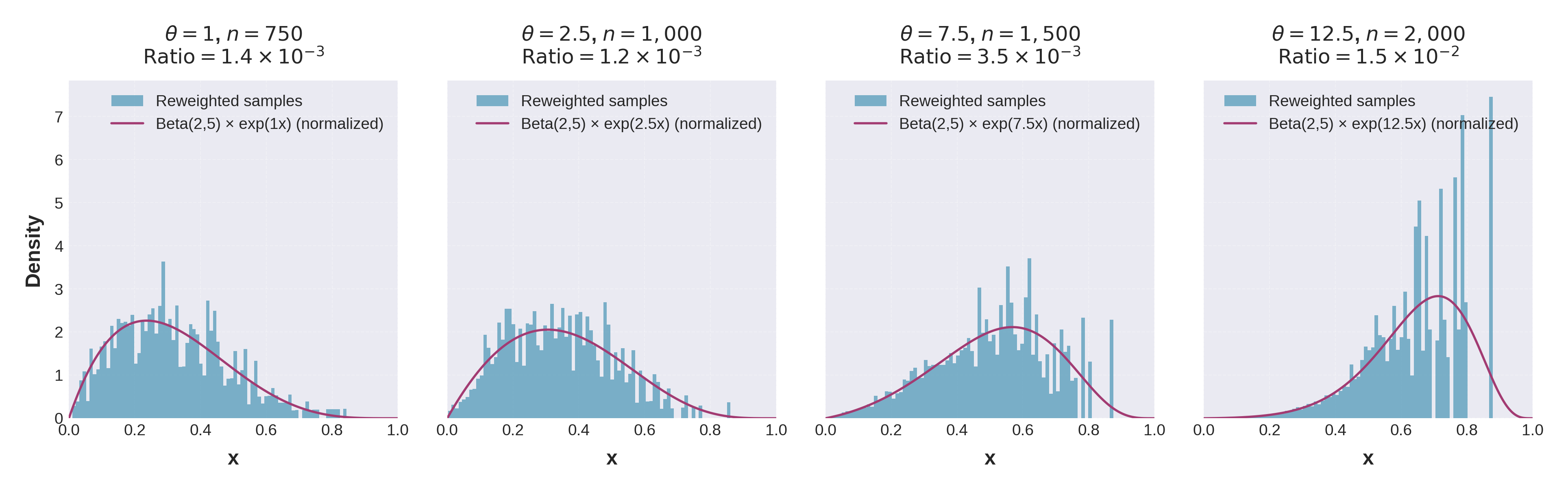}}
\caption{Exponential tilting of Beta(2, 5) distribution with a sequence of $\left(\theta_i, n_i\right)$ s.t. $M_\theta / n \not\rightarrow 0$.}
\label{exp2}
\end{figure}

As Figure~\ref{exp2} shows, if $\frac{M_{\theta_n}}{n} \to \infty$ then the sample maximizer receives higher empirical weight than the other samples, demonstrating Theorem~\ref{thm:slemp1D}(c).

\begin{figure}[H]
\centering
\fbox{\includegraphics[width=0.95\textwidth]{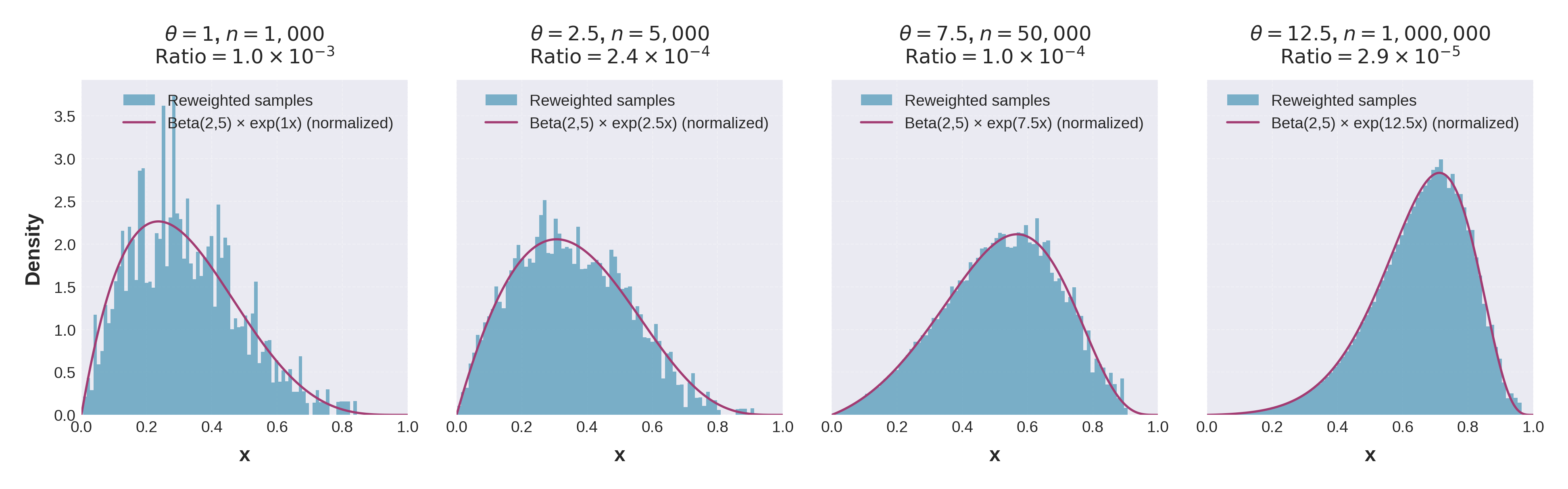}}
\caption{Exponential tilting of Beta(2, 5) distribution with a sequence of $\left(\theta_i, n_i\right)$ s.t. $M_\theta / n \rightarrow 0$.}
\label{exp3}
\end{figure}

\begin{figure}[H]
\centering
\fbox{\includegraphics[width=0.95\textwidth]{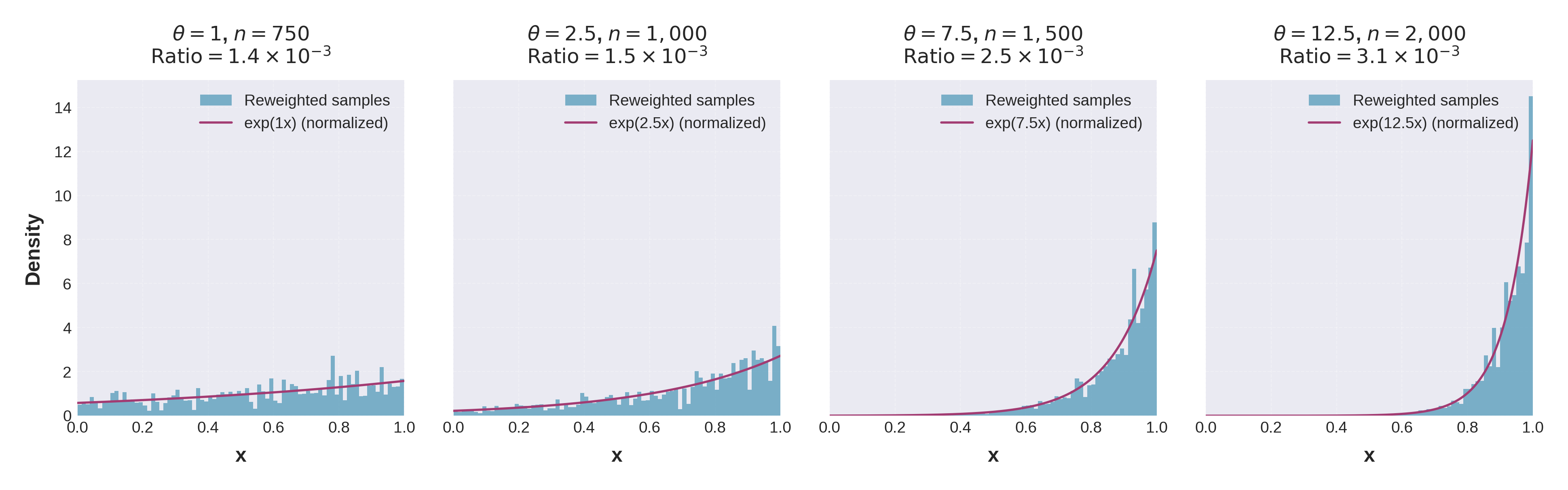}}
\caption{Exponential tilting of Uniform(0, 1) distribution with a sequence of $\left(\theta_i, n_i\right)$ s.t. $M_\theta / n \not\rightarrow 0$.}
\label{exp4}
\end{figure}

\begin{figure}[H]
\centering
\fbox{\includegraphics[width=0.95\textwidth]{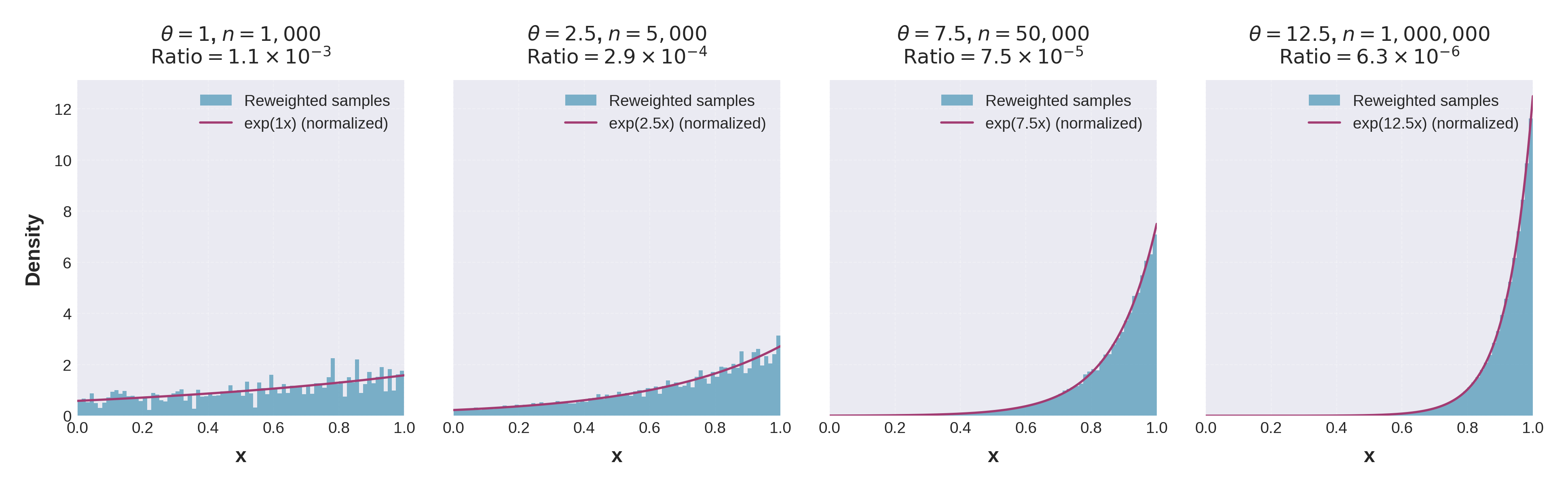}}
\caption{Exponential tilting of Uniform(0, 1) distribution with a sequence of $\left(\theta_i, n_i\right)$ s.t. $M_\theta / n \rightarrow 0$.}
\label{exp5}
\end{figure}

Conversely, if $\frac{M_{\theta_n}}{n} \to 0$, then twisting does take place accurately in the Weibull regime in Fig~\ref{exp3}, demonstrating Theorem~\ref{thm:slemp1D}(a). The next two figures, Fig~\ref{exp4} and Fig~\ref{exp5} demonstrate the same phenomena for exponentially tilting a uniform random variable.

\begin{figure}[H]
\centering
\fbox{\includegraphics[width=0.95\textwidth]{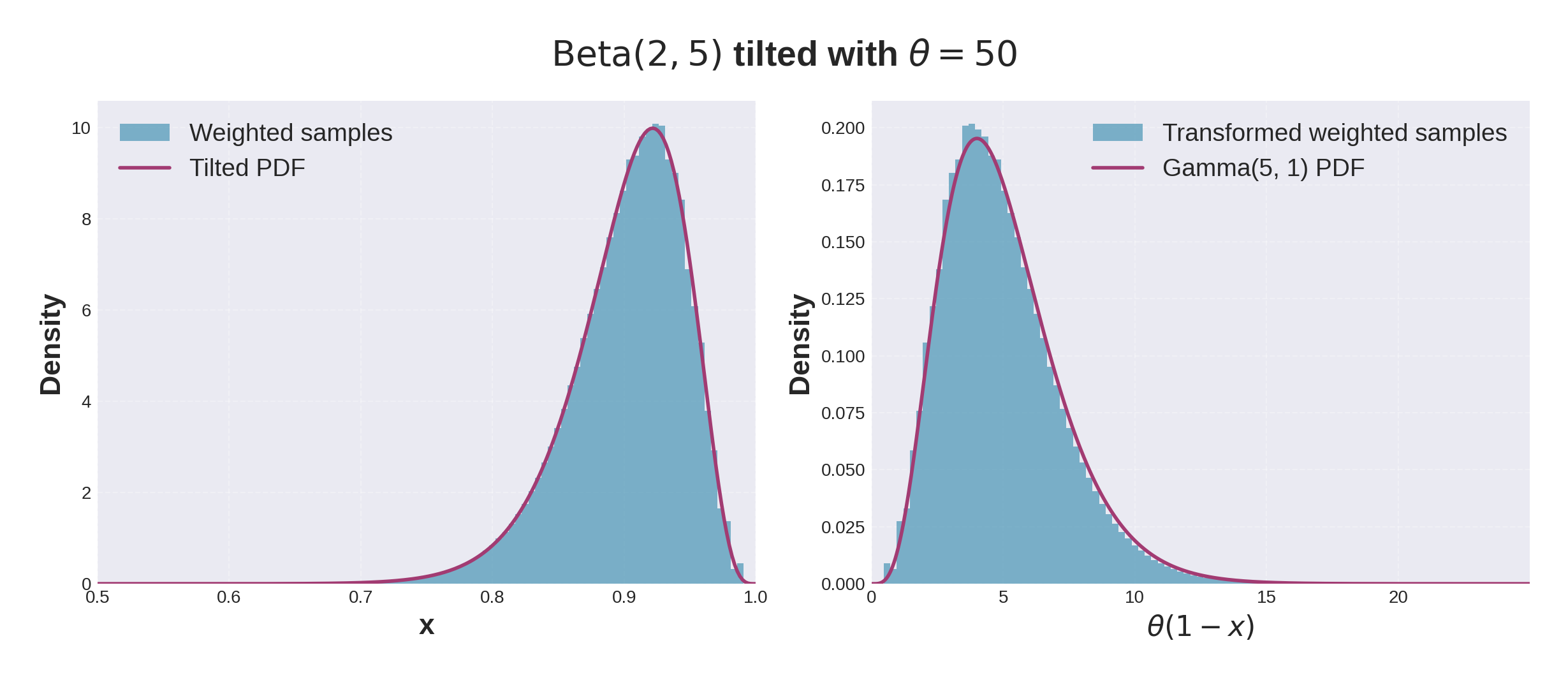}}
\caption{Exponential tilting of Beta(2, 5) distribution with $\theta = 50$. The PDF of the samples and true distribution is on the left. The transformed samples, given in theorem~\ref{thm:sltrue1D}, with the PDF of $\Gamma(5, 1)$, the corresponding scaling limit.}
\label{exp6}
\end{figure}

Finally, we demonstrate the Gamma scaling limit of the empirical random variable in Fig~\ref{exp6}, thereby adding credence to Theorem~\ref{thm:slemp1D}(a).

\section{CONCLUSION}\label{conc}

In this work, we discussed the asymptotic efficiency of a distribution estimator obtained by exponentially tilting the empirical distribution in different regimes. We provided scaling limits and characterized the asymptotic sample complexity needed to accurately twist distributions in the Weibull regime. We also showed the fundamental hardness of exponentially tilting unbounded distributions. As we observe in the introduction, our work builds upon the existing literature on SNIS estimators.

In the conclusion, we summarize our work's relationship with SNIS literature, highlighting our novel contributions to the area. In our future work we will also be testing these estimators in the context of generative neural networks.

However, we did not address unbounded random variables and vectors, particularly heavy tailed ones which lie in the Frechet regime. We believe this is a genuinely hard problem, since such distributions can only be subject to polynomial tilting.

\bibliography{arXivFINAL29122025}
\bibliographystyle{alpha}

\clearpage

\thispagestyle{empty}

\begin{appendices}

\section{Proofs for Section 2}\label{appendix:a}

In this section, we will prove Theorem~\ref{thm:bd} with the assistance of techniques from empirical process theory and inequalities for the suprema of Gaussian processes. 

We begin by recalling notation from Section~\ref{fixedtheta}. Let $X$ be a continuous random variable, $n \geq 1$ be a fixed positive integer, $\theta \in \mathbb{R}$ a fixed positive real number and $g : \R \to \R$ a strictly increasing, continuous function, which in particular ensures that $e^{\theta g(x)}$ is strictly increasing in $x$.

Throughout this section, we set the following notation. Let $X_1,\ldots,X_n$ be independent and identically distributed as $X$. Let $Y= e^{\theta g(X)}$ and $Y_i = e^{\theta g(X_i)}$.

In the first part of this section, we will prove the first part of Theorem~\ref{thm:bd} i.e. the existence of the limiting Gaussian process $\mathcal G_{\theta}$ with some specified covariance. In the second part, we will prove the bounds on the limiting random variable.

\subsection{Proof of Gaussian limit $\mathcal G_{\theta}$ in Theorem~\ref{thm:bd}}

Recall that the true and empirical distributions $X_{\theta}$ and $R_{n,\theta}$ are given by \eqref{twisted1D} and \eqref{empirical1D} respectively, and the $\KS$ distance is given by \eqref{KS}. Let $F_{n,\theta}(x)$ be the CDF of $R_{n,\theta}$, and for convenience define $G_n(x), x \in (-\infty,+\infty]$ by
\begin{gather}
G_n(x) = \frac 1n\sum_{i=1}^n Y_i\1_{X_i \leq x}. \label{gn}
\end{gather}
Note that $G_n(+\infty) = \frac 1n \sum_{i=1}^n Y_i$. By the law of large numbers, the random vector $(G_n(x), G_n(+\infty))$ converges a.s. to $(\E[Y\1_{X \leq x}], \E[Y])$ for each fixed $x\in \R$. By the continuous mapping theorem, we have
$$
F_{n,\theta}(x) = \frac{G_n(x)}{G_n(+\infty)} \to \frac{\E[Y\1_{X \leq x}]}{\E[Y]} = F_{\theta}(x)
$$
a.s. for each $x \in \R$. However, the quantity we are studying is
\begin{align}
\KS(F_{n,\theta}, F_{\theta})  =& \sup_{x} \left|F_{n,\theta}(x) - F_{\theta}(x)\right| \nonumber \\ =& \sup_{x} \left|\frac{G_n(x)}{G_n(+\infty)} - \frac{\E[Y\1_{Y \leq e^{\theta g(x)}}]}{\E[Y]}\right|.\label{eq:Y}
\end{align}
Therefore, we require uniform control in the convergence rates of each of these random variables; in particular, we must look at the entire stochastic process $(G_{n}(x))_{x \in (-\infty,+\infty]}$ and show that it converges somewhere in some uniform fashion. However, we have the machinery to demonstrate not just uniform convergence, but the Gaussian limit as well. To set up the appropriate notation, let
\begin{equation}\label{F}
\mathcal F = \left\{f_x : f_x(y) = y \1_{y \leq e^{\theta g(x)}}, y \in \R_+, x \in (-\infty,+\infty)\right\}\cup \{f_{\infty}\},
\end{equation}
where $f_{\infty}(y) = y$. By \eqref{gn}, $$
G_n(x) = \frac 1n \sum_{i=1}^n f_x(Y_i) \quad,\quad \E[G_n(x)] = \E[f_x(Y_i)] = \E[Y 1_{X \leq x}].
$$

In order to apply results from empirical process theory, we first quickly verify the equation above \cite[equation (2.1.1)]{vandervaart1996weak}. Observe that for any $y \in \R$, $$
\sup_{f_x \in \mathcal{F}} |f_x(y)-\E[f_x(Y)]| \leq \sup_{f_x \in \mathcal{F}} |f_x(y)| + \sup_{f_x \in \mathcal{F}} \E[f_x(Y)] \leq y + \E[Y].
$$
The above inequality establishes that if for a class of functions $\mathcal F'$ on a normed space, we define \begin{equation}\label{linfty}
l^{\infty}(\mathcal{F'}) = \left\{h : \mathcal{F'} \to \mathbb R, \sup_{f \in \mathcal{F'}} |h(f)| < \infty\right\},\quad \|h\|_{\mathcal{F'}} = \sup_{f \in \mathcal{F'}} |h(f)|,
\end{equation} 
then $G_n$ can be seen as a random element of $l^{\infty}(\mathcal F)$ by letting $f_x$ be the input instead of $x$ i.e. $G_n(f_x) = \frac 1n \sum_{i=1}^n f_x(Y_i)$. However, we abuse notation and let $G_n(x)$ denote the random element itself.

The following lemma shows that the fluctuations $G_n- \E[G_n]$ converge weakly to a Gaussian process on $l^{\infty}(\mathcal F)$. For the definition of weak convergence (which we denote by $\Rightarrow$), we refer the reader to \cite[Section 1.3]{vandervaart1996weak}.

\begin{lemma}\label{lem:conv}
We have $$
\sqrt{n}(G_n - \E[G_n]) \Rightarrow \mathcal G'_{\theta}
$$
for a centered Gaussian process $\mathcal G'_{\theta}$ on $(-\infty,+\infty]$ with covariance functional given by $$
\Cov(\mathcal{G}'_{\theta}(x), \mathcal G'_{\theta}(y)) = \E[f_x(Y)f_y(Y)] - \E[f_x(Y)]\E[f_y(Y)] = \Cov(Y1_{X \leq x}, Y \1_{X \leq y}).
$$
Here, the weak convergence takes place on the space $l^{\infty}(\mathcal F)$.
\end{lemma}

Before we begin the proof of the lemma, some definitions are in order. Given two functions $u,l : (-\infty,+\infty] \to \R$, define the "bracket" class of functions $[u,l]$ by 
\begin{equation}\label{bracket}
[u,l] = \{f : (-\infty,+\infty] \to \R : u \leq f \leq l\}.
\end{equation}
Given $\epsilon>0$, an $\epsilon$-bracket is a class of functions $[u,l]$ such that $\sqrt{\E[(l(Y)-u(Y))^2]}< \epsilon$. We call a class $\mathcal{S}$ of $\epsilon$-brackets as an $\epsilon$-bracket cover of $\mathcal{F}$, if the union of all $\epsilon$-brackets in $\mathcal{S}$ contains $\mathcal{F}$. The smallest size of such a cover $\mathcal{S}$ will be denoted by $N_{[]}(\epsilon, \mathcal{F}, L^2(Y))$, the \emph{$\epsilon$-bracket covering number} of $\mathcal{F}$.

We now begin the proof of Lemma~\ref{lem:conv}.

\begin{proof}[Proof of Lemma~\ref{lem:conv}] 

By the discussion at the start of \cite[Section 2.5.2]{vandervaart1996weak}, it suffices to prove that \begin{equation}\label{entropy}
\int_{0}^{\infty} \sqrt{\log N_{[]}(\epsilon, \mathcal{F}, L^2(Y))} d \epsilon < \infty.
\end{equation}

For this, we need an upper bound on $N_{[]}(\epsilon, \mathcal{F}, L^2(Y))$. This is obtained by finding at least one $\epsilon$-bracket cover for each $\epsilon >0$, and finding its size.

Instead, we propose a sequence of classes containing brackets, which are naturally inspired by the structure of $\mathcal{F}$. For each $n \geq 1$, let $x_1<x_2<\ldots<x_n$ be real numbers chosen such that if $x_{n+1} = +\infty$ and $x_0 = -\infty$ then 
\begin{equation}\label{brac}
\E[(f_{x_i}(Y)- f_{x_{i-1}}(Y))^2] = \E[Y^21_{X \in (x_{i-1},x_i)}] = \frac {\E[Y^2]}{n+1}.
\end{equation}
for all $i=1,2,\ldots,n+1$. Note that such a choice requires the continuity of $Y$ and $X$ (and hence $X$ alone).

Let \begin{equation}\mathcal{S}_n = \{[f_{x_{i}}, f_{x_{i+1}}] : i=1,\ldots,n\}\cup\{[0, f_{x_1}]\}.\label{sn}\end{equation} Then, $\mathcal{S}_n$ covers $\mathcal{F}$ by \eqref{bracket} and \eqref{F}. Furthermore, by \eqref{brac}, $\mathcal{S}_n$ is a $\sqrt{\frac{\E[Y^2]}{n+1}}$-bracket cover of $\mathcal{F}$. Finally, $|\mathcal{S}_n| = n+1$ by \eqref{sn}.

Thus, let $\epsilon>0$ be arbitrary. Observe that if $\epsilon>\sqrt{\frac{\E[Y^2]}{n+1}}$, then $\mathcal{S}_n$ is also an $\epsilon$-bracket cover of $\mathcal{F}$. By the definition of $N_{[]}(\epsilon, \mathcal{F}, L^2(Y))$, \begin{equation}\label{nbd}
\sqrt{\frac{\E[Y^2]}{n+1}}<\epsilon < \sqrt{\frac{\E[Y^2]}{n}} \implies \sqrt{\log N_{[]}(\epsilon, \mathcal{F}, L^2(Y))} \leq \sqrt{\log(n+1)}.
\end{equation}
Finally, if $\epsilon> \sqrt{\E[Y^2]}$, then $\mathcal{S} =\{[0,f_{+\infty}]\}$ is an $\epsilon$-bracket cover of $\mathcal{F}$. Thus, $\log N_{[]}(\epsilon, \mathcal{F}, L^2(Y)) = 0$. We can now upper bound the left hand side of \eqref{entropy} using \eqref{nbd} and this observation.
\begin{align*}
&\int_{0}^{\infty} \sqrt{\log N_{[]}(\epsilon, \mathcal{F}, L^2(Y))} d \epsilon \\
= & \int_{0}^{\sqrt{\E[Y^2]}} \sqrt{\log N_{[]}(\epsilon, \mathcal{F}, L^2(Y))} d \epsilon \\
\leq & \sum_{n=1}^{\infty} \left(\sqrt{\frac{\E[Y^2]}{n}}  - \sqrt{\frac{\E[Y^2]}{n+1}}\right) \sqrt{\log(n+1)} \\
\leq & \sqrt{\E[Y^2]} \sum_{n=1}^{\infty} \left(\frac{1}{\sqrt{n}} - \frac{1}{\sqrt{n+1}}\right)\sqrt{\log(n+1)} < C \sqrt{\E[Y^2]} < \infty,
\end{align*}
where $$C = \sum_{n=1}^{\infty} \left(\frac{1}{\sqrt{n}} - \frac{1}{\sqrt{n+1}}\right)\sqrt{\log(n+1)} < \infty.$$

Thus \eqref{entropy} has been proved, and consequently we have 
\begin{equation}\label{conv}
\sqrt{n}\left(\frac 1n\sum_{i=1}^n f_x(Y_i) -\E[f_x(Y_i)]\right) \Rightarrow \mathcal G'_{\theta}
\end{equation}
for some Gaussian process $\mathcal G'_{\theta}$ on $(-\infty,+\infty]$. The covariance functional in the lemma follows from \cite[equation (2.1.2)]{vandervaart1996weak}.
\end{proof}

An application of the $\delta$-method should help us interpolate between the process $(G_n(x))_{x \in (-\infty,+\infty]}$ and $(F_{n,\theta}(x))_{x \in (-\infty,+\infty)}$.

\begin{proposition}\label{lem:gauss}
We have \begin{equation*}
\sqrt{n}(F_{n,\theta}(x) - F_{\theta}(x)) \Rightarrow \mathcal G_{\theta}
\end{equation*}
where $\mathcal G_{\theta}$ is a centered Gaussian process on $\mathbb R$ with covariance functional $$
Cov(\mathcal{G}_{\theta}(x_1), \mathcal{G}_{\theta}(x_2)) = \frac 1{y_3^2} \left(v_{12} - v_{23}\frac{y_1}{y_3} - v_{13}\frac{y_2}{y_3} + v_{33}\frac{y_1y_2}{y_3^2}\right) \quad x_1,x_2 \in \R.
$$
Here, $x_3 = +\infty$, $v_{ij} = \Cov(Y\1_{X\leq x_i} , Y\1_{X \leq x_j})$ and $y_i = \mathbb E[Y \1_{X \leq x_i}]$.

(Note : This weak convergence takes place on the space
$l^{\infty}(\mathcal F \setminus \{f_{\infty}\})$, which we recall from \eqref{linfty}.)

Furthermore, recalling the $\KS$ distance from $\eqref{KS}$, $$
\sqrt{n}\KS(F_{n,\theta}, F_{\theta}) \Rightarrow \mathcal \sup_{x} |\mathcal G_{\theta}(x)|
$$
in distribution (in the usual sense i.e. as real valued random variables).
\end{proposition}

\begin{proof}
Define a map $\phi: l^{\infty}(\mathcal{F}) \to l^{\infty}(\mathcal F \setminus \{f_{\infty}\})$ by $$
[\phi(h)](f_x) = \frac{h(f_x)}{h(f_{\infty})}.
$$
Observe that $\phi(G_n) = F_{n,\theta}$ by \eqref{gn} and \eqref{empirical1D}, and similarly $\phi(\E[G_n]) = F_{\theta}$, where $\E[G_n]$ is just a constant functional on $\mathcal{F}$. We wish to apply \cite[Theorem 3.9.4]{vandervaart1996weak}, for which it is sufficient to check that $\phi$ is Hadamard differentiable (see page 372, \cite[Section 3.9.1]{vandervaart1996weak} for the definition of Hadamard differentiability) on $h = \E[G_n]$,  where we note that $\E[G_n(+\infty)] \neq 0$.

This is a rather standard exercise which we include for the sake of completeness. Let $j_n$ be a sequence converging to $j$ in $l^{\infty}(\mathcal F)$, with $t_n \to 0$. For $f \in \mathcal{F}$ we must assess the limit of the quantity \begin{align}
\left[\frac{\phi(h+t_nj_n) - \phi(h)}{t_n}\right](f) =& \frac 1{t_n} \left[\frac{(h+t_nj_n)(f)}{(h+t_nj_n)(f_{\infty})} - \frac{h(f)}{h(f_{\infty})}\right] \nonumber\\
=& \frac 1{t_n}\left[(h+t_nj_n)(f)\left(\frac{1}{(h+t_nj_n)(f_{\infty})} - \frac {1}{h(f_{\infty})}\right)+\frac{t_nj_n(f)}{h(f_{\infty})}\right] \nonumber\\
& = -\frac{(h+t_nj_n)(f)j_n(f_{\infty})}{(h+t_nj_n)(f_{\infty})h(f_{\infty})} + \frac{j_n(f)}{h(f_{\infty})} \nonumber \\
& \to -\frac{h(f)j(f_{\infty})}{h(f_{\infty})^2} + \frac{j(f)}{h(f_{\infty})} = [[\phi'(h)](j)](f).\label{diff}
\end{align}
Thus, $\phi$ is Hadamard differentiable at $h$ with the above derivative. By \cite[Theorem 3.9.4]{vandervaart1996weak}, 
\begin{equation}\label{last}
\sqrt{n}(F_{n,\theta} - F_{\theta}) \Rightarrow [\phi'[h]](\mathcal G'_{\theta}) = \mathcal G_{\theta}.
\end{equation}

 Next, we unravel the notation on the right hand side to ascertain the true limit. Indeed, since $\phi'[h]$ is a linear map, $\mathcal G_{\theta}$ is also a centered Gaussian process. We have, by definition, that $h(f_x) = \E[Y 1_{X \leq x}]$ and $h(f_{\infty}) = \E[Y]$. Therefore, by \eqref{diff}, $$
\mathcal G_{\theta}(f_x) = \frac{\mathcal G'_{\theta}(f_x)}{\E[Y]}-\frac{\E[Y1_{X \leq x}]\mathcal G'_{\theta}(f_{\infty})}{\E[Y]^2}.
$$

All we must do now is find the covariance using the above expression and Lemma~\ref{lem:conv}. For any $x, y \in \mathbb R$, if we treat $\mathcal G_{\theta}(f_x) \equiv \mathcal G_{\theta}(x)$ as a stochastic process on $\R$ instead of a random element of $\mathcal{F} \setminus \{f_{\infty}\}$, we have
\begin{align*}
    &\Cov(\mathcal G_{\theta}(x), \mathcal G_{\theta}(y)) \\ =& \frac{\E[Y1_{X \leq x}]\E[Y\1_{X \leq y}]}{\E[Y]^4}\Var\left(\mathcal G'_{\theta}(\infty)\right) - \Cov (\mathcal G'_{\theta}(x), \mathcal G'_{\theta}(\infty))\frac{\E[Y \1_{X \leq y}]}{\E[Y]^3} \\ -&  \Cov (\mathcal G'_{\theta}(y), \mathcal G'_{\theta}(\infty))\frac{\E[Y \1_{X \leq x}]}{\E[Y]^3} + \frac{\Cov(\mathcal G'_{\theta}(x),\mathcal G'_{\theta}(y))}{\E[Y]^2} \\
    =& \frac{\E[Y1_{X \leq x}]\E[Y\1_{X \leq y}]\Var(Y)}{\E[Y]^4}-  \frac{\E[Y \1_{X \leq y}]\Cov(Y, Y\1_{X \leq x})}{\E[Y]^3} \\
    -& \frac{\E[Y \1_{X \leq x}]\Cov(Y, Y\1_{X \leq y})}{\E[Y]^3} +  \frac{\Cov(Y\1_{Y \leq x}, Y \1_{Y \leq y})}{\E[Y]^2} \\
\end{align*}
which, after observation, is of the form given in the proposition, as desired.

The convergence of $\sqrt{n} \KS(F_{n,\theta} , F_{\theta})$ follows, since $$\KS(F_{n,\theta} , F_{\theta}) = \sup_{x}|F_{n,\theta} - F_{\theta}| = \|F_{n,\theta}-F_{\theta}\|_{l^{\infty}(\mathcal{F})}$$ by \eqref{linfty}. The norm is a continuous functional, therefore by the continuous mapping theorem \cite[Theorem 1.3.6]{vandervaart1996weak}, it follows that $$
\sqrt{n} \KS(F_{n,\theta} , F_{\theta}) \overset{d}{\to} \sup_x |\mathcal G_{\theta}(x)|,
$$
as desired.
\end{proof}

Note that Proposition~\ref{lem:gauss} has completed a part of the proof of Theorem~\ref{thm:bd}. The bounds on the expectation of the limiting random variable will be shown in the next subsection.

\subsection{Proof of bounds in Theorem~\ref{thm:bd}}\label{sec:prop}

We shall first state some notions required to prove the required bounds. Indeed, both the upper and lower bound can be given in terms of covering numbers; however, the notion of distance and cover will be simpler than that considered in the previous subsection. 

Let $\mathcal{G}_{\theta}$ be as in Proposition~\ref{lem:gauss} and define a distance on $\mathbb R$ by 
\begin{equation}\label{eq:dist}
d(s,t) = \sqrt{\mathbb E[|\mg_{\theta}(s) - \mg_{\theta}(t)|^2]}.
\end{equation}
For any $x \in \mathbb R, r>0$ let $B_d(x,r) = \{y \in \mathbb R : d(x,y) \leq r\}$. We can now define, for any $\epsilon>0$, the covering number 
\begin{equation}\label{eq:pn}
N(\epsilon) = \inf\{n : \exists x_1,x_2,\ldots,x_n \in \mathbb R, \mathbb R = \cup_{i=1}^n B_d(x_i,\epsilon)\}.
\end{equation}
That is, $N(\epsilon)$ is the smallest size of a set satisfying the property that every real number is at most $\epsilon$ away (in distance $d$) from some point in this set. 

We will also define some constants and set some notation that will make analysis easier in the upcoming sections. Define
\begin{equation}\label{eq:K}
K = \sup_{s,t \in \mathbb R} d(s,t).
\end{equation} 
$K$ is the diameter of the metric space $(\mathbb R,d)$. We need estimates on $K$, which are difficult to obtain without having a simpler formula for $d$. The following lemma obtains such a formula.

\begin{lemma}\label{lem:rew}
For any $x_1<x_2$ (and symmetrically for $x_2\leq x_1$), if $A = \{X \in [x_1,x_2]\}$, then $$
d^2(x_1,x_2) = \frac 1{\mathbb E[Y]^4} \left(\mathbb E [Y^2\1_A]\mathbb E [Y\1_{A^c}]^2 + \mathbb E [Y^2\1_{A^c}] \mathbb E [Y\1_A]^2\right).
$$
\end{lemma}

Assuming the above lemma, we can obtain bounds on $K$. Recall $M_{\theta}$ from \eqref{m2tmt2}.

\begin{proposition}\label{thm:K}
We have $$
\sqrt{M_{\theta}} \geq K \geq \frac 12 \sqrt{M_{\theta}}.
$$
\end{proposition}

Finally, the following lemma will be used to control $d(x,y)$ when $x,y$ are (in the usual metric) close to each other.

\begin{lemma}\label{lem:est}
Let $A = \{X \in [s,t]\}$. If $\frac{\mathbb E[Y\1_A]}{\mathbb E[Y]}  = L\leq \frac 12$, then
$$
d^2(s,t) \leq L^2\frac{\mathbb E[Y^2]}{\mathbb E[Y]^2} = L^2 M_{\theta}.
$$
\end{lemma}

The three results above will all be proved in the order they appear in the last subsection. We shall assume these and now proceed to the proof of the lower bound.

\subsubsection{Lower bound}

The following theorem(\cite[Theorem A.2.5]{vandervaart1996weak}) will be used to lower bound $\mathbb E \sup_x |\mg_{\theta}(x)|$.

\begin{theorem}[Sudakov-Fernique]\label{thm:sf}
For any $\epsilon>0$,$$
\mathbb E\left[\sup_x |\mg_{\theta}(x)|\right] \geq \frac 13 \epsilon \sqrt{\log N(\epsilon)}.
$$
\end{theorem}

The lower bound in Theorem~\ref{thm:bd} is now an easy consequence of the above theorem.

\begin{proof}[Proof of lower bound in Theorem~\ref{thm:bd}]
Let $\epsilon = K/3$. We claim that $N(\epsilon) \geq 2$. Indeed, if $N(\epsilon) = 1$, then there is a point $y \in \mathbb R$ such that $d(y,x) < \frac K3$ for all $x \in \mathbb R$. However, by the triangle inequality, for any $x,z \in \R$, 
$$d(x,z) \leq d(x,y)+d(y,z) \leq \frac{2K}{3},$$ contradicting the fact that $K$ is the diameter of the space. By Theorem~\ref{thm:sf} and Proposition~\ref{thm:K}, it follows that 
$$
\mathbb E\left[\sup_x |\mg_{\theta}(x)|\right]  \geq \frac{K}{9} \sqrt{\log 2} \geq C\sqrt{M_{\theta}}.
$$ 
for some constant $C = \frac{\sqrt{\log 2}}{18}$ independent of $\mg_{\theta}$.
\end{proof}

\subsubsection{Upper bound}

We shall now derive an upper bound on $\mathbb E[\sup_x |\mg(x)|]$. For this, we use the following inequality, which is \cite[Corollary 2.2.8]{vandervaart1996weak}, augmented with the relation between covering and packing numbers under \cite[Definition 2.2.3, page 98]{vandervaart1996weak}.

\begin{theorem}\label{thm:sud}
Recall the covering number $N(\epsilon)$ from \eqref{eq:pn}. There is a universal constant $\kappa$ such that for any $t_0 \in \R$, we have $$
\mathbb E\left[\sup_x |\mg_{\theta}(x)|\right] \leq \E |\mg_{\theta}(t_0)| + \kappa \int_0^{\infty} \sqrt{\log N\left(\frac 12\epsilon\right)} d \epsilon.
$$
\end{theorem}

We will now prove the upper bound in Theorem~\ref{thm:bd}. The proof is very similar to the proof of Lemma~\ref{lem:conv} but they differ mildly in the manner in which the cover is constructed.

\begin{proof}[Proof of upper bound in Theorem~\ref{thm:bd}]
We will look to bound the right hand side of Theorem~\ref{thm:sud}. In order to derive an upper bound for $N(\epsilon)$ for any $\epsilon>0$, we must find at least one set of points $S \subset \mathbb R$ such that for every $x \in \mathbb R$, $B_d(x,\epsilon) \cap S \neq \emptyset$. We can then say that $|S| \geq N(\epsilon)$. 

Suppose $n \in \mathbb N$. By continuity of $Y$, we can find points $x_1,x_2,\ldots,x_n$ such that if $x_0=-\infty, x_{n+1} = +\infty$ and $A_n =\{X\in [x_n,x_{n+1}]\}$, then \begin{equation}\label{eq:ineq}
\frac{\mathbb E[Y\1_{A_n}]}{\mathbb E[Y]} = \frac {1}{n+1}.
\end{equation}

We claim that if $\epsilon > \frac{\sqrt{M_{\theta}}}{n+1}$ then $\cup_{i=1}^n B_d(x_i,\epsilon) = \mathbb R$. Indeed, suppose that $y \in \mathbb R$. Then, $y \in A_i$ for some $i$. Let $x_k$ be any one of the endpoints of $A_i$ (or the only real endpoint, if the other endpoint is infinite). The interval $[y,x_k]$ (or $[x_k,y]$ depending upon the order between 
the two) is contained in $A_i$. Therefore, if $L = \frac{\mathbb E[Y\1_{X \in [y,x_k]}]}{\mathbb E[Y]}$, then $L \leq \frac{1}{n+1}$. By Lemma~\ref{lem:est} and \eqref{eq:ineq},
$$
d(y,x_k)\leq L \sqrt{M_{\theta}} \leq \frac{\sqrt{M_{\theta}}}{n+1}< \epsilon.
$$
This completes the proof of the claim since $y \in \mathbb R$ was arbitrary.

Having obtained this, we are now ready to estimate the integral on the right hand side in Theorem~\ref{thm:sud}. Observe that $N\left(\frac 12 \epsilon\right) \geq 2$ if $K \geq \frac{\epsilon}{2} > \frac{\sqrt{M_{\theta}}}{3}$, where $K$ is defined in \eqref{eq:K}. On the other hand, we have shown that for any $n \geq 3 \in \mathbb N$ that
$$
\frac{\sqrt{M_{\theta}}}{(n+1)} < \frac{\epsilon}{2} \leq  \frac{\sqrt{M_{\theta}}}{n} \implies N\left(\frac{\epsilon}{2}\right) \leq n.
$$
Therefore, \begin{equation}\label{eq:sude}
\int_0^{K} \sqrt{\log N\left(\frac{\epsilon}{2}\right)} d \epsilon \leq \left(2K - \frac{2\sqrt{M_{\theta}}}{3}\right) \sqrt{\log 2} + \sum_{n=3}^{\infty} \sqrt{2\log(n)}\left(\frac{\sqrt{M_{\theta}}}{n} - \frac{\sqrt{M_{\theta}}}{n+1}\right)
\end{equation}
It only remains to show that the right hand side is bounded by a constant multiple of $\sqrt{M_{\theta}}$, completing the proof.  However, $K \leq \sqrt{M_{\theta}}$ by Proposition~\ref{thm:K}. Therefore, $$
\left(K - \frac{\sqrt{M_{\theta}}}{3}\right) \sqrt{\log 2} \leq \left(\frac{2}{3} \sqrt{\log 2}\right)\sqrt{M_{\theta}} = C \sqrt{M_{\theta}}.
$$
for some $C$ independent of $M$. Applying this, the right hand side of \eqref{eq:sude} is bounded by $$
\sqrt{M_{\theta}}\left(C + \sum_{n=3}^{\infty} \sqrt{\log(n)}\left(\frac{1}{n} - \frac{1}{n+1}\right)\right) \leq C' \sqrt{M_{\theta}}
$$
for some constant $C'>0$, since the series within the bracket converges.

Combining Theorem~\ref{thm:sud}, \eqref{eq:sude} and the above estimate, we obtain for any $t_0 \in \R$ that
\begin{equation}\label{ineq}
\mathbb E\left[\sup_{x} |\mathcal{G}_{\theta}(x)|\right] \leq \E[|\mathcal{G}_{\theta}(t_0)|]+ C_4 \sqrt{M_{\theta}}.
\end{equation}

Now, $\mathcal G_{\theta}$ is a mean zero Gaussian random variable with variance specified by Proposition~\ref{lem:gauss}. As $t_0 \to -\infty$, observe that $Y1_{X \leq t_0} \to 0$ in $L^2$ by the dominated convergence theorem. Therefore, $\Var(\mathcal G_{\theta}(t_0)) \to 0$ as $t_0 \to -\infty$. However, noting that $$\E[|\mathcal G_{\theta}(t_0)|] = \sqrt{\Var(\mathcal G_{\theta}(t_0))}\sqrt{2/\pi}$$
it follows that this particular term can be made as small as necessary. In particular, let $t_0$ be chosen so that $\E[|\mathcal G_{\theta}(t_0)|] \leq \frac 12$. Then, by \eqref{ineq}, 
\begin{align*}
\mathbb E\left[\sup_{x} |\mg_{\theta}(x)|\right] \leq \frac 12 + C_4 \sqrt{M_{\theta}}
\leq\left(\frac 12+C_4\right) \sqrt{M_{\theta}}
\end{align*}
since $M_{\theta}\geq 1$. This completes the proof.
\end{proof}

\subsubsection{Proof of Lemma~\ref{lem:rew}, Proposition~\ref{thm:K} and Lemma~\ref{lem:est}}

Having completed the proof of the bounds in  Theorem~\ref{thm:bd}, we now present the proofs of the lemmas and theorems which we used to prove them, starting with Lemma~\ref{lem:rew}.

\begin{proof}[Proof of Lemma~\ref{lem:rew}]
Recall the notation used to define the covariance in Proposition~\ref{lem:gauss}. We will use this notation in our proof. That is, let $x_1,x_2 \in \R$, and $x_3 = +\infty$. We reuse the notation $v_{ij} = \Cov(Y\1_{X \leq x_i} ,Y \1_{X \leq x_j})$ and $y_i = \mathbb E[Y \1_{X \leq x_i}]$. With this, if $Z_i = \mathcal{G}_{\theta}(x_i)$, then \begin{gather*}
\Cov(Z_1,Z_1) = \frac 1{y_3^2} \left(v_{11} - 2v_{13}\frac{y_1}{y_3} + v_{33}\frac{y_1^2}{y_3^2}\right). \\
\Cov(Z_2,Z_2) = \frac 1{y_3^2} \left(v_{22} - 2v_{23}\frac{y_2}{y_3} + v_{33}\frac{y_2^2}{y_3^2}\right) \\
-2\Cov(Z_1,Z_2) = \frac {-2}{y_3^2} \left(v_{12} - v_{23}\frac{y_1}{y_3} - v_{13}\frac{y_2}{y_3} + v_{33}\frac{y_1y_2}{y_3^2}\right).
\end{gather*}
Adding up the following and multiplying by the common factor $y_3^2$, we obtain \begin{align*}
y_3^2d^2(x_1,x_2) &= (v_{11}+v_{22} - 2v_{12})\\
&- 2(v_{23} - v_{13})\left(\frac{y_2-y_1}{y_3}\right)\\
&+ v_{33}\left(\frac{(y_2-y_1)^2}{y_3^2}\right)
\end{align*}
Now, $v_{11} + v_{22} - 2v_{12} = \displaystyle \Var(Y \1_{X \in [x_1,x_2]})$. Then, $$
(v_{23} - v_{13})\left(\frac{y_2-y_1}{y_3}\right) = \Cov\left(Y\1_{X \in [x_1,x_2]}, Y\left(\frac{y_2-y_1}{y_3}\right)\right).
$$
Finally, $$
v_{33}\left(\frac{(y_2-y_1)^2}{y_3^2}\right) = \Var\left(Y\left(\frac{y_2-y_1}{y_3}\right)\right).
$$
Combining these three terms and adjusting the $y_3^2$ term (and using the identity $\displaystyle \Var(V+W) = \displaystyle \Var(V) + \displaystyle \Var(W) - 2\Cov(V,W)$ for any two random variables $V,W$), we obtain $$
\mathbb E[Y]^2 d^2(x_1,x_2)  = \Var\left(Y\left(\1_{X \in [x_1,x_2]} - \frac{\mathbb E[Y \1_{X \in [x_1,x_2]}]}{\mathbb E[Y]}\right)\right).
$$

Let $A = \{X \in [x_1,x_2]\}$. Since the random variable on the right has expectation zero, its variance is equal to the square of its expectation i.e. \begin{align}\label{eq:st1}
\mathbb E[Y]^2d^2(x_1,x_2) &= \mathbb E\left[Y^2 \left(\1_{A} - \frac{\mathbb E[Y \1_A]}{\mathbb E[Y]}\right)^2\right] \nonumber\\
&= \mathbb E\left[Y^2\1_{A}\right] + \frac{\mathbb E\left[Y^2\right]\mathbb E[Y \1_A]^2}{\mathbb E[Y]^2} - 2 \mathbb E[Y^2\1_A]\frac{\mathbb E[Y \1_A]}{\mathbb E[Y]}\nonumber\\
& = \frac 1{\mathbb E[Y]^2} \left(\mathbb E[Y^2\1_A] \mathbb E[Y]^2 + \mathbb E[Y^2] \mathbb E[Y\1_A]^2 - 2 \mathbb E[Y^2\1_A]\mathbb E[Y\1_A]\mathbb E[Y]\right).
\end{align}

Note that by writing $\mathbb E[Y] = \mathbb E[Y\1_A] + \mathbb E[Y\1_{A^c}]$ and  $\mathbb E[Y^2] = \mathbb E[Y^2\1_A] + \mathbb E[Y^2\1_{A^c}]$, we may break each of the terms into ones involving $A$ and $A^c$. Doing this for the first term,
\begin{equation}\label{eq:st2}
\mathbb E[Y^2\1_A] \mathbb E[Y]^2 = \mathbb E [Y^2\1_A](\mathbb E[Y\1_A]^2 + \mathbb E[Y\1_{A^c}]^2 +2\mathbb E[Y\1_A]\mathbb E[Y\1_{A^c}]).
\end{equation}
For the second term, we have
\begin{equation}\label{eq:st3}
\mathbb E[Y^2] \mathbb E[Y\1_A]^2 = (\mathbb E[Y^2\1_{A}] + \mathbb E[Y^2\1_{A^c}]) \mathbb E[Y\1_A]^2 .
\end{equation}
Finally, for the third equation we have 
\begin{equation}\label{eq:st4}
 2 \mathbb E[Y^2\1_A]\mathbb E[Y\1_A]\mathbb E[Y] = 2 \mathbb E[Y^2\1_{A}]\mathbb E[Y\1_{A}]^2 + 2 \mathbb E[Y^2\1_{A}]\mathbb E[Y\1_A]\mathbb E[Y\1_{A^c}].
\end{equation}

The right hand side of \eqref{eq:st1} equals $\frac 1{\mathbb E[Y]^2} \times ($\eqref{eq:st2} $+$\eqref{eq:st3}$-$\eqref{eq:st4}$)$. After multiple cancellations in the expanded expressions, we obtain$$
\mathbb E[Y]^2d^2(x_1,x_2) = \frac 1{\mathbb E[Y]^2} \left(\mathbb E [Y^2\1_A]\mathbb E [Y\1_{A^c}]^2 + \mathbb E [Y^2\1_{A^c}] \mathbb E [Y\1_A]^2\right).
$$ 
Dividing by $\mathbb E[Y]^2$ gives the conclusion.
\end{proof}

This is followed by the proof of Proposition~\ref{thm:K}.

\begin{proof}[Proof of Proposition~\ref{thm:K}]

By \eqref{eq:K} and Lemma~\ref{lem:rew},
\begin{equation}\label{eq:k1}
K^2 = \sup_{x_1,x_2 \in \mathbb R} \frac 1{\mathbb E[Y]^4} \left(\mathbb E [Y^2\1_A]\mathbb E [Y\1_{A^c}]^2 + \mathbb E [Y^2\1_{A^c}] \mathbb E [Y\1_A]^2\right),
\end{equation}

We will first prove the upper bound. Note that for any set $A$, we have $\mathbb E[Y\1_A] \leq \E[Y]$. Directly applying this bound in the right hand side of \eqref{eq:k1},
$$
K^2 \leq \sup_{x_1,x_2 \in \mathbb R} \frac 1{\mathbb E[Y]^2} \left(\mathbb E [Y^2\1_A]+ \mathbb E [Y^2\1_{A^c}]\right) = M_{\theta}.
$$

We now obtain a lower bound. Since $Y$ is a continuous random variable, by the intermediate value theorem there exist points $x_1,x_2$ such that $\mathbb E(Y\1_{A}) = \frac 12 \mathbb E[Y]$, where we recall that $A = \{X \in [x_1,x_2]\}$. Then, $\mathbb E(Y\1_{A^c}) = \frac 12 \mathbb E[Y]$. This implies that 
$$
d^2(x_1,x_2) = \frac 1{4\mathbb E[Y]^2} \left(\mathbb E [Y^2\1_A]+ \mathbb E [Y^2\1_{A^c}]\right) = \frac {M_{\theta}}4,
$$
and by \eqref{eq:K}, $$
K \geq d(x_1,x_2) = \frac 12\sqrt{M_{\theta}}.
$$
\end{proof}

Finally, we conclude this section with the proof of Lemma~\ref{lem:est}.

\begin{proof}[Proof of Lemma~\ref{lem:est}] 
Suppose that $\frac{\mathbb E[Y\1_A]}{\mathbb E[Y]}  = L \leq \frac 12$. Then, $L \leq 1-L$. The proof now follows from Lemma~\ref{lem:rew} :
\begin{align*}
d^2(x_1,x_2) &= \frac 1{\mathbb E[Y]^4} \left(\mathbb E [Y^2\1_A]\mathbb E [Y\1_{A^c}]^2 + \mathbb E [Y^2\1_{A^c}] \mathbb E [Y\1_A]^2\right) \\
&= \frac{\mathbb E [Y^2\1_A]}{\mathbb E[Y]^2} (1-L)^2 + \frac{\mathbb E [Y^2\1_{A^c}]}{\mathbb E[Y]^2} L^2 \\
& \leq L^2\left(\frac{\mathbb E [Y^2\1_A]}{\mathbb E[Y]^2} + \frac{\mathbb E [Y^2\1_{A^c}]}{\mathbb E[Y]^2}\right) \\
& \leq L^2M_{\theta}. 
\end{align*}
\end{proof}

\subsection{Proof of the variance bound \eqref{dekhna}}

The key inequality which will be used in the proof of the concentration bound \eqref{dekhna} is the Borell-Tsirelson-Ibragimov-Sudakov (Borell-TIS) inequality (cf. \cite[Theorem 2.1.1]{Adler2010}).

\begin{theorem}[Borell-TIS]\label{btis}
Let $\mathscr{G}_t$ be an $R^d$-valued centered Gaussian process with separable index set $T$ which is a.s. bounded. Let $\|f\| = \sup_{t} |\mathscr{G}_t|$. Then, for all $u>0$,
$$
\mathbb P(\|f\| - \E\|f\| > u) \leq e^{-u^2/2\sigma^2},
$$
where $$
\sigma^2 = \sup_{t \in T} \E[|\mathscr{G}_t|^2] = \sup_{t \in T} \Var(\mathscr{G}_t).
$$
\end{theorem}

We are now ready to complete the proof of the concentration bound.

\begin{proof}[Proof of concentration bound in Theorem~\ref{thm:bd}]

We shall now apply this theorem to the process $\mathcal{G}_{\theta}(t)$, whose covariance functional is as in Proposition~\ref{lem:gauss}. Recall that $\sigma^2 = \sup_{t \in \R} \Var(\mathcal G_{\theta}(t))$. However, observe that if $d$ is the distance defined in \eqref{eq:dist} then we have $$
\Var(\mathcal G_{\theta}(t)) = \lim_{x \to \infty} d(t,x).
$$
This is a direct consequence of the dominated convergence theorem and the formula for the covariance functional. By Lemma~\ref{lem:rew} it follows that $\sigma^2 \leq M_{\theta}$.

Therefore, by Theorem~\ref{btis}, if $Z= \sup_{t} |\mathcal G_{\theta}(t)|$, then for any $u>0$,
$$
\Prob(|Z - \E[Z]|>u) \leq e^{-u^2/\sigma^2} \leq e^{-u^2/M_{\theta}}.
$$
This is the statement we wanted to prove.
\end{proof}

\section{Proofs for Section 3}\label{appendix:b}

In this section, we shall prove Theorem~\ref{thm:genbd}. Recall that $X$ is any random variable and $g : \R \to \R$ is a measurable function such that $\E[e^{\eta g(X)}] < \infty$ for all $\eta>0$. For a sequence $\theta_n$, $n \geq 1$ let $F_{\theta_n}$ be the distribution of $X_{\theta_n}$ given by \eqref{twisted1D}. Recall $R_{n,\theta_n}$ from \eqref{empirical1D} and its CDF $F_{n,\theta_n}$.

We will begin with the proof of Proposition~\ref{prop:numden}, which is a standard exercise in unwrapping the definitions of $F_{\theta_n}$ and $F_{n, \theta_n}$.

\begin{proof}[Proof of 
Proposition~\ref{prop:numden}]

Note that for any $q \geq p>0$ and $a,b>0$ we have the inequality $$
\left|\frac{p}{q} - \frac{a}{b}\right| \leq \frac{|p-q|+|a-b|}{b}.
$$
This follows from routine algebra :
\begin{multline*}
\left|\frac{p}{q} - \frac{a}{b}\right| = \left|\frac{pb-aq}{qb}\right| =\left|\frac{p\frac{b}{q}-a}{b}\right| =  \left|\frac{p\frac{b}{q}-p+p-a}{b}\right|\\
= \left|\frac{\frac{p}{q}(b-q)+p-a}{b}\right| \leq \frac{p}{q}\frac{|b-q|}{b}+\frac{|p-a|}{b} \leq \frac{|b-q|+|p-a|}{b}.
\end{multline*}
    
    For any sequence $\theta_n, n \geq 1$, $t>0$ and $x \in \R$, applying this to $|F_{n,\theta_n}(x) - F_{\theta_n}(x)|$ gives 
    \begin{align*}
&\left\{\sup_{x} |F_{n,\theta_n}(x) - F_{\theta_n}(x)| \geq t\right\} \\ =& \left\{\sup_{x}\left|\frac{\frac 1n \sum_{i=1}^n e^{\theta_n g(X_{i})}\1_{X_{i} \leq x}}{\frac 1n \sum_{i=1}^n e^{\theta_n g(X_{i})}} - \frac{\mathbb E[e^{\theta_n g(X)} \1_{X \leq x}]}{\mathbb E[e^{\theta_n g(X)}]}\right| \geq t\right\}\\
\subset & \left\{\sup_{x}\frac{ \left|\frac 1n \sum_{i=1}^n e^{\theta_n g(X_{i})}\1_{X_{i} \leq x}-\mathbb E[e^{\theta_n g(X)} \1_{X \leq x}]\right| + \left|\frac 1n \sum_{i=1}^n e^{\theta_n g(X_{i})}-\mathbb E[e^{\theta_n g(X)}]\right|}{\mathbb E[e^{\theta_n g(X)}]}\geq t\right\}\\
\end{align*}
We must now separate the above set into two parts, one of which involves only the numerator comparison, and the other the denominator comparison. To do this, we note that if $a+b \geq t$ then either $a \geq t/2$ or $b \geq t/2$. Using this fact,
\begin{align*}
&\left\{\sup_{x}\frac{ \left|\frac 1n \sum_{i=1}^n e^{\theta_n g(X_{i})}\1_{X_{i} \leq x}-\mathbb E[e^{\theta_n g(X)} \1_{X \leq x}]\right| + \left|\frac 1n \sum_{i=1}^n e^{\theta_n g(X_{i})}-\mathbb E[e^{\theta_n g(X)}]\right|}{\mathbb E[e^{\theta_n g(X)}]}\geq t\right\} \\
\subset & \left\{\sup_{x}\frac{ \left|\frac 1n \sum_{i=1}^n e^{\theta_n g(X_{i})}\1_{X_{i} \leq x}-\mathbb E[e^{\theta_n g(X)} \1_{X \leq x}]\right|}{\mathbb E[e^{\theta_n g(X)}]}\geq t/2\right\} \\ \cup &\left\{\frac{\left|\frac 1n \sum_{i=1}^n e^{\theta_n g(X_{i})}-\mathbb E[e^{\theta_n g(X)}]\right|}{\mathbb E[e^{\theta_n g(X)}]}\geq t/2\right\} \\
\subset & \left\{\sup_{x}\left|\frac 1n \sum_{i=1}^n e^{\theta_n g(X_{i})}\1_{X_{i} \leq x}-\mathbb E[e^{\theta_n g(X)} \1_{X \leq x}]\right| \geq \frac t2 \mathbb E[e^{\theta_n g(X)}]\right\}\\ \cup& \left\{\left|\frac 1n \sum_{i=1}^n e^{\theta_n g(X_{i})}-\mathbb E[e^{\theta_n g(X)}]\right| \geq \frac t2 \mathbb E[e^{\theta_n g(X)}] \right\}
\end{align*}

Combining the above two set containments and using the union bound, \begin{align*}
&\mathbb P\left[\sup_{x} |F_{n,\theta_n}(x) - F_{\theta_n}(x)| \geq t\right] \\ \leq & \mathbb P\left[\sup_{x}\left|\frac 1n \sum_{i=1}^n e^{\theta_n g(X_{i})}\1_{X_{i} \leq x}-\mathbb E[e^{\theta_n g(X)} \1_{X \leq x}]\right| \geq \frac t2 \mathbb E[e^{\theta_n g(X)}]\right]\\+&\mathbb P \left[\left|\frac 1n \sum_{i=1}^n e^{\theta_n g(X_{i})}-\mathbb E[e^{\theta_n g(X)}]\right| \geq \frac t2 \mathbb E[e^{\theta_n g(X)}] \right]
\end{align*}
This is precisely the statement of the proposition with $T_n = \frac{t}{2}\E[e^{\theta_n g(X)}]$.
\end{proof}

Next, we prove Theorem~\ref{ledvdg}. 

\begin{proof}[Proof of Theorem~\ref{ledvdg}]
By the remark straddled between pages 3 and 4 of \cite{1bc2448d-c21e-3f24-8a81-d9751180f50d}, it is sufficient to prove that if $Y_1,\ldots,Y_n$ are iid random variables with distribution $Y$ such that $\E[Y^2]<\infty$, $x_1,x_2,\ldots,x_m \in \R$ are arbitrary points, and 
\begin{equation}\label{zmn}
Z_{m,n} = \sup_{j=1,\ldots,m} \left|\frac 1n \sum_{i=1}^n Y_i 1_{Y_i \leq x_j} - \mathbb E[Y 1_{Y \leq x_j}]\right|,
\end{equation}
then for every $\epsilon, x>0$,
\begin{equation}\label{tbp}
\mathbb P\left(Z_{m,n} \geq (1+\epsilon)\E[Z_{m,n}]+x\right) \leq \frac{c_{\epsilon}\sqrt{M}}{x \sqrt n}
\end{equation}
where $M = \E[Y^2]$.

Our attempt, for the rest of this proof, is to match our notation to \cite{1bc2448d-c21e-3f24-8a81-d9751180f50d} and use \cite[Corollary 3.1]{1bc2448d-c21e-3f24-8a81-d9751180f50d}. Let \begin{equation}\label{zi}
Z(j) = (Z_i(j))_{i=1,\ldots,n} = \left(Y_i\1_{Y_i \leq x_j} - \E[Y 1_{Y \leq x_j}]\right)_{i=1,\ldots,n}
\end{equation}
be an $n$-dimensional vector depending upon $j$. Let
\begin{equation}\label{sigma}
\sigma = \max_{1 \leq j \leq m} \sqrt{\frac 1n \sum_{i=1}^n \mathbb E[Z_i(j)^2]} = \max_{1 \leq j \leq m} \sqrt{\Var(Y\1_{Y \leq x_j})}.
\end{equation}
Note that $\sigma$ is independent of $n$. Furthermore, by \eqref{zmn} and \eqref{zi},
\begin{equation}\label{zmnzi}
Z_{m,n} = \max_{1 \leq j \leq m}\left|\frac 1n \sum_{i=1}^n Z_i(j)\right|
\end{equation}
is in line with \cite[(c)]{1bc2448d-c21e-3f24-8a81-d9751180f50d}. Define $\mathcal{E}_i =Y_i+\E[Y]$ for $1 \leq i \leq n$, and observe that by \eqref{zi} and \eqref{zmnzi},
$$
|Z_i(j)| \leq \mathcal{E}_i, \quad \mathbb E[\mathcal{E}_i^2] = \E[Y^2]+3\E[Y]^2 \leq M_0^2,
$$
where $M_0 = 2\sqrt{\E[Y^2]}$. Thus,  $\mathcal{E}_i$ and $M_0$ together satisfy \cite[(4)]{1bc2448d-c21e-3f24-8a81-d9751180f50d} with $p=2$.

We have everything in place to apply \cite[Corollary 3.1]{1bc2448d-c21e-3f24-8a81-d9751180f50d} with $l=1$ and $p=2$ on the right hand side. Doing so, we have for every $\epsilon, x>0$ that 
$$
\mathbb P\left(Z_{m,n} \geq (1+\epsilon)\E[Z_{m,n}]+x\right) \leq \frac{\left(\frac{64}{\epsilon}+7+\epsilon\right)\frac{M_0}{\sqrt{n}} + \frac{4}{n\sigma}}{x}.
$$
Note that as $n \to \infty$, $\frac{M_0/\sqrt{n}}{4/(n \sigma)} \to +\infty$. Therefore, there exists a constant $c_n>0$ such that $\frac{M_0}{\sqrt{n}} \geq C \frac{4}{n \sigma}$ for all $n\geq 1$, which implies along with $M_0 = 2 \sqrt{\E[Y^2]}$ that $$
\mathbb P\left(Z_{m,n} \geq (1+\epsilon)\E[Z_{m,n}]+x\right) \leq c_{\epsilon}\frac{\sqrt{\E[Y^2]}}{x\sqrt{n}},
$$
which matches \eqref{tbp} that was to be proved.
\end{proof}

Finally, we conclude the proof of Theorem~\ref{thm:genbd}.

\begin{proof}[Proof of Theorem~\ref{thm:genbd}]

By definition of convergence in probability and the $\KS$ distance \eqref{KS}, we must prove that for every $t>0$, \begin{equation}\label{inprob}
s_n\mathbb P\left[\sup_{x} |F_{n,\theta_n}(x) - F_{\theta}(x)| \geq t\right] \to 0
\end{equation}
as $n \to \infty$, if $s_n^2 \frac{M_{\theta_n}}{n} \to 0$. Ignoring the $s_n$ for now, we directly apply Proposition~\ref{prop:numden} to the rest of the left hand side :
\begin{align}
&\mathbb P\left[\sup_{x} |F_{n,\theta_n}(x) - F_{\theta_n}(x)| \geq t\right]  \\ \leq & \mathbb P\left[\sup_{x}\left|\frac 1n \sum_{i=1}^n e^{\theta_n g(X_{i})}\1_{X_{i} \leq x}-\mathbb E[e^{\theta_n g(X)} \1_{X \leq x}]\right| \geq \frac t2 \mathbb E[e^{\theta_n g(X)}]\right]\nonumber\\+&\mathbb P \left[\left|\frac 1n \sum_{i=1}^n e^{\theta_n g(X_{i})}-\mathbb E[e^{\theta_n g(X)}]\right| \geq \frac t2 \mathbb E[e^{\theta_n g(X)}] \right]\label{propbd}
\end{align}
At this point we just need to estimate the two terms on the right hand side above. The second one is much easier and can be bounded by Chebyshev's inequality : \begin{equation}
\mathbb P \left[\left|\frac 1n \sum_{i=1}^n e^{\theta_n g(X_{i})}-\mathbb E[e^{\theta_n g(X)}]\right| \geq \frac{t}{2} \mathbb E[e^{\theta_n g(X)}] \right] \nonumber
\leq  \frac{\Var\left(\frac 1n \sum_{i=1}^n e^{\theta_n g(X_{i})}\right)}{\frac {t^2}4\mathbb E[e^{\theta_n g(X)}]^2}\nonumber \leq  \frac {4M_{\theta_n}}{n t^2}.\label{cheby}
\end{equation}

For the first term on the right hand side of Proposition~\ref{prop:numden}, note that $\mathbb E[Z_n] \to 0$ since $Z_n \overset{d}{\to} 0$ by the Glivenko Cantelli lemma. Thus, we may apply Theorem~\ref{ledvdg} with $Y = e^{\theta_n g(X)}$, $\epsilon=1$ and $x = \frac{t}{2}\E[Y] - 2\mathbb{E}[Z_n]$ (which is positive for $n$ large enough), where $Z_n$ is as in the theorem. This gives $$
\mathbb P\left(Z_n \geq \frac{t}{2}\E[Y]\right) \leq \frac{c}{\left( \frac{t}{2}\E[Y] - 2\mathbb{E}[Z_n]\right)\sqrt{n}} \sqrt{\E[Y^2]} = \frac{c}{\sqrt{n}} \sqrt{\frac{\E[Y^2]}{\left( \frac{t}{2}\E[Y] - 2\mathbb{E}[Z_n]\right)^2}}
$$
for some constant $c_n>0$ independent of $X,\theta$ and $n$.

Now, let $s_n$ be a sequence satisfying \eqref{decay}. Combining the above inequality with \eqref{cheby} and \eqref{propbd}, and multiplying by $s_n$ on both sides, \begin{equation*}
s_n\mathbb P\left[\sup_{x} |F_n(x) - F_{\theta_n}(x)| \geq t\right]  \leq s_n\frac{4M_{\theta_n}}{nt^2} + s_n\frac{c}{\sqrt{n}} \sqrt{\frac{\E[Y^2]}{\left( \frac{t}{2}\E[Y] - 2\mathbb{E}[Z_n]\right)^2}}.
\end{equation*}

 The first term above converges to $0$ by \eqref{decay}, since $\frac{M_{\theta_n}}{n} \to 0$. Consequently, \begin{multline*}\sqrt{\frac{\E[Y^2]}{\left( \frac{t}{2}\E[Y] - 2\mathbb{E}[Z_n]\right)^2}} \times \frac 1{\frac{2}{t} \sqrt{M_{\theta_n}}} \to 1 \\ \implies s_n\frac{c}{\sqrt{n}} \sqrt{\frac{\E[Y^2]}{\left( \frac{t}{2}\E[Y] - 2\mathbb{E}[Z_n]\right)^2}} \times \frac {\sqrt{n}}{\frac 2t s_n c \sqrt{M_{\theta_n}}} \to 1\end{multline*}
 as $n \to \infty$. Now, since $s_n$ satisfies \eqref{decay}, $$\frac {\sqrt{n}}{\frac 2t s_n c \sqrt{M_{\theta_n}}} \to +\infty,$$
which implies that the other term must converge to $0$ if the product is to converge to $1$.

 Thus, \eqref{inprob} has been proved for arbitrary $t>0$, completing the proof.
\end{proof}

\section{Auxiliary material for the proofs in Section~\ref{1DWeibull}, \ref{mdweibull} and \ref{unbdd}}\label{appendix:c}

In this section, we will collect results that apply to all the upcoming sections. We divide these into three parts. The first part will be used to prove Lemma~\ref{KStoSL} and Lemma~\ref{KStoSLhd}. The second will be used to prove the positive results in each regime i.e. when $\frac{M_{\theta_n}}{n} \to 0$. In the final part we will define Poisson random measures and prove some key lemmas which will allow us to isolate the common aspects of the proofs in the regimes where $\frac{M_{\theta_n}}{n} \not \to 0$.

Throughout this section, let $X$ be a random vector in $\R^d, d \geq 1$ (note : we include the case $d=1$ where $X$ is just a random variable).

\subsection{Strengthening convergence of random vectors}

We begin with the following lemma which is key to establishing the proofs of Lemma~\ref{KStoSL} and Lemma~\ref{KStoSLhd}. Recall the $\textsf{KS}$ distance from \eqref{KS} and the notion of uniform convergence on compacts from Lemma~\ref{KStoSLhd}. The first lemma strengthens convergence in distribution to convergence in these distances, under the assumption that the limiting distribution is continuous. The second asserts that two random variables are different if and only if their CDFs differ at a common point of continuity.

\begin{lemma}\label{assist}\begin{enumerate}[label = (\alph*)]
\item Suppose that $Y_n$ is a sequence of random vectors which converge, in distribution, to a continuous random vector $Z$. Then, if $Z$ is continuous, we have $\textsf{KS}(Y_n,Z) \to 0$ (if $d=1$), and uniform convergence over compacts of $F_{Y_n}$ to $F_Z$ (if $d >1$). 
\item If $Z_1, Z_2$ are random vectors and if $F_{Z_1}(c) = F_{Z_2}(c)$ at every $c \in \R^d$ at which both $F_{Z_1},F_{Z_2}$ are continuous, then $Z_1\overset{d}{=}Z_2$.
\end{enumerate}
\end{lemma}
\begin{proof}
We begin with the proof of part (a). For the KS distance, the argument is as in \cite[page 3]{Resnick1987}. We can adapt this argument to uniform convergence on compacts as well. 

Let $\mathcal R \subset \mathbb R^d$ be any axis-aligned compact hypercube. Since $Z$ is continuous, we have $F_{X_n}(c) \to F_{Z}(c)$ for every $c \in \mathcal{R}$. We claim that, in fact, the convergence above is uniform over all $c\in \mathcal{R}$.

To prove this, note that $\mathcal{R}$ is compact, hence $F_Z$ is uniformly continuous over $\mathcal{R}$. Given $\epsilon>0$, we can divide $\mathcal{R}$ into finitely many axis-aligned hypercubes $\mathcal{R}_i, 1 \leq i \leq M$ such that if $x,y \in \mathcal{R}_i$ then $|F_Z(x) - F_Z(y)| < \epsilon$. Let $x_i,y_i \in \mathcal{R}_i$ be the lower left and upper right endpoints of $\mathcal{R}_i$.

Let $N$ be chosen large enough so that if $n>N$, then $|F_{X_n}(z) - F_{Z}(z)| < \epsilon$ for all $z = x_i,y_i, 1 \leq i \leq M$. Now, if $t \in \mathcal{R}$, then $t\in \mathcal{R}_i$ for some $i$. If $n>N$ then \begin{align*}
&\|F_{X_n}(y_i) - F_{X_n}(x_i)\|\\ \leq& \|F_{X_n}(y_i) - F_Z(y_i)\| + \|F_Z(y_i) - F_Z(x_i)\| + \|F_Z(x_i) - F_{X_n}(x_i)\| \\ <& 3\epsilon.
\end{align*}
Observe that $F_{X_n}(t_i) \in [F_{X_n}(x_i), F_{X_n}(y_i)]$ and $F_{Z}(t_i) \in [F_Z(x_i), F_Z(y_i)]$. Therefore, \begin{gather*}\|F_{X_n}(t_i) - F_{X_n}(x_i)\| \leq \|F_{X_n}(y_i) - F_{X_n}(x_i)\| < 3\epsilon \\ \|F_{Z}(t_i) - F_{Z}(x_i)\| \leq \|F_{Z}(y_i) - F_{Z}(x_i)\| < \epsilon\end{gather*}
Therefore we have
\begin{align*}
&\|F_{X_n}(t_i) - F_{Z}(t_i)\| \\ \leq &\|F_{X_n}(t_i) - F_{X_n}(x_i)\| + \|F_{X_n}(x_i) - F_{Z}(x_i)\| + \|F_Z(x_i) - F_Z(t_i)\| \\ <& 5\epsilon
\end{align*}
whenever $n>N$. Thus, $F_{X_n} \to F_Z$ uniformly over $\mathcal{R}$. Since any compact set is contained in a large enough compact axis-aligned hypercube, the proof is complete.

Part (b) directly follows from a particular kind of right-continuity that multivariate CDFs possess. Indeed, if $x_n\in \R^d$ converges to $x$ and $x_n$ is coordinatewise bigger than $x$ for all $n$, then $F(x_n) \to F(x)$ for any multivariate CDF $F$ (see \cite[Theorem 4.25]{KAL}). It follows from right-continuity and monotonicity of CDFs that their discontinuities form a set of measure zero. That is, given any two random vectors $Z_1,Z_2$, their CDFs $F_{Z_1},F_{Z_2}$ are simultaneously continuous outside a set of measure zero.

Thus, for any $c\in \R^d$ we can find points $c_n \to c$ such that $F_{Z_1}, F_{Z_2}$ are both continuous at $c_n$ for all $n$, and $c_n$ is coordinatewise greater than $c$. If we assume that $F_{Z_1}(c) = F_{Z_2}(c)$ at all points of continuity, then we have $F_{Z_1}(c_n) = F_{Z_2}(c_n)$ and hence by right continuity we get $F_{Z_1}(c) = F_{Z_2}(c)$, as desired.
\end{proof}

\subsection{A key lemma in establishing accuracy in the regime $\frac{M_{\theta_n}}{n} \to 0$}

In this section we will outline the key proposition which will be used to prove part(a) of Theorems~\ref{thm:slemp1D}, \ref{thm:emphd} and \ref{empunbdd}. As explained prior to the statement of Theorem~\ref{thm:slemp1D}, the idea is that for $\theta_n$ small, the numerators of \eqref{twisted} and \eqref{empirical} are asymptotically equivalent, likewise the denominators.

We make this precise in the following proposition. Recall $M_{\theta}$ from \eqref{m2tmt2}. We assume that $g(x)=x$ throughout this section, since Lemmas~\ref{weibtoweib1D} and \ref{weibtoweibhd} allow us to replace $X$ with $g(X)$ in all arguments that follow.

\begin{proposition}\label{posassist}
Let $X$ be a random vector, and suppose $\theta_n$ be a sequence of vectors such that $\frac{M_{\theta_n}}{n} \to 0$. Let $B_n$ be a sequence of measurable sets such that 
\begin{equation}\label{hyp}
\limsup_{n \to \infty} \frac{\E[e^{\theta_n^T X}]}{\E[e^{\theta_n^T X}\1_{X \in B_n}]} <\infty.
\end{equation}
Then, $$
\frac{\mathbb P(R_{n,\theta_n} \in B_n)}{\mathbb P(X_{\theta_n} \in B_n)} \to 1
$$
in probability, as $n \to \infty$. In particular, if 
\begin{enumerate}[label = (\alph*)]
\item There exist scalars/vectors $a_n$ and scalars $b_n$ such that $\frac{X_n-a_n}{b_n} \overset{d}{\to} Z$, and 
\item For every continuity point $x$ of $F_Z$, the sets $B_n = \{X \leq b_nx+a_n\}$ satisfy the hypothesis \eqref{hyp}. (Here, $\leq$ is component-wise for vectors)
\end{enumerate}
Then, $\frac{R_{n,\theta_n} - a_n}{b_n} \overset{d}{\to} Z$. 
\end{proposition}
\begin{proof}
Observe that by \eqref{twisted} and \eqref{empirical}, 
\begin{equation}\label{combine}
\frac{\mathbb P(R_{n,\theta_n} \in B_n)}{\mathbb P(X_{\theta_n} \in B_n)} = \frac{\frac 1n\sum_{i=1}^n e^{\theta_n^TX_{i}}\1_{X_{i} \in B_n}}{\E[e^{\theta_n^T X}\1_{X\in B_n}]}\frac{\E[e^{\theta_n^T X}]}{\frac 1n\sum_{i=1}^n e^{\theta_n^TX_{i}}}.
\end{equation}
We will prove that each ratio converges to $1$ in probability, following which their product will also do so. Fix $\epsilon>0$ and $n \geq 1$. In the following argument, let $A_n = B_n$ for all $n \geq 1$ or $A_n = \R^d$ for all $n \geq 1$. By Chebyshev's inequality , \begin{align*}
&\mathbb P\left[\left|\frac{\frac 1n\sum_{i=1}^n e^{ \theta_n^TX_{i}}\1_{X_{i} \in A_n}}{\E[e^{ \theta_n^T X}\1_{X\in A_n}]}-1\right| > \epsilon\right] \\
= &\mathbb P\left[\left|\frac 1n \sum_{i=1}^n e^{ \theta_n^TX_{i}}\1_{X_{i} \in A_n} -\E[e^{\theta_n^T X}\1_{X\in A_n}]\right| > \epsilon\E[e^{\theta_n^T X}\1_{X\in A_n}]\right] \\
\leq & \frac{\Var\left(e^{\theta_n^T X}\1_{X\in A_n}\right)}{n \epsilon^2 \E\left[e^{\theta_n^T X}\1_{X\in A_n}\right]^2} \\
\leq &  \frac{\Var\left(e^{\theta_n^T X}\1_{X\in A_n}\right)}{n \epsilon^2 \E\left[e^{\theta_n^T X}\1_{X\in A_n}\right]^2} \\
\leq & \frac{\E\left[e^{2\theta_n^T X}\right]}{n \epsilon^2 \E\left[e^{ \theta_n^T X}\1_{X\in A_n}\right]^2} \\
\leq & \frac 1{\epsilon^2}\frac{\E[e^{2\theta_n^T X}]}{n\E[e^{\theta_n^T X}]^2} \frac{\E[e^{\theta_n^T X}]^2}{\E[e^{\theta_n^T X}\1_{X \in A_n}]^2} = \frac{1}{\epsilon^2} \frac{M_{\theta_n}}{n} \frac{\E[e^{\theta_n^T X}]^2}{\E[e^{\theta_n^T X}\1_{X \in A_n}]^2}.
\end{align*}
By our assumptions, it follows that this quantity converges to $0$ as $n \to \infty$ for all $\epsilon>0$, which implies that $$
\frac{\frac 1n\sum_{i=1}^n e^{ \theta_n^TX_{i}}\1_{X_{i} \in A_n}}{\E[e^{ \theta_n^T X}\1_{X\in A_n}]} \to 1
$$ 
in probability. Applying this result for $A_n=B_n$ and then $A_n=\mathbb R^d$, and multiplying the resulting statements together proves that the terms in \eqref{combine} converge to $1$ in probability, as desired.

That $\frac{X_n-a_n}{b_n} \overset{d}{\to} Z$ in distribution is tantamount to saying that $F_{\frac{X_n-a_n}{b_n}}(x) \to F_Z(x)$ at all continuity points $x$ of $F_Z$, which is true by assumption (a). Observe that $$
\frac{X_n - a_n}{b_n} \leq x \iff X_n \leq b_nx+a_n
$$
However, assumption (b) guarantees along with the previous part that $$
\frac{\mathbb P(X_n \leq b_nx+a_n)}{\mathbb P(R_{n,\theta_n} \leq b_nx+a_n)} \to 1
$$
in probability. The conclusion follows easily from here.
\end{proof}

\subsection{Point Processes, Poisson random measures and the regimes $\frac{M_{\theta_n}}{n}\not \to 0$}

Next, we introduce point processes and the Poisson random measure, which will be used to explicitly find the limiting random variables in part (b) in each of Theorems~\ref{thm:slemp1D}, \ref{thm:emphd} and \ref{empunbdd}. It turns out that they will be essential to the proofs in part (c) of these theorems as well. We will establish some key lemmas which will be used to prove these results together, but the proofs of the theorems themselves will be relegated to future sections. We take our material from \cite[Chapter 3]{Resnick1987}.

Let $\delta_x$ denote the Dirac measure $\delta_x(A) = \1_{x \in A}$. A point measure $M$ on $E \subset \R^d,d \geq 1$ is a measure of the form $M= \sum_{i} \delta_{x_i}$ where $\{x_i\} \subset E$  is an at-most countable collection of points (which can also be finite e.g. the Dirac measure, and can also contain the same point multiple times e.g. twice the Dirac measure). Let $\mathcal M_p(E)$ denote the set of all point measures on $E$. It is equipped with the coarsest topology such that for all Borel sets $A \subset E$, the evaluation map $e_A : \mathcal M_p(E) \to [0,\infty]$ given by $e_A(M) = M(A)$ is measurable. It turns out that this topology coincides with vague convergence of measures, restricted to $\mathcal M_p(E)$, which we recall from the discussion before Assumption~\ref{ass:mvrv}.

Any random element taking values in $\mathcal M_{p}(E)$ is called a point process on $E$. It can thus be thought of as a random collection of points on $E$.

\begin{definition}\label{PRM}
Let $E\subset \mathbb R^d$ and $\nu$ be a Radon measure on $E$. A point process $M$ taking values on $E$ is called a Poisson random measure, or PRM with intensity measure $\mu$ if 
\begin{enumerate}[label = (\alph*)]
\item For all $A \subset E$, $M(A)$ is distributed as a Poisson random variable with parameter $\nu(A)$. (Note : If $\nu(A)=0$ then $M(A)=0$ and if $\nu(A) = \infty$ then $M(A) = \infty$ a.s.)
\item For $A_1,A_2,\ldots,A_n \subset E$ mutually disjoint, $\{M(A_i)\}_{i=1,2,\ldots,n}$ forms an independent collection of random variables.
\end{enumerate}
\end{definition}

By \cite[Proposition 3.6]{Resnick1987}, a PRM on any Borel subset of $\R^d$ with any intensity measure $\nu$ that is Radon, always exists and is unique up to distribution. We denote a PRM with intensity measure $\nu$ by $PRM(\nu)$, and integration with respect to $PRM(\nu)$ will be denoted by $dPRM(\nu)$ with the integration variable being omitted (e.g. $\int f(x) dPRM(\nu)$).

The convergence of a sequence of point processes $M_{n} \in \mathcal M_p(E)$ to a point process $M \in \mathcal M_p(E)$ will be considered in the weak sense here, see \cite[Section 3.5]{Resnick1987} for more details. We use $\Rightarrow$ to denote this convergence. Note that this is not an abuse of notation; the same notion of convergence was used for  Gaussian process limits as well earlier.

The following results will now be used in the proofs of various point-process related convergences in our work. We begin with the following vague convergence result.

\begin{lemma}\label{pppc}
Let $X_i$ be iid random vectors taking values in a Borel $E \subset \mathbb R^d, d \geq 1$ with the same distribution $X$. Suppose that for some sequence of vectors $a_n$ and scalars $b_n$ we have 
$$
n \P\left(\frac{X-a_n}{b_n} \in \cdot\right) \to \nu\left(\cdot\right)
$$
vaguely (see Assumption~\ref{ass:mvrv} for the definition of vague convergence) for some Radon measure $\nu$ on $E$. Then, $$
\sum_{k=1}^n \delta_{(X_k-a_n)/b_n} \Rightarrow PRM(\nu)
$$
on $\mathcal M_p(E)$.
\end{lemma}
\begin{proof}
This is a simple application of the "warm-up" exercise in the proof of \cite[Proposition 3.21]{Resnick1987} with $X_{n,j} = \frac{X_j-a_n}{b_n}, n,j \geq 1$.
\end{proof}

We will now setup a common lemma which converts point process convergence to convergence in distribution. Since this may require case-by-case analysis, we set the tone. 

\begin{lemma}\label{dekhliyo}
Let $X$ be a random vector in $\R^d, d \geq 1$ and $\theta \in \R^d$ be a unit vector. Let $c_n$ be a sequence of positive real numbers and $\theta_n = c_n \theta$. Suppose that, for some positive real sequence $a_n$ and vectors $x_n$, bounded Borel set $D \in \R^d$ and Radon measure $\nu$ on $\R^d$, all the following conditions hold.
\begin{enumerate}[label = (\alph*)]
\item $c_na_n \to C \in (0,\infty)$, and $c_n \to \infty$. (Hence, $a_n \theta_n \to C\theta$)
\item $n \mathbb P\left(\frac{(x_n-X)}{a_n} \in \cdot\right) \to \nu(\cdot)$ vaguely.
\item $\frac{x_n-A_n}{a_n} \to D$ as sets, as $n \to \infty$.
\end{enumerate}
Then,
\begin{multline*}
\left(\sum_{i=1}^n e^{-\theta_n^T(x_n-X_i)}\1_{X_i \in A_n}, \sum_{i=1}^n e^{-\theta_n^T(x_n-X_i)}\right)\\ \overset{d}{\to} \left(\int e^{-C\theta^T y}\1_{y\in D} dPRM(\nu),\int e^{-C\theta^T y}\1_{y\in D} dPRM(\nu)\right).
\end{multline*}
\end{lemma}

\begin{proof}
We will show that this convergence holds using the Laplace transform (see \cite[Theorem 6.3]{KAL}). In order to do this, let $s_1,s_2<0$ be arbitrary and consider the quantity
$$
\mathbb E\left[\exp\left(s_1\sum_{i=1}^n e^{-\theta_n^T(x_n-X_i)}1_{A_n}+s_2\sum_{i=1}^n e^{-\theta_n^T(x_n-X_i)}\right)\right]
$$
Begin from the left hand side :
\begin{align*}
&\mathbb E\left[\exp\left(s_1\sum_{i=1}^n e^{-\theta_n^T(x_n-X_i)}1_{A_n}+s_2\sum_{i=1}^n e^{-\theta_n^T(x_n-X_i)}\right)\right] \\
=& \mathbb E\left[\exp\left(\sum_{i=1}^n \left (s_1e^{-\theta_n^T(x_n-X_i)}1_{A_n}+s_2e^{-\theta_n^T(x_n-X_i)}\right)\right)\right] \\
 = &\mathbb E\left[\exp\left(s_1e^{-\theta_n^T(x_n-X_1)}1_{A_n}+s_2e^{-\theta_n^T(x_n-X_1)}\right)\right]^n \\=& \left(1- \mathbb E\left[1-\exp\left(s_1e^{-\theta_n^T(x_n-X_1)}1_{A_n}+s_2e^{-\theta_n^T(x_n-X_1)}\right)\right]\right)^n.
\end{align*}
The term inside the expectation goes to zero as $n \to \infty$, since $\|\theta_n\| =c_n \to \infty$ by assumption (a). Therefore, \begin{equation}\label{thr}
\lim_{n \to \infty} \frac{\left(1- \mathbb E\left[1-\exp\left(s_1e^{-\theta_n^T(x_n-X_1)}1_{A_n}+s_2e^{-\theta_n^T(x_n-X_1)}\right)\right]\right)^n}{\exp\left(-n \mathbb E\left[1-\exp\left(s_1e^{-\theta_n^T(x_n-X_1)}1_{A_n}+s_2e^{-\theta_n^T(x_n-X_1)}\right)\right]\right)} = 1.
\end{equation}
It suffices, therefore, to focus on the limit of $$
n \mathbb E\left[1-\exp\left(s_1e^{-\theta_n^T(x_n-X_1)}1_{A_n}+s_2e^{-\theta_n^T(x_n-X_1)}\right)\right].
$$
Let $R>0$ be arbitrary, and $$Y_n = \frac{(x_n-X_1)}{a_n}$$ for notational convenience. By assumption (b), we have $n \mathbb P(Y_n \in \cdot) \to \nu(\cdot)$ vaguely. Now,
\begin{align}
    &n \mathbb E\left[1-\exp\left(s_1e^{-\theta_n^T(x_n-X_1)}1_{A_n}+s_2e^{-\theta_n^T(x_n-X_1)}\right)\right] \nonumber\\=& 
    n \mathbb E\left[1-\exp\left(s_1e^{-a_n\theta_n^T y}1_{y \in \frac{x_n-A_n}{a_n}}+s_2e^{-a_n\theta_n^Ty}\right)\right]  \label{zero1}
\end{align}

Let $R>0$ be an arbitrary real number. Note that
\begin{align}
    &n \mathbb E\left[1-\exp\left(s_1e^{-a_n\theta_n^T y}1_{y \in \frac{x_n-A_n}{a_n}}+s_2e^{-a_n\theta_n^Ty}\right)\right] \nonumber\\
    \geq & n \int_{C_1\theta^Ty\leq R, \|y\| \leq R}\left[1-\exp\left(s_1e^{-a_n\theta_n^T y}1_{y \in \frac{x_n-A_n}{a_n}}+s_2e^{-a_n\theta_n^Ty}\right)\right] d\mathbb P_{Y_n}(y) \label{zero2}
\end{align}

For $R$ fixed, this term converges due to vague convergence of $Y_n$, and the fact that we are restricted to a compact set (a small approximation argument is required, whose basic idea is derived from the end of the proof of \cite[Proposition 3.21]{Resnick1987}.)

\begin{multline}\label{firstone}
 n \int_{C_1 \theta^T y\leq R}\left[1-\exp\left(s_1e^{-a_n\theta_n^T y}1_{y \in \frac{x_n-A_n}{a_n}}+s_2e^{-a_n\theta_n^Ty}\right)\right] d\mathbb P_{Y_n}(y) \\ \to \int_{C_1 \theta^T y\leq R, \|y\| \leq R}\left[1-\exp\left(s_1e^{-C_1\theta^T y}1_{y \in D}+s_2e^{-C_1\theta^Ty}\right)\right] d\nu(y).
\end{multline}
Combining this statement and \eqref{zero2} we have as $R \to \infty$ that 
\begin{align}
&\liminf_{n \to \infty} n \mathbb E\left[1-\exp\left(s_1e^{-a_n\theta_n^T y}1_{y \in \frac{x_n-A_n}{a_n}}+s_2e^{-a_n\theta_n^Ty}\right)\right]\nonumber \\ \geq & \int \left[1-\exp\left(s_1e^{-C_1\theta^T y}1_{y \in D}+s_2e^{-C_1\theta^Ty}\right)\right] d\nu(y).\label{lowerbd}
\end{align}

On the other hand, suppose that $R_n$ is any sequence of real numbers increasing to $+\infty$. We have 
\begin{align*}
& n \mathbb E\left[1-\exp\left(s_1e^{-a_n\theta_n^T y}1_{y \in \frac{x_n-A_n}{a_n}}+s_2e^{-a_n\theta_n^Ty}\right)\right] \\ 
 =& \int_{\theta_n^T y\leq R_n} \left[1-\exp\left(s_1e^{-a_n\theta_n^T y}1_{y \in \frac{x_n-A_n}{a_n}}s_2e^{-a_n\theta_n^Ty}\right)\right] d\nu(y)\\ +& \int_{\theta_n^T y> R_n} \left[1-\exp\left(s_1e^{-a_n\theta_n^T y}1_{y \in \frac{x_n-A_n}{a_n}}+s_2e^{-a_n\theta_n^Ty}\right)\right] d\nu(y) \\
 \leq & \int \left[1-\exp\left(s_1e^{-a_n\theta_n^T y}1_{y \in \frac{x_n-A_n}{a_n}}+s_2e^{-a_n\theta_n^Ty}\right)\right] d\nu(y) \\ +& n \mathbb P(\theta_n^T y> R_n) (-s_1-s_2)e^{-a_nR_n},
\end{align*}
where we used the inequality
\begin{align*}1-\exp\left(s_1e^{-a_n\theta_n^T y}1_{y \in \frac{x_n-A_n}{a_n}}+s_2e^{-a_n\theta_n^Ty}\right) \leq &-\left(s_1e^{-a_n\theta_n^T y}1_{y \in \frac{x_n-A_n}{a_n}}+s_2e^{-a_n\theta_n^Ty}\right)\\ \leq & -(s_1+s_2)e^{-a_n\theta_n^T y}\end{align*}
Let $R_n$ be a sequence chosen so large that $n \mathbb P(\theta_n^T Y > R_n) \to 0$. Then, it follows by letting $n \to \infty$ above that \begin{multline*}
\limsup_{n \to \infty} n \mathbb E\left[1-\exp\left(s_1e^{-a_n\theta_n^T y}1_{y \in \frac{x_n-A_n}{a_n}}+s_2e^{-a_n\theta_n^Ty}\right)\right] \\ \geq \int \left[1-\exp\left(s_1e^{-C_1\theta^T y}1_{y \in D}+s_2e^{-C_1\theta^Ty}\right)\right] d\nu(y).
\end{multline*}
Combining this with the lower bound \eqref{lowerbd} we have \begin{align}\label{try}
&\lim_{n \to \infty} n \mathbb E\left[1-\exp\left(s_1e^{-a_n\theta_n^T y}1_{y \in \frac{x_n-A_n}{a_n}}+s_2e^{-a_n\theta_n^Ty}\right)\right]\nonumber \\ =& \int \left[1-\exp\left(s_1e^{-C_1\theta^T y}1_{y \in D}+s_2e^{-C_1\theta^Ty}\right)\right] d\nu(y).
\end{align}

It now suffices to look at the Laplace transform at $(-s_1,-s_2)$ of the left hand side of the lemma. That is, 
\begin{multline*}
\mathbb E\left[\exp\left(s_1 \int e^{-C_1\theta^Ty}1_{y \in D} d PRM(\nu) + s_2 \int e^{-C_1\theta^Ty} dPRM(\nu)\right)\right] \\ = \mathbb E\left[\exp\left( \int (s_1e^{-C_1\theta^Ty}1_{y \in D}+s_2e^{-C_1\theta^T y}) d PRM(\nu)\right)\right].
\end{multline*}
However, this is precisely the Laplace functional of $PRM(\nu)$ evaluated at the function $f(s) = -(s_1e^{-C_1y}1_{y \in D}+s_2e^{-C_1y})$ (see \cite[Section 3.2]{Resnick1987} for the definition). This admits an explicit formula for $PRM(\nu)$ by \cite[Proposition 3.6(ii)]{Resnick1987}. Applying that formula,
$$
\mathbb E\left[\exp\left(- \int f(s) d PRM(\nu)\right)\right] = \exp\left(-\int_{0}^{\infty} (1-e^{-f(s)})d\nu(s)\right)
$$
which matches exactly with the right hand side of \eqref{try}, once we use \eqref{thr}. Thus, the proof is complete.
\end{proof}

Having proved this important lemma, we shall now prove a result that will help us with part (c) of Theorems~\ref{thm:slemp1D} and ~\ref{thm:emphd}. We proceed to motivate the lemma, state it and then prove it. Suppose for ease of clarity that out of the iid samples $X_1,X_2,\ldots,X_n$ of $X$, that $X_1$ is the sample maximizer of the quantity $\theta^T x$, for some fixed unit vector $\theta$. Consider the quantity \begin{align*}
\frac{\sum_{i=1}^n e^{\theta_n^T X_{i}}}{e^{\theta_n^T X_1}} &= \sum_{i=1}^n e^{\theta_n^T (X_{i}-X_1)}
\end{align*}

Suppose that a particular functional of $X$ converges to a PRM, for instance let $x\in \R^d$ and $a_n$ be a sequence such that $\sum_{i=1}^n \delta_{\frac{x-X_i}{a_n}} \Rightarrow PRM( \nu)$ for some Radon measure $\nu$, and $\theta_n$ be a sequence of vectors such that $a_n\theta_n \to \theta$ for some vector $\theta$. Observe that $$
\sum_{i=1}^n e^{\theta_n^T (X_{i}-X_1)} = \sum_{i=1}^n e^{a_n\theta_n^T\frac{(x - X_1)}{a_n} - a_n\theta_n^T\frac{(x - X_i)}{a_n}}
$$

The term on the exponential is a random measure we have already studied, minus its very first term. Since $a_n \theta_n \to \theta$, we expect that the rest will converge to a PRM minus its first term, which is a phenomena general enough to warrant its own lemma. This is the lemma that we will now proceed to prove. A simple point process is one which, when evaluated at a singleton $\{a\}$ for any $a \in \R$, results in either $0$ or $1$. 
 
\begin{lemma}\label{pppconv}
Let $\eta_n$ be a sequence of simple point processes on $\mathbb R^d$ such that $\eta_n$ has $n$ points. Suppose $\eta_n \Rightarrow \eta$. For any simple point process $\mu$ on $\R^d$ with finitely many points, let $s_{\mu} = \argmin \theta^T y$ be the sample minimizer of $\theta^T y$ in $\mu$. Suppose $s_{\mu}$ is unique and isolated a.s. in $\mu$. then, $\eta_n - s_{\eta_n} \Rightarrow \eta - s_{\eta}$.
\end{lemma}
\begin{proof}
We will show that the map $\mu \to \mu - s_\mu$ is continuous on the subspace of simple point processes. Then, this result follows directly from the continuous mapping theorem. In order to do this, we must first prove that $\mu \to s_{\mu}$ is a continuous real-valued function on the space of all $\mu$ with a unique minimizer of $\theta^T y$.

Let $\mu_n$ converge to $\mu$ vaguely. We will prove that $s_{\mu_n}$ converges to $s_{\mu}$. Let $R_1 = \theta^T s_{\mu}$ and $R_2 = \|s_{\mu}\|$. Note that for any $r_1,r_2>0$, we have $$
\mu_n(\{\theta^Ty < r_1, \|y\| \leq r_2\}) \to \mu(\{\theta^Ty < r_1, \|y\| \leq r_2\})
$$
Let $\epsilon>0$ be chosen small enough so that $\theta^Ts_{\mu} < \inf_{\mu\setminus s_{\mu}} \theta^Ty - \epsilon$. Then, for any $\epsilon'>0$, 
$$
\mu(\{\theta^Ty < R-\epsilon, \|y\| \leq R_2-\epsilon'\}) = 0,\mu(\{\theta^Ty < R+\epsilon, \|y\| \leq R_2+\epsilon'\}) = 1.
$$
Thus, for large enough $n$ we have by vague convergence that 
$$
\mu_n(\{\theta^Ty < R-\epsilon, \|y\| \leq R_2\}) =0,\mu_n(\{\theta^Ty < R+\epsilon, \|y\| \leq R_2\}) \to 1,
$$
which implies that $\theta^T s_{\mu_n} \in (R_1-\epsilon, R_1+\epsilon)$ and $\|s_{\mu_n}\| \in (R_2-\epsilon,R_2+\epsilon)$. By the choice of $\epsilon$, it follows that $s_{\mu_n} = s_{\mu}$ for $n$ large enough.

Now, suppose that $A$ is a set such that $\partial A = 0$. We must prove that $(\mu_n - s_{\mu_n})(A) \to (\mu - s_{\mu})(A)$. This is equivalent to $$
\mu_n(A + s_{\mu_n}) \to \mu(A+s_{\mu}).
$$
For any $\epsilon>0$, consider the open set $$A^{-\epsilon} = \{x\in A+s_{\mu} : \|x-y\|< \epsilon \text{ for all } y \in \partial (A+s_{\mu})\}$$
and the compact set
$$
A^{+\epsilon} = \{x\in A+s_{\mu} : \|x-y\| \leq \epsilon \text{ for some } y \in A+s_{\mu}\}
$$
Then, $A^{-\epsilon} \subset A+s_{\mu} \subset A^{+\epsilon}$, and since $s_{\mu_n} \to s_{\mu}$, we have that $\1_{A+s_{\mu_n}} \to \1_{A+ s_{\mu}}$ pointwise as functions on $\R^d$. In particular, this implies that for some $N \in \mathbb N$, $A^{-\epsilon}\subset A + s_{\mu_n} \subset A^{+\epsilon}$ for all $n>N$. Thus, $$
\mu_n(A^{-\epsilon}) \leq \mu_n(A+s_{\mu_n}) \leq \mu_n(A^{+\epsilon})
$$
for all $n>N$.

By \cite[Proposition 3.12]{Resnick1987}, $\limsup_{n} \mu_n(A^{+\epsilon}) \leq \mu(A^{+\epsilon})$ and $\liminf_{n} \mu_n(A^{-\epsilon}) \geq \mu(A^{-\epsilon})$. Thus, we have \begin{multline*}
\mu(A^{-\epsilon}) \leq \liminf_{n} \mu_n(A^{-\epsilon}) \leq \liminf_{n} \mu_n(A+s_{\mu_n}) \leq \limsup_{n} \mu_n(A+s_{\mu_n}) \\ \leq \limsup_{n} \mu_n(A^{+\epsilon}) \leq \mu(A^{+\epsilon}).
\end{multline*}
This is true for any $\epsilon>0$. The result now follows by letting $\epsilon \to 0$.

Having shown that $\mu \to \mu - s_{\mu}$ is continuous on the space of simple point processes, we conclude by the continuous mapping theorem (see \cite[Theorem 5.1]{BIL}) since $\eta$ is a simple point process a.s.
\end{proof}

We shall now proceed to prove the theorems in each of the remaining sections in their order of appearance. It will be apparent that the results in the forthcoming sections will utilize the preceding results heavily.

\section{Proofs for Section~\ref{1DWeibull}}\label{appendix:d}

In this section, we will prove all the results in Section~\ref{1DWeibull}. Let $X$ be a random variable in the Weibull regime (see Assumption~\ref{weibull}), with maximum value $\mx>0$ and tail index $\alpha>0$. 

\subsection{Proof of Lemma~\ref{KStoSL} and Weibull regime examples}

We begin by proving Lemma~\ref{KStoSL} and the examples of random variables in the Weibull regime.

\begin{proof}[Proof of Lemma~\ref{KStoSL}]

Recalling that $b_n>0$, note that for any $n \geq 1$, \begin{align}& \sup_{x \in \R} |\Prob(X_{1,n} \leq x) - \Prob(X_{2,n} \leq x)|\nonumber \\ =& \sup_{x\in \R} \left|\Prob\left(\frac{X_{1,n}-a_n}{b_n} \leq \frac{x-a_n}{b_n}\right)-\Prob\left(\frac{X_{2,n}-a_n}{b_n} \leq \frac{x-a_n}{b_n}\right)\right| \nonumber\\
& = \sup_{y \in \R} \left|\Prob\left(\frac{X_{1,n}-a_n}{b_n} \leq y\right)-\Prob\left(\frac{X_{2,n}-a_n}{b_n} \leq y\right)\right|\label{touch} \end{align}
since the map $x \mapsto y=\frac{x-a_n}{b_n}$ is a bijective mapping on $\R$. Now recalling the definition \eqref{KS} of $\KS$ we have \begin{align*}
\KS(X_{1,n},X_{2,n}) =& \sup_{x} \left|\Prob\left(\frac{X_{1,n}-a_n}{b_n} \leq x\right)-\Prob\left(\frac{X_{2,n}-a_n}{b_n} \leq x\right)\right| \\ =& \KS\left(\frac{X_{1,n}-a_n}{b_n}, \frac{X_{2,n}-a_n}{b_n}\right)
\end{align*}
For part (a), we may apply Lemma~\ref{assist} since $Z_1=Z_2=Z$ is continuous. As a consequence, we obtain that $\textsf{KS}\left(\frac{X_{i,n}-a_n}{b_n}, Z\right) \to 0$. The result now follows from the above equality and the triangle inequality for the $\textsf{KS}$ distance.

For part (b), since $Z_1 \neq Z_2$ in distribution, by Lemma~\ref{assist} there is a point $c$ at which both the CDF of $Z_1$ and the CDF of $Z_2$ are continuous, such that $F_{Z_1}(c) \neq F_{Z_2}(c)$. But $\Prob\left(\frac{X_{i,n}-a_n}{b_n} \leq c\right) \to F_{Z_i}(c)$ by convergence in distribution. That the $\KS$ distance between the scaled random variables doesn't go to zero, now follows directly from this observation and \eqref{touch}.
\end{proof}

Next, we will sketch the proof of the list of examples in Lemma~\ref{examples1DWeibull}.

\begin{proof}[Proof of Lemma~\ref{examples1DWeibull}]

For (a), let $X = Beta(a,b)$ for some $a,b>0$. Then, $X$ has the density $cx^{a-1}(1-x)^{\mx-1}$ for $x \in [0,1]$ and a normalizing constant $c$. Clearly its maximum value is $\mx=1$, and by the chain rule, the function $1-F_X(1-x)$ has the derivative $c(1-x)^{a-1}x^{\mx-1}$. 

Note that $c(1-x)^{a-1}$ is a function which is continuous and non-zero at $0$, hence regularly varying with index $0$ at $0$. On the other hand $x^{\mx-1}$ is regularly varying of order $\mx-1$. By the product rule \cite[Proposition 1.5.7(iv)]{BinghamGoldieTeugels1987} it is readily seen that $c(1-x)^{a-1}x^{\mx-1}$ is regularly varying of order $\mx-1$ at $0$.

By \cite[Proposition 1.5.8]{BinghamGoldieTeugels1987}, its integral i.e. $1-F_X(1-x)$ is regularly varying at $0$ of index $b>0$. Hence, it lies in the Weibull regime with parameter $b$, completing the proof of the first example.

For (b), suppose that $X$ is either a normal or exponential random variable, and let $M>0$. It is easily checked that if $f$ is the density of the truncated random variable $X_M = X\1_{X \leq M}$, then $\lim_{x \to M^-} f(x) = f(M)$ is a positive number. Therefore, it is clear that $f$ is regularly varying at $M$ or order $0$. By an argument similar to the one made for $Beta(a,b)$, it follows that $X$ is in the Weibull regime with $1-F_X(M-x)$ regularly varying at $0$ with index $1$.

We shall now address (c) in detail. Let $X,Y$ be in the Weibull regime with parameters $\alpha,\beta$ respectively and maximum values $\mx_X,\mx_Y$. Then, $X+Y$ has maximum value $\mx_X+\mx_Y$. Let $F_X$ be the CDF of $X$ and $F_Y$ the CDF of $Y$. Note that for any $u>0$, \begin{align}\mathbb P(X+Y> (\mx_X+\mx_Y-u)) =& \int_{\mx_X-u}^{\mx_X} (1-F_Y (\mx_X+\mx_Y-u-x)) dF_X(x) \nonumber\\
=& \int_{0}^{u} (1-F_Y(\mx_Y-u+x)) d(1-F_X(\mx_X - x))) \nonumber \\
=& \int_{0}^{u} (1-F_X(\mx_X - x))d(1-F_Y(\mx_Y-u+x)) \label{one}\end{align}
where we used integration by parts and noted that $\lim_{x \to 0^+} 1-F_X(\mx_X-x) = 0$ and $\lim_{x \to u^-} 1-F_Y(\mx_Y-u+x) = 0$. 

Similarly, for any $t,u>0$,
\begin{align}
\mathbb P(X+Y> (\mx_X+\mx_Y-tu)) = & \int_{0}^{tu} (1-F_Y(\mx_Y-tu+x)) d(1-F_X(\mx_X - x)) \nonumber\\= & t\int_{0}^{u} (1-F_Y(\mx_Y-tu+tx)) d(1-F_X(\mx_X - tx)) \nonumber \\
=& t\int_0^u \mathbb (1-F_X(\mx-tx)) d(1-F_Y(\mx_Y-tu+tx)),\label{two}
\end{align}
where we used integration by parts and noted that $\lim_{x \to 0^+} 1-F_X(\mx_X-tx) = 0$ and $\lim_{x \to u^-} 1-F_X(\mx_X-tu+tx) = 0$. 

The representations so far clearly hint at the appearance of regularly varying terms in the integrand and integrator of \eqref{one} and \eqref{two}. Thus, we expect the ratio of these two terms to behave regularly, which is what will be proved now. We shall now fix $t>1$ and study the ratio of the left hand sides of \eqref{one} and \eqref{two} as $u \to 0$. An analogous argument will apply for $t<1$.

Let $A>1,\delta>0$ be arbitrary. By the definition in \eqref{weibull} of the Weibull regime and Potter's theorem \cite[Theorem 1.5.6(iii)]{BinghamGoldieTeugels1987}, there exists $X>0$ such that for $p,q \leq X$, \begin{gather}
A^{-1} q^{\alpha-\delta} \leq \frac{(1-F_Y (\mx_Y-qp))}{(1-F_Y(\mx_Y-p))} \leq A q^{\alpha+\delta} \label{three}\\
A^{-1} q^{\beta-\delta} \leq \frac{(1-F_X(\mx_X-qp))}{(1-F_X(\mx_X-p))} \leq A q^{\beta+\delta}\label{four}
\end{gather}

Now, applying \eqref{three} to every $u < X/t$ and integrating these inequalities from $0$ to $u$ with respect to $1-F_X(\mx_X-x)$, 
\begin{equation}\label{five}
A^{-1}t^{\alpha-\delta}\leq \frac{\int_{0}^u (1-F_Y(\mx_Y-tu+tx))d(1-F_X(\mx-x))}{\int_{0}^u (1-F_Y(\mx_Y-u+x))d(1-F_X(\mx-x))} \leq At^{\alpha+\delta}
\end{equation}
Similarly, applying \eqref{four} to every $u<X/t$ and  integrating these inequalities from $0$ to $u$ with respect to $1-F_Y(\mx_Y-tu+tx)$, 
\begin{equation}\label{six}
A^{-1}t^{\beta-\delta}\leq \frac{\int_{0}^u (1-F_X(\mx-tx))d(1-F_Y(\mx_Y-tu+tx))}{\int_{0}^u (1-F_X(\mx-x))d(1-F_Y(\mx_Y-tu+tx))} \leq At^{\beta+\delta}
\end{equation}

Multiplying \eqref{five} and \eqref{six}, while also noting the equalities in \eqref{one} and \eqref{two}, $$
A^{-2} t^{\beta + \alpha-2\delta}\leq \frac{\mathbb P(X+Y> (\mx_X+\mx_Y-tu))}{\mathbb P(X+Y> (\mx_X+\mx_Y-u))} \leq A^2 t^{\beta+\alpha+2\delta}
$$
for all $u,t>0$ such that $tu<X$. Now, fixing $t$ and letting $u \to 0$ above gives that \begin{align*}
A^{-2} t^{\beta + \alpha-2\delta} \leq & \liminf_{u \to 0} \frac{\mathbb P(X+Y> (\mx_X+\mx_Y-tu))}{\mathbb P(X+Y> (\mx_X+\mx_Y-u))}\\ \leq & \limsup_{u \to 0}\frac{\mathbb P(X+Y> (\mx_X+\mx_Y-tu))}{\mathbb P(X+Y> (\mx_X+\mx_Y-u))} \leq  A^2 t^{\beta+\alpha+2\delta}.
\end{align*}
Since this is true for all $\delta>0$ and $A>1$, letting $A \downarrow 1$ and $\delta \downarrow 0$ gives $$
\lim_{u \to 0} \frac{\mathbb P(X+Y> (\mx_X+\mx_Y-tu))}{\mathbb P(X+Y> (\mx_X+\mx_Y-u))} = t^{\alpha+\beta}.
$$
The same argument as here can be repeated with $t<1$, showing that the above statement holds for any $t>0$ (for $t=1$ it is immediate).  Consequently, the result has been proved.
\end{proof}

We will next prove Lemma~\ref{weibtoweib1D} i.e. that the site-specific tilt $g$ retains the Weibull assumption.

\begin{proof}[Proof of Lemma~\ref{weibtoweib1D}]
Since $g$ is strictly increasing, the maximum value of $g(X)$ is $g(\mx)$. Let $F_X$ be the CDF of $X$, and $F_{g(X)}$ be the CDF of $g(X)$. Consider the quantity \begin{align*}
1 - F_{g(X)}(g(\mx)-x) &= \Prob(g(X) > g(\mx)-x) \\ &= \Prob(X > g^{-1}(g(\mx)-x))\\& = 1-F_X(g^{-1}(g(\mx)-x)).
\end{align*}
Here, $g^{-1}$ is the set theoretic inverse of $g$. Note that since $g$ is strictly increasing, $g^{-1}(y) < x$ if and only if $y< g(x)$. 

We must prove that $1 - F_{g(X)}(g(\mx)-x)$ is regularly varying at $0$. However, we have written this function as a composition of two functions $$
 1-F_{g(X)}(g(\mx)-x) = h(j(x))
$$
where $h(z) = 1 - F_X(\mx-z)$ and $j(y) = \mx-g^{-1}(g(\mx)-y)$.

Observe that $j(y) \to 0$ as $y \to 0$ by the continuity of $f$. Furthermore, $h$ is regularly varying at $0$ with index $\alpha$. If we can show that $j$ is regularly varying at $0$ with index $\frac 1{\beta}$, then by \cite[Proposition 1.5.7]{BinghamGoldieTeugels1987}, $x \mapsto 1-F_{g(X)}(g(\mx)-x)$ is regularly varying at $0$ of order $\frac{\alpha}{\beta}$, completing the proof. 

However, observe that the definition of $j(y)$ implies that $g(\mx) - g(\mx-j(y)) = y$. In other words, $j$ is the set theoretic inverse of the function $x \mapsto g(\mx) - g(\mx-x)$. By \cite[Theorem 1.5.12]{BinghamGoldieTeugels1987}, $j$ is regularly varying at $0$ of order $\frac 1{\beta}$.
\end{proof}

Thus, general twists are absorbed into our setup. We will therefore continue to assume that $X$ is in the Weibull regime with index $\alpha>0$ and maximum value $B$. Let $\theta_n, n \geq 1$ be a sequence, and recall $M_{\theta_n}$ from \eqref{m2tmt2}. Also recall $X_{\theta_n}$ from \eqref{twisted1D}.

We will now prove Lemma~\ref{lem:nu}.

\begin{proof}[Proof of Lemma~\ref{lem:nu}]
For part (a), suppose that $A \subset \mathbb R^d$ is compact and satisfies $\nu(A) = 0$. Then, by Assumption~\ref{ass:mvrv}, \begin{align}
\nu(cA) &= \lim_{t \to 0} \frac{1}{U(t)} \mathbb P\left(\frac{x_{\theta}-X}{t}\in cA\right)\nonumber\\
&= \lim_{t \to 0} \frac{1}{U(t)} \mathbb P\left(\frac{x_{\theta}-X}{ct}\in A\right)\nonumber\\
&= \lim_{t \to 0} \frac{c^{\alpha}}{U(ct)} \mathbb P\left(\frac{x_{\theta}-X}{ct}\in A\right) \nonumber\\
& = c^{\alpha} \nu(A).\label{scaling}
\end{align}

To prove part (b), let $\nu\{0 \leq \theta^Ty \leq 1\} =C\in (0,\infty)$, and define the measure $\mu$ on $\mathcal{S} = \{\theta^T y = 1\}$ by $$
\mu(B) =\nu\left\{0 < \theta^T y < 1 , \frac{y}{\theta^T y} \in B\right\}.
$$
Observe that $\mu(\mathcal{S}) = C$. Now, observe that for any $c,x>0$ and $B \subset \mathcal{S}$, by part (a) we have
\begin{equation}\label{deter}
\nu\left(\left\{\theta^T y \in (0,c], \frac{y}{\theta^T y} \in B\right\}\right) = c^{\alpha} \left(\left\{\theta^T y \in (0,1], \frac{y}{\theta^T y} \in B\right\}\right) = c^{\alpha}\mu(B).
\end{equation}
Since sets of the form $$\left\{\theta^T y \in (0,c], \frac{y}{\theta^T y} \in B\right\}$$ generate $\{\theta^T y > 0\}$, it follows that $\nu$ is uniquely determined by the scaling property, and if $a,b>0$ and $B \subset \mathcal{S}$ are arbitrary then \begin{align*}
&\nu\left(\left\{\theta^T y \in (a,b], \frac{y}{\theta^T y} \in B\right\}\right) \\ =& \nu\left(\left\{\theta^T y \in (0,b], \frac{y}{\theta^T y} \in B\right\}\right) - \nu\left(\left\{\theta^T y \in (0,c], \frac{y}{\theta^T y} \in B\right\}\right)\\ =& (b^{\alpha} - a^{\alpha}) \mu(B)
\end{align*}   
by \eqref{deter}.

To prove part(c), observe that $$
\int_{\theta^T y >0} e^{-\theta^T y } d \nu(y) < \sum_{n=1}^{\infty} e^{-n} \nu\left(\left\{n-1 < \theta^T y \leq n\right\}\right) \leq C\sum_{n=1}^{\infty} e^{-n} (n^{\alpha} - (n-1)^{\alpha}) < \infty,
$$
where we used part (b).
\end{proof}

We will now dedicate a section to the asymptotic behavior of $M_{\theta_n}$ and the proof of Theorem~\ref{thm:m2tmt21D}.

\subsection{Proof of Theorem~\ref{thm:m2tmt21D}}

We first prove Lemma~\ref{lem:asy} about the growth of $\E[e^{\theta_n X}]$ as $\theta_n \to \infty$. Throughout this section let $F_X$ denote the CDF of $X$.

\begin{proof}[Proof of Lemma~\ref{lem:asy}]
Write \begin{align}
\E[e^{\theta_n X}] &=\int_{-\infty}^\mx e^{\theta_n y} dF_X(y)  \nonumber\\
& = e^{\theta_n \mx}\int_{-\infty}^\mx e^{\theta_n (y-\mx)} dF_X(y) \nonumber\\
& = -e^{\theta_n \mx} \int_{0}^{\infty} e^{-\theta_n u} dF_X(\mx-u)\nonumber\\
& = e^{-\theta_n \mx} \int_0^{\infty} e^{-\theta_n u} d G_X(u),\label{help1}
\end{align}
 where $G_X(u) = 1-F_X(\mx-u)$. Now, $G_X(u)$ is regularly varying at $0$ of order $\alpha>0$. By the Karamata Tauberian theorem (see \cite[Theorem 3, Section XIII.5]{feller1971}) it follows that $$
\frac{\int_{0}^{\infty} e^{-\theta_n u} dG_X(u)}{\Gamma(\alpha+1)\left(1-F_X\left(\mx-\frac{1}{\theta_n}\right)\right)} \to 1.
$$
Our result directly follows from the above equation, \eqref{help1} and algebraic rearrangements.
\end{proof}

The proof of Theorem~\ref{thm:m2tmt21D} is a corollary of this.

\begin{proof}[Proof of Theorem~\ref{thm:m2tmt21D}]
If $\mx$ is the maximum value of $X$, then we have by the previous lemma that \begin{gather*}
\frac{\mathbb E[e^{2 \theta_n X}]}{e^{2 \theta_n \mx} \left(1-F_X\left(\mx-\frac{1}{2\theta_n}\right)\right)} \to \Gamma(1+\alpha) \\
\frac{e^{2 \theta_n \mx}\left(1-F_X\left(\mx-\frac{1}{\theta_n}\right)\right)^2}{\mathbb E[e^{\theta_n X}]^2} \to \frac 1{\Gamma(1+\alpha)^2} 
\end{gather*}
Multiplying these together and noting that $\frac{\left(1-F_X\left(\mx-\frac{1}{\theta_n}\right)\right)}{\left(1-F_X\left(\mx-\frac{1}{2\theta_n}\right)\right)} \to 2^{\alpha}$ by Assumption~\ref{weibull}, we obtain the result.
\end{proof}

In the next section, we will prove the scaling limit of $X_{\theta}$ i.e. Theorem~\ref{thm:sltrue1D}.

\subsection{Proof of Theorem~\ref{thm:sltrue1D}}

The following lemma covers the asymptotics of event probabilities under large twists. We will use it to prove Theorem~\ref{thm:sltrue1D}.

\begin{lemma}\label{help}
Let $\theta_n$ be any sequence converging to infinity, and $r_n$ be any sequence increasing to a real number $\mx$. Suppose $\theta_{n}(\mx-r_n) \to C \in (0,\infty)$. Then, for any random variable $X$ belonging to the Weibull regime with parameter $\alpha>0$ and maximum value $\mx$, we have $$\frac{\mathbb E[e^{\theta_n X}\1_{X > r_n}]}{\mathbb E[e^{\theta_n X}]} \to \frac{\int_{0}^C t^{\alpha-1}e^{-t}dt}{\Gamma(\alpha)} \in (0,1).$$
\end{lemma}
\begin{proof}
As usual, let $F_X$ be the CDF of $X$. We begin by noting that $$
\E[e^{\theta_n X}\1_{X > r_n}] = \int_{r_n}^{\mx} e^{\theta_n y} dF_X(y).
$$
We perform the following manipulations, taking $G_X(u) = 1-F_X(\mx-u)$. Since this is regularly varying at $0$ by Assumption~\ref{weibull}, by \cite[Theorem 1.4.1]{BinghamGoldieTeugels1987} we have $G_X(u) = L(u)u^{\alpha}$ for some slowly varying function $L$ (i.e. regularly varying with index $0$) at zero. 

\begin{align}
&\int_{r_n}^{\mx} e^{\theta_n y} dF_X(y) \nonumber \\ = &e^{\theta_n \mx} \int_{r_n}^\mx e^{\theta_n(y-\mx)}dF_X(y)\nonumber\\
=& -e^{\theta_n\mx}\int_0^{\mx-r_n} e^{-\theta_n u} dF_X(\mx-u)\nonumber\\
=&  e^{\theta_n\mx}\int_0^{\mx-r_n} e^{-\theta_n u} dG_X(u) \nonumber\\
\overset{IBP}{=}& e^{\theta_n \mx}\left[[e^{-\theta_n u}G_X(u)]^{\mx-r_n}_0 + \int_{0}^{\mx-r_n} G_X(u)\theta_ne^{-\theta_nu} du\right]  \nonumber \\
=&e^{\theta_n \mx}\left[e^{-\theta_n(\mx-r_n)}G_X(\mx-r_n) + \int_{0}^{\theta_n(\mx-r_n)} L\left(\frac{t}{\theta_n}\right)\left(\frac{t}{\theta_n}\right)^{\alpha}e^{-t} dt\right] \nonumber \\
= & e^{\theta_n \mx}\theta_n^{-\alpha}L\left(\frac 1{\theta_n}\right)e^{-\theta_n(\mx-r_n)}L(\mx-r_n)(\mx-r_n)^{\alpha}L^{-1}\left(\frac{1}{\theta_n}\right)\theta_n^{\alpha}\nonumber \\ +&e^{\theta_n \mx}\theta_n^{-\alpha}L\left(\frac 1{\theta_n}\right)\int_{0}^{\theta_n(\mx-r_n)} \frac{L(\frac{t}{\theta_n})}{L\left(\frac {1}{\theta_n}\right)}t^{\alpha}e^{-t} dt\label{eqn}
\end{align}
We used the fact that $G_X(\mx) = 0$ above. We will now consider the limit of both quantities inside the bracket on the last line.

Clearly, the first term is easy to handle : \begin{equation}\label{firstterm}e^{-\theta_n(\mx-r_n)}L(\mx-r_n)(\mx-r_n)^{\alpha}L^{-1}\left(\frac{1}{\theta_n}\right)\theta_n^{\alpha} \to e^{-C}C^{\alpha},\end{equation} since $\theta_n(\mx-r_n) \to C$ and $L$ is slowly varying. 

For the second term in \eqref{eqn}, note that $\theta_n(\mx-r_n) \to C$. Thus, it is bounded, and we assume that $\theta_n(\mx-r_n)<M$ for all $N>0$. Pointwise, on $[0,M]$ we have $\frac{L(\frac{t}{\theta_n})}{L\left(\frac {1}{\theta_n}\right)} \to 1$. Thus, we have $$
\frac{L(\frac{t}{\theta_n})}{L\left(\frac {1}{\theta_n}\right)}t^{\alpha}e^{-t}\1_{[0,\theta_n(\mx-r_n)]} \to t^{\alpha}e^{-t}\1_{[0,C]} 
$$ 
as $n$ converges to infinity, pointwise on $[0,M]$. To apply the dominated convergence theorem, we note that $L$ is slowly varying at $0$, and therefore $L(s)s^{\alpha/4} \to 0$ as $s \to 0$, and $L(s)s^{-\alpha/4} \to +\infty$ as $s \to 0$ (\cite[Proposition 1.5.1]{BinghamGoldieTeugels1987}). Therefore, there is a constant K such that $Ks^{\alpha/4} \leq L(s) \leq K s^{-\alpha/4}$ for $s \in [0,M]$, which implies that $$
\frac{L(\frac{t}{\theta_n})}{L\left(\frac {1}{\theta_n}\right)}t^{\alpha}e^{-t}\1_{[0,\theta_n(\mx-r_n)]} \leq K^2 t^{3\alpha/4}e^{-t}\1_{[0,M]}
$$
which is integrable. By the dominated convergence theorem,
$$
\int_{0}^{\theta_n(\mx-r_n)} \frac{L(\frac{t}{\theta_n})}{L\left(\frac {1}{\theta_n}\right)}t^{\alpha}e^{-t} dt \to \int_{0}^{C}t^{\alpha}e^{-t} dt.
$$
Combining the above with \eqref{firstterm}, \eqref{eqn} and Lemma~\ref{lem:asy},
$$
\frac{\E[e^{\theta_n X}\1_{X > r_n}]}{\E[e^{\theta_n X}]} \to \frac{e^{-C}C^{\alpha} + \int_0^C t^{\alpha}e^{-t}dt}{\Gamma(\alpha+1)} = \frac{\int_0^{C} t^{\alpha-1}e^{-t} dt}{\Gamma(\alpha)},
$$
where the numerator simplifies by integration-by-parts, and the denominator satisfies $\Gamma(\alpha+1) = \alpha \Gamma(\alpha)$. This completes the proof.
\end{proof}

With this, one can complete the proof of Theorem~\ref{thm:sltrue1D}.

\begin{proof}[Proof of Theorem~\ref{thm:sltrue1D}]

Let $C>0$ be arbitrary. For a sequence $\theta_n$ increasing to infinity, let $s_n = \frac{C}{\theta_n}$. We have following some algebraic rearrangement, 
$$
 \{\theta_n(\mx - X_{\theta_n}) \leq C\}=\left\{X_{\theta_n} \geq \mx - \frac{C}{\theta_n}\right\}= \{X_{\theta_n} \geq \mx-s_n\}.
$$
Now, by the definition of tilting,
$$
\mathbb P(X_{\theta_n} \geq \mx-s_n) = \frac{\E[e^{\theta_n X}\1_{X \geq \mx - s_n}]}{\E[e^{\theta_n X}]} =  \frac{\E[e^{\theta_n X}\1_{X \geq \mx - s_n}]}{\E[e^{\theta_n X}]}.
$$
 By Lemma~\ref{help} applied with $X$ and $r_n = \mx - s_n$, 
 $$
\mathbb P(X_{\theta_n} \geq \mx - s_n)  \to  \frac{\int_0^C t^{\alpha-1}e^{-t}dt}{\Gamma\left(\alpha\right)} = \mathbb P\left(\Gamma\left(\alpha,1\right)\leq C\right).
$$
It follows that $$ \Prob(\theta_n(\mx -X_{\theta_n}) \leq C) \to  \mathbb P\left(\Gamma\left(\alpha,1\right)\leq C\right)$$
as $n \to \infty$. This is sufficient to show the theorem, since we only need to consider $C \in (0,\infty)$ for non-negative random variables to converge in distribution.
\end{proof}

Finally, recall the estimator random variable $R_{n,\theta}$ defined by $$
\mathbb P(R_{n,\theta} \in A) = \frac{\sum_{i=1}^n e^{\theta f(X_{i})}\1_{X_{i} \in A}}{\sum_{i=1}^n e^{\theta f(X_{i})}}.
$$
We are now ready to start proving Theorem~\ref{thm:slemp1D}. We will prove its parts in the order in which they were stated.

\subsection{Proof of Theorem~\ref{thm:slemp1D}}

We begin with the proof of part (a). Note that the rough idea was already discussed above the statement of the theorem.

\begin{proof}[Proof of Theorem~\ref{thm:slemp1D}(a)]
Let $a_n,b_n$ be such that $\frac{X-a_n}{b_n} \to \Gamma(\alpha,1)$ as in Theorem~\ref{thm:sltrue1D}. Note that for any $x \in \R$, that $B_n = \{y \leq a_n+b_nx\}$ satisfies \eqref{hyp} follows from Lemma~\ref{help}. Furthermore, $\frac{M_{\theta_n}}{n} \to 0$ by assumption. Thus, by Proposition~\ref{posassist} and Theorem~\ref{thm:sltrue1D}, $X_n$ and $R_{n,\theta_n}$ have the same scaling limit, completing the proof.
\end{proof}

Next, we prove Theorem~\ref{thm:slemp1D}(b).

In this section, we explicitly find the limiting random variable in Theorem~\ref{thm:slemp1D}(b) using Lemma~\ref{dekhliyo}.

\begin{proof}[Proof of Theorem~\ref{thm:slemp1D}](b)
Suppose that $\frac{M_{\theta_n}}{n} \to C\in (0,\infty)$. Then, by Theorem~\ref{thm:m2tmt21D} we know that $\theta_n(\mathcal{M} - F_X^{-1}(1-\frac 1n)) \to C_1$ for some $C_1>0$, and hence $\theta_n \to \infty$. Let $C_2 \in (0,\infty)$ and $A_n = [\mx- C_2(\mx-F_X^{-1}(1-\frac 1n)), \mx]$.

Consider $\mathbb P(R_{n,\theta} \in A_n)$, and divide the top and bottom by $e^{\theta \mx}$ to obtain 
\begin{align*}
\mathbb P(R_{n,\theta} \in A_n) =& \frac{\sum_{i=1}^n e^{\theta X_{i}}\1_{X_{i} \in A_n}}{\sum_{i=1}^n e^{\theta X_{i}}}\\ =&  \frac{\sum_{i=1}^n e^{\theta (X_{i}-\mx)}\1_{X_{i} \in A_n}}{\sum_{i=1}^n e^{\theta (X_{i}-\mx)}} \\
=& \Phi_n\left(\sum_{i=1}^n e^{\theta_n(\mx-X_i)}1_{A_n}, \sum_{i=1}^n e^{\theta_n(\mx-X_i)}\right),
\end{align*}
where 
$$
\Phi(y,z) = \frac{y}{z}.
$$

At this point, we verify the hypotheses of Lemma~\ref{dekhliyo}. We take $\theta = 1$, $C_n = \theta_n$. Then, $X_{\theta_n} \to \mx$ in probability. Let $$x_n = \mathcal{M}, a_n = \mx - F_X^{-1}\left(1-\frac 1n\right) , C_1 = C_1,  \nu(dy) = \alpha y^{\alpha-1} dy \text{ and } D = [0,C_2].$$
Hypotheses (a) and (c) of the lemma are easily verified, while hypothesis (b) follows from the Weibull regime being an instance of multivariate regular variation (see the examples in Section~\ref{mdweibull}).

Since $\Phi$ is a continuous mapping, and $\int e^{-C_1y} dPRM(\nu) \neq 0$ with probability $1$ by the definition of a PRM, by Lemma~\ref{dekhliyo} it would follow that \begin{align*}
\mathbb P(R_{n,\theta} \in A_n) = & \Phi_n\left(\sum_{i=1}^n e^{-\theta_n(\mx-X_i)}1_{A_n}, \sum_{i=1}^n e^{-\theta_n(\mx-X_i)}\right) \\ \overset{d}{\to} & \frac{\int e^{-C_1y}1_{y \leq C_2} d PRM(\nu)}{\int e^{-C_1y} dPRM(\nu)} := \mathbb P(Z \leq C_2)
\end{align*}

for some random variable $Z$. Finally, note that \begin{equation*}\mathbb P(R_{n,\theta} \in A_n) = \mathbb P\left(\theta_n(\mx-R_{n,\theta}) \leq C_2 \theta_n\left(\mx-F_X^{-1}\left(1-\frac 1n\right)\right)\right) \to \mathbb P(Z \leq C_2)
\end{equation*}
 Combining the two statements above, 
 $$
\mathbb P(\theta_n(\mx-R_{n,\theta}) \leq C_1C_2) \to \mathbb P(Z\leq C_2) = \mathbb P(C_1Z \leq C_1C_2).
$$
Thus, the random variable $Z_{C_1,PRM} = C_1Z$ is the desired limit. Note that this random variable depends upon the Poisson random measure, while the other scaling limit $\Gamma(\alpha,1)$ does not. It follows that these two random variables are not the same, which concludes the proof.

We remark that the limiting random variable $Z_{PRM}$ is continuous, but do not prove this here.
\end{proof}

Now we will prove the last part i.e. Theorem~\ref{thm:slemp1D}(c).

\begin{proof}[Proof of Theorem~\ref{thm:slemp1D}(c)]

We claim that in probability, $$
\frac{\sum_{i=1}^n e^{\theta_n X_{i}}}{e^{\theta_n \max_i X_{i}}} \to 1
$$
if $\theta_n(B - F_X^{-1}(1-\frac 1n)) \to +\infty$ (which, by Theorem~\ref{thm:m2tmt21D} is implied by $\frac{M_{\theta_n}}{n} \to +\infty$). In particular, this implies that 
\begin{equation}
\frac{\sum_{i=1}^n e^{\theta_n X_{i}}\1_{X_{i} \neq \max_j X_j}}{e^{\theta_n \max_i X_{i}}} \to 0.
\end{equation}

 By Lemma~\ref{pppc}, $\sum_{i=1}^n \delta_{\frac{\mx-X_{i}}{\mx-F^{-1}(1-\frac 1n)}} \Rightarrow PRM( \nu)$ where $d\nu(y) = \alpha y^{\alpha-1}$. Now, note that the smallest point of the former point process is attained when $X_{i} = \max_j X_j$, which without loss of generality we assume is attained at $X_1$. Therefore, we get by Lemma~\ref{pppconv} (with $\theta = 1$) that 
$$
\sum_{i=1}^n \delta_{\frac{X_1-X_{i}}{\mx-F^{-1}(1-\frac 1n)}} \Rightarrow PRM( \nu) - s_{PRM( \nu)}.
$$

Let $M>0$ be an arbitrary but fixed parameter. Let $f(x) = e^{-Mx}\1_{x <M}$. Note that $f$ is continuous and compactly supported on $[0,\infty)$. Thus, by the continuous mapping theorem (see \cite[Section 3.5]{Resnick1987}) we have \begin{equation}\label{eq:conv}
\sum_{i=1}^n e^{-M\frac{X_1-X_{i}}{\mx-F_X^{-1}(1-\frac 1n)}} \xrightarrow{d} \int e^{-Mx}\1_{x<M}\left[PRM( \nu) - s_{PRM( \nu)}\right].
\end{equation}

Now, it is clear that for any $\epsilon>0$,
\begin{align*}
&\lim_{M \to \infty} \mathbb P\left[\sum_{i=1}^n e^{-\theta_n\left(X_1-X_i\right)}\1_{X_1-X_i < (\mx-F_X^{-1}(1-\frac 1n)M)} > 1+\epsilon\right] \\ = &\mathbb P\left[\sum_{i=1}^n e^{-\theta_n(X_i-X_1)} > 1+\epsilon\right].
\end{align*}

For any $\epsilon>0$ and $M>0$, employing \eqref{eq:conv} gives \begin{align}
&\limsup_{n} \Prob\left[\sum_{i=1}^n e^{-\theta_n\left(\mx-F_X^{-1}(1-\frac 1n)\right)\left(\frac{X_1-X_i}{B - F_X^{-1}(1-\frac 1n)}\right)}\1_{X_1-X_i < (\mx-F_X^{-1}(1-\frac 1n)M)} > 1+\epsilon\right] \nonumber \\
\leq & \limsup_{n \to \infty}\Prob\left[\sum_{i=1}^n e^{-M\left(\frac{X_1-X_i}{B - F_X^{-1}(1-\frac 1n)}\right)}\1_{X_1-X_i < (\mx-F_X^{-1}(1-\frac 1n)M)} > 1+\epsilon\right] \nonumber\\ = & \Prob\left[\int e^{-Mx}\1_{x<M} (PRM(\nu) - s_{PRM(\nu)}) > 1+\epsilon\right]. \label{Seven}
\end{align}

Note that $e^{-Mx}\1_{x<M}$ is dominated by $e^{-x}$ for $M>1$. Furthermore, as $M \to \infty$, we have $e^{-Mx}\1_{x<M} \to 1_{\{0\}}$ pointwise. In order to use the dominated convergence theorem, we need to prove that $e^{-x}$ is integrable under $PRM(\nu) - s_{PRM(\nu)}$ a.s. : we will prove the stronger statement that $e^{-x}$ is integrable under $PRM(\nu)$ a.s.

In order to prove this, observe that $(1-e^{-x}) \leq x$ for all $x>0$. Thus, 
\begin{equation*}
\int_{0}^{\infty} (1-e^{-e^{-x}})d\nu(x) \leq \int_{0}^{\infty} e^{-x} d\nu(x) = \alpha\int_{0}^{\infty} e^{-x}x^{\alpha-1}dx<\infty.
\end{equation*}
The left hand side is equal to $-\ln \E[e^{-\int_{0}^{\infty} e^{-x} dPRM(\nu)}]$ by \cite[Proposition 3.6(ii)]{Resnick1987}. It follows that $\int_{0}^{\infty} e^{-x} dPRM(\nu)$ is finite a.s.

By the dominated convergence theorem, $\int e^{-Mx}\1_{x<M} (PRM(\nu) - s_{PRM(\nu)}) \to 1$ a.s., which implies that the probability in \eqref{Seven} converges to $0$. In particular, it follows that 

$$
\limsup_{n} \Prob\left[\sum_{i=1}^n e^{-\theta_n\left(\mx-F_X^{-1}(1-\frac 1n)\right)\left(\frac{X_1-X_i}{B - F^{-1}(1-\frac 1n)}\right)} > 1+\epsilon\right] = 0,
$$
which instantly proves the claim, since the sum must always exceed $1$ (by taking the term $i=1$), therefore it was sufficient to check only exceedance by $\epsilon$ to confirm convergence in probability.

Having proved the claim, recall that 
\begin{align*}
\mathbb P(R_{n,\theta_n} \in A_n) &= \frac{\sum_{i=1}^n e^{\theta_n X_{i}}\1_{X_{i} \in A_n}}{\sum_{i=1}^n e^{\theta_n X_{i}}} \\
&= \frac{\sum_{i=1}^n e^{\theta_n X_{i}}\1_{X_{i} \in A_n}}{e^{\theta_n \max_i X_{i}}} \frac{e^{\theta_n \max_i X_{i}}}{\sum_{i=1}^n e^{\theta_n X_{i}}}
\end{align*}
The second term converges in probability to $1$ by our claim, while the first term can be written as 
$$
\frac{\sum_{i=1}^n e^{\theta_n X_{i}}\1_{X_{i} \in A_n}}{e^{\theta_n \max_i X_{i}}} = \1_{\max_{i} X_{i} \in A_n} + \sum_{i=1}^n e^{\theta_n (X_{i} - \max_i X_{i})}\1_{X_{i} \in A_n, X_{i} \neq \max_j X_j}.
$$
Again, by the claim, the second term converges to $0$ in probability. It now follows that $$
\mathbb P(R_{n,\theta_n} \in A_n) - \mathbb P(\max_i X_{i} \in A_n) \to 0
$$
in probability. In particular, for any fixed positive $C$ take $A_n = B - C(\mx-F^{-1}(1-\frac 1n))$. Then, clearly $$
\mathbb P(\max_i X_{i} \in A_n) = \mathbb P\left(\frac{\mx-\max_i X_{i}}{\mx- F_X^{-1}(1-\frac 1n)} \geq C\right) \to \mathbb P(-W_{\alpha} \geq C).
$$
Therefore, we have $\mathbb P(R_{n , \theta_n} \in A_n) \to \mathbb P(-W_{\alpha} \geq C)$, which rearranges itself to $$
P\left(\frac{\mx-R_{n,\theta_n}}{\mx- F_X^{-1}(1-\frac 1n)} \geq C\right) \to \mathbb P(-W_{\alpha} \geq C)
$$
as desired. This completes the proof.
\end{proof}

\section{Proofs for Section~\ref{mdweibull}}\label{appendix:e}

In this section, we prove the results from Section~\ref{mdweibull}. Recall that $X$ is now a vector, $\theta\in \R^d$ a fixed unit vector, $c_n>0$ a sequence tending to infinity, and the distributions of $X_{c\theta}$ and $R_{n,c\theta}$ are given by \eqref{twisted} and \eqref{empirical}.

We begin with the proof of Lemma~\ref{KStoSLhd}, which is essentially the same as the proof of Lemma~\ref{KStoSL}. 

\begin{proof}
Let $a_n$ be a sequence of vectors and $b_n>0$. Then, for any axis-aligned hypercube $R$, observe that $\{\frac{r-a_n}{b_n} : r \in R\}$ is still an axis-aligned rectangle. Furthermore, this map is invertible, and thus is a $1-1$ mapping on $\mathcal{R}$, the set of all axis-aligned hypercubes on $\mathcal{R}$. Therefore, for any $n\geq 1$, \begin{equation}
\sup_{R} |\mathbb P(X_{1,n} \in R)-\mathbb P(X_{2,n} \in R)| = \sup_{R} \left|\mathbb P\left(\frac{X_{1,n}-a_n}{b_n} \in R\right)-\mathbb P\left(\frac{X_{2,n}-a_n}{b_n} \in R\right)\right|.\label{touchtwo}
\end{equation}
For part(a), if $Z_1=Z_2=Z$ is continuous, then by Lemma~\ref{assist} we know that the distribution of $\frac{X_{i,n}-a_n}{b_n}$ converges uniformly on compacts to $F_{Z}$. Given this, part(a) is clear by \eqref{touchtwo}.

On the other hand, for part (b) we note that if $Z_1 \neq Z_2$ then one can find a continuity point $c$ of $Z_1,Z_2$ at which their CDFs differ, by Lemma~\ref{assist}. However, if $\mathcal{R}_c = \{x : x_i \leq c_i, 1 \leq i \leq d\}$ then the difference between the CDFs at $c$ equals $|\mathbb P(Z_1 \in \mathcal{R}_c)-\mathbb P(Z_{2} \in \mathcal{R}_c)| = \epsilon>0$, which is non-negative. Let $c'$ be such that $F_{Z_1}(c'), F_{Z_2}(c') < \epsilon/2$. Then, it follows that $$
|\mathbb P(Z_1 \in \mathcal{R}_c \setminus \mathcal{R}_{c'})-\mathbb P(Z_{2} \in \mathcal{R}_c \setminus \mathcal{R}_{c'})| > \epsilon-\epsilon/2-\epsilon/2>0.
$$
Since $X_{i,n} \to Z_i$ in distribution, it follows that their CDFs cannot uniformly converge to the same limit on $\mathcal{R}_c \setminus \mathcal{R}_{c'}$, completing the proof.
\end{proof}

Next, we sketch a proof of Lemma~\ref{examplesmdweibull}, the list of examples of random vectors satisfying Assumption~\ref{ass:mvrv}.

\begin{proof}[Proof of Lemma~\ref{examplesmdweibull}]

We begin with the one-dimensional example. Let $X$ be a random variable with maximum value $B$ which is in the Weibull regime with parameter $\alpha>0$. Let $F_X$ be the CDF of $X$.

We will take $x_{\theta} = \mx$. Let $\theta>0$ be arbitrary and $$A = [0,c]$$ for some arbitrary $c>0$. Let $U(t) = 1-F_X(\mx-t)$, which is regularly varying with index $\alpha>0$ by Assumption~\ref{weibull}. Then, $$
\frac 1{U(t)} \mathbb P\left(\frac{\mx-X}{t} \in A\right) = \frac 1{U(t)} \mathbb P(X \geq \mx-tc)  = \frac{1-F_X(\mx-tc)}{1-F_X(\mx-t)} = c^{\alpha}.
$$
Thus, it follows that $\nu([0,c]) = c^{\alpha}$ for all $c>0$, which uniquely determines $\nu$. It is sufficient to consider intervals to prove vague convergence, hence $\nu$ is the vague limit under consideration. 

For the integrability condition, note that $\theta^T y = \theta y$ is uniquely maximized at $y=\mx$. Then, $\nu(\{0<\theta y\leq 1\}) < \infty$ holds since $$
\nu(\{0<\theta y\leq 1\}) = \int_{0}^{\frac 1{\theta}} e^{-y} d \nu(y) < \Gamma(\alpha)<\infty.
$$
We have completed the proof.

For the second example, let $(X_1,X_2,\ldots,X_d)$ be a bounded random vector with independent components, and let $x = (\mx_1,\ldots,\mx_d)$ be the vector of maximal values of the $X_i$. Suppose $X_i$ is in the Weibull regime with index $\rho_i$. We let $\theta$ be any vector with positive entries so that $\theta^Ty$ is uniquely maximized at $x$.

By an approximation argument, it is sufficient to find a function $U(t)$ and a measure $\nu$ such that $$
\frac{1}{U(t)} \mathbb P(x-X \in tA) \to \nu(A)
$$
for all hypercubes $A$ of the form $A = \times_{i=1}^n [a_i,b_i]$ where $a_i \leq b_i$ for $1 \leq i \leq d$. However, in this case, $$
x \in x- tA \text{ if and only if } X_i \in [\mx_i - ta_i, \mx_i - tb_i] \text{ for all } 1 \leq i \leq d.
$$
Therefore, $$
\mathbb P(x-X \in tA)  = \prod_{i=1}^d \mathbb P(X_i \in [\mx_i-ta_i, \mx_i - tb_i]).
$$
Let $U_i(t) = \frac {1}{1-F_{X_i}(\mx-t)}$. Then,
$$\frac 1{U_i(t)}\mathbb P\left(X_i \in [\mx_i-ta_i, \mx_i - tb_i]\right) = \frac 1{U_i(t)}\mathbb P\left(\frac{\mx-X_i}{t} \in \left[a_i,b_i\right]\right) \to \nu_i([a_i,b_i])$$
by Assumption~\ref{weibull}, where $\nu_i([0,y]) = y^{\rho_i}$. Thus, taking the product of the above convergence over all $i \in 1,\ldots,d$,
$$
\frac 1{\prod_{i=1}^d U_i(t)} \mathbb P(x-X \in tA) \to \nu(A)
$$
where $\nu = \otimes_{i=1}^n \nu_i$ is the product measure. Since the product of regularly varying functions is still regularly varying, the result follows as long as we can demonstrate the integrability condition.

However, recall that$\theta_i>0$ for all $1 \leq i \leq d$. We now have $$
\nu(\{0 < \theta^Ty < 1\}) \leq \nu\left(\times_{i=1}^d \left[0,\frac{1}{\theta_i}\right]\right) = \prod_{i=1}^d \int_0^{\frac{1}{\theta_i}} e^{-y}d \nu_i(y) < \prod_{i=1}^d \Gamma\left(\frac 1{\theta_i}\right)< \infty
$$
by the definition of $\nu_i$. This proves the integrability result.

We argue (d) first, since its analysis will essentially resolve (c). Suppose that $X$ is a truncated normal $X = N(0,I)\1_{N(0,I) \leq r}$ for some $r>0$. Fix $y \in \R^d, \|y\| = r$ and let $\theta = \frac 1r y$. Note that $\theta^Tz$ is maximized at $z = y$ within the support of $X$. Furthermore, note that $N(0,1)$ possesses a density $f$ which is continuous at $y$, and $f(y)>0$. Thus, for any $t>0$ and compact set $A$ with measure-zero boundary, 
\begin{align*}
\mathbb P\left(\frac{y-X}{t} \in A\right) &= \mathbb P(X \in y-tA) \\
&= \int_{(y-tA) \cap \{\|x\| \leq r\}} f(z) dz\\
&= t^d\int_{A \cap\{\|y-tz\| \leq r\}} f(y-tz) dz.
\end{align*}

As $t\to 0$ we have $A \cap \{\|y-tz\| \leq r\} \to A \cap \{\theta^T z >0\}$, and $f(y-tz) \to f(y)$. By the dominated convergence theorem, $$
\frac 1{f(y)t^d}\mathbb P\left(\frac{y-X}{t} \in A\right)  \to |A \cap \{\theta^T z >0\}|.
$$
Thus, $\nu(A) = |A \cap \{\theta^T z >0\}|$ is the Lebesgue measure restricted to $\{\theta^T z >0\}$. In this case, the integrability condition \eqref{intcond} does \emph{not} hold, however, since $\nu\{\{0<\theta^t y<1\}\} = +\infty$. 

Now, if $X$ is a uniform random vector on any polytope $P$, then we may argue, similarly to the previous case, that $U(t) = t^d$ and $\nu(A)$ is the Lebesgue measure on $A$, not restricted to $\mathbb R^d$ but rather to the conical hull of finitely many "extreme" vectors $v_1,\ldots,v_d$ originating from $0$. In this case, note that the intersection of this cone with $\{0<\theta^Ty<1\}$ has finite measure, which is tantamount to $\nu(\{0<\theta^Ty<1\}) < \infty$, as desired.
\end{proof}

We now prove Lemma~\ref{weibtoweibhd}.

\begin{proof}[Proof of Lemma~\ref{weibtoweibhd}]
    Let $X$ satisfy Assumption~\ref{ass:mvrv} at $x_{\theta}$ with function $U(t)$ and vague limit $\nu$. For the sake of clarity, we assume $(x_{\theta})_i \neq 0$ for all $1 \leq i \leq d$ : the proof is easily modified if one of these doesn't hold. This implies that for every bounded rectangle $A$ we have
    $$
\frac 1{U(t)}\mathbb P\left(\frac{x_{\theta}-X}{t} \in A\right) \to\nu(A).
    $$
    Let $A = \times_{i=1}^d [a_i,b_i]$ with $a_i \leq b_i$ and suppose $g(x) = (x_i^{\alpha_i})_{1 \leq i \leq d}$. Now, we have 
    \begin{align*}\left\{\frac{g(x_{\theta})-g(X)}{t} \in A\right\} = \left\{g(X) \in g(x_{\theta}) - tA\right\} =&\{X \in g^{-1}(g(x_{\theta})  - tA)\} \\=& \left\{\frac{x-X}{t} \in \frac{x_{\theta}-g^{-1}(g(x_{\theta})  - tA)}{t}\right\}
    \end{align*}

Expanding out the definition of $g$ and $g^{-1}$, $$
\frac{x_{\theta}-g^{-1}(g(x_{\theta})  - tA)}{t} = \times_{i=1}^d \left[\frac{x_{\theta_i} - (((x_{\theta})_i)^{\alpha_i}- ta_i)^{\frac 1{\alpha_i}}}{t}, \frac{x_{\theta_i} - (((x_{\theta})_i)^{\alpha_i}-tb_i)^{\frac 1{\alpha_i}}}{t}\right]
$$
Observe that as $t \to 0$, the above set converges : $$
\times_{i=1}^d \left[\frac{x_{\theta_i} - (((x_{\theta})_i)^{\alpha_i}- ta_i)^{\frac 1{\alpha_i}}}{t}, \frac{x_{\theta_i} - (((x_{\theta})_i)^{\alpha_i}-tb_i)^{\frac 1{\alpha_i}}}{t}\right] \to \times_{i=1}^d \frac{\left[a_i,b_i\right]}{\alpha_ix_{\theta_i}^{\alpha_i-1}}.
$$
Thus, it follows that $$
\frac 1{U(t)}\mathbb P\left(\frac{g(x_{\theta})-g(X)}{t} \in A\right) \to \mu(B)
$$
where if $B = \times_{i=1}^d [a_i,b_i]$ with $a_i \leq b_i$ then $\mu(B) = \nu\left(\times_{i=1}^d \frac{\left[a_i,b_i\right]}{\alpha_ix_{\theta_i}^{\alpha_i-1}}\right)$. That is, $\mu$ is just a rescaling of $\nu$. This completes the proof.
\end{proof}

Unlike the previous appendix section, we prove Theorem~\ref{m2tmt2hd} and Theorem~\ref{sltruemd} in the same subsection, since their proofs are extremely similar. The three parts of Theorem~\ref{thm:emphd} will then be covered later.

\subsection{Proof of Theorem~\ref{m2tmt2hd} and Theorem~\ref{sltruemd}}

The following lemma will be extremely helpful in this subsection. Let $X$ satisfy Assumption~\ref{ass:mvrv} at $x_{\theta}$ with limiting measure $\nu$ and regularly varying function $U(t)$. 

\begin{lemma}\label{genforhd}
Let $B\subset \mathbb R^d$ be such that $\nu(\partial B) = 0$, where $\partial B$ is the topological boundary of $B$. We have $$
\frac{\mathbb E[e^{c_n \theta^T X}\1_{X \in x_{\theta} - \frac{1}{c_n}B}]}{e^{c_n \theta^T x_{\theta}} U(\frac 1{c_n})} \to \int_{B} e^{-\theta^Ty}d\nu(y)
$$
as $c_n \to\infty$.
\end{lemma}
\begin{proof}
We have \begin{align*}
&\E[e^{c_n \theta^T X}\1_{X \in x_{\theta} - \frac{1}{c_n}B}] \\
=& e^{c_n \theta^T x_{\theta}}\int e^{-c_n \theta^T (x_{\theta}-X)}\1_{X \in x_{\theta} - \frac{1}{c_n}B} d \Prob(X) \\
=& e^{c_n \theta^T x_{\theta}}\int e^{-\theta^Ty}\1_{y \in B} d \Prob\left(c_n(x_{\theta}-X)\right) \\
=& e^{c_n \theta^T x_{\theta}}U\left(\frac 1{c_n}\right)\frac{\int e^{-\theta^Ty}\1_{y \in B} d \Prob\left(c_n(x_{\theta}-X)\right)}{U(\frac 1{c_n})}.
\end{align*}
Observe that $e^{-\theta^T y}\1_{y \in B}$ is discontinuous only on $\partial B$, but $\nu(\partial B) = 0$. Thus, as $c_n \to\infty$, the latter term converges to $\int_B e^{-\theta^T y} d \nu(y)$ by Assumption~\ref{ass:mvrv}, completing the proof.
\end{proof}

The proof of the theorem follows from the above lemma rather easily.

\begin{proof}[Proof of Theorem~\ref{m2tmt2hd}]
By Lemma~\ref{genforhd}, \begin{gather*}
\frac{\E[e^{2c_n\theta^TX}]}{U(1/2c_n)e^{2c_n\theta^T x_\theta}} \to \int e^{-\theta^Ty}d\nu(y) \\
\frac{U^2(1/c_n)e^{2c_n\theta^T x_\theta}}{\E[e^{c_n\theta^TX}]^2} \to \frac 1{\left(\int e^{-\theta^Ty}d\nu(y)\right)^2}
\end{gather*}
Multiplying these together, $$
\frac{\frac{\E[e^{2c_n\theta^TX}]}{\E[e^{c_n \theta^T X}]^2}}{U^{-2}(1/c_n)U(1/2c_n)\int e^{-\theta^Ty}d\nu(y)} \to 1.
$$
Note that $$
\frac{U(1/{2c_n})}{U(1/c_n)} \to 2^{\alpha}
$$
as $c_n \to\infty$. Hence, $$
U(1/c_n)\frac{\E[e^{2c_n\theta^TX}]}{\E[e^{c_n \theta^T X}]^2} \to \frac{2^{-\alpha}}{\int e^{-\theta^Ty}d\nu(y)},
$$
which was to be proved.
\end{proof}

The above lemma is also enough to provide a simple proof for the scaling limit of $X_{c_n\theta}$.

\begin{proof}[Proof of Theorem~\ref{sltruemd}]
It suffices to prove that for every set $B \subset \mathbb R^d$ such that $\nu(\partial B) =0$, we have $$
\mathbb P\left(c_n(x_{\theta}-X_{c_n\theta}) \in B\right) \to  \mathbb P\left(Z\in B\right).
$$
However, observe that \begin{align*}
&\mathbb P\left(c_n(x_{\theta}-X_{c_n\theta}) \in B\right) \\
=&\mathbb P\left(X_{c_n\theta} \in x_{\theta} - \frac{1}{c_n}B\right) \\
& = \frac{\E[e^{c_n \theta^T X\1_{X \in x_{\theta} - \frac{1}{c_n}B}}]}{\E[e^{c_n \theta^T X}\1_{X \in \mathbb R^d}]} \\
& \to \frac{\int_{B} e^{-\theta^Ty}d\nu(y)}{\int e^{-\theta^Ty}d\nu(y)},
\end{align*}
where we applied Lemma~\ref{genforhd} once for the numerator and once for the denominator with $B$ replaced by $\mathbb R^d$. This completes the proof.
\end{proof}

In the next subsection, we will prove Theorem~\ref{thm:emphd}.

\subsection{Proof of Theorem~\ref{thm:emphd}}

The below proof of Theorem~\ref{thm:emphd}(a) is extremely similar to the proof of Theorem~\ref{thm:slemp1D}(a).

\begin{proof}[Proof of Theorem~\ref{thm:emphd}(a)]
We will use Proposition~\ref{posassist}. Observe that $\frac{M_{\theta_n}}{n} \to 0$. Let $B$ be any Borel set such that $\nu(B)>0$ and $\nu(\partial B) = 0$, and let $B_n = x_{\theta} - \frac{B}{c_n}$. By Theorem~\ref{sltruemd} it is clear that $\lim_{n \to \infty} \frac{\E[e^{\theta_n^T X}]}{\E[e^{\theta_n^T X}\1_{X \in B_n}]} < \infty$. Therefore, applying Proposition~\ref{posassist} together with Theorem~\ref{sltruemd} furnishes the proof immediately.
\end{proof}

Next, we will prove Theorem~\ref{thm:emphd}(b).

\begin{proof}[Proof of Theorem~\ref{thm:emphd}(b)]
Suppose that $\frac{M_{c_n \theta}}{n} \to C\in (0,\infty)$. Then, by Theorem~\ref{m2tmt2hd} we know that $c_nU^{-1}(\frac 1n) \to C_1$ for some $C_1>0$, and hence $\theta_n \to \infty$. Let $D$ be a Borel set, and let $A_n = x_{\theta} - U^{-1}(\frac 1n)D$.

Consider $\mathbb P(R_{n,c_n\theta} \in A_n)$, and divide the top and bottom by $e^{\theta^T x_{\theta}}$ to obtain 
\begin{align*}
\mathbb P(R_{n,c_n\theta} \in A_n) =& \frac{\sum_{i=1}^n e^{c_n\theta^T X_{i}}\1_{X_{i} \in A_n}}{\sum_{i=1}^n e^{c_n\theta^T X_{i}}}\\ =&  \frac{\sum_{i=1}^n e^{c_n\theta^T(X_{i}-x_{\theta})}\1_{X_{i} \in A_n}}{\sum_{i=1}^n e^{c_n\theta^T(X_{i}-x_{\theta})}} \\
=& \Phi_n\left(\sum_{i=1}^n e^{c_n\theta_n^T(x_{\theta}-X_i)}1_{A_n}, \sum_{i=1}^n e^{c_n\theta_n^T (x_{\theta}-X_i)}\right),
\end{align*}
where 
$$
\Phi(y,z) = \frac{y}{z}.
$$

At this point, we verify the hypotheses of Lemma~\ref{dekhliyo}. We take $c_n = \|\theta_n\|$ and $\theta$ as in the lemma itself. Then, $X_{\theta_n} \to x_{\theta}$ in probability. Let $x_n = x_{\theta}$, $a_n =U^{-1}(\frac 1n)$, $C_1$ be as above, $\nu$ be as in \eqref{ass:mvrv} and $D$ as chosen. All the hypotheses of Lemma~\ref{dekhliyo} are easily verified.

Since $\Phi$ is a continuous mapping, and $\int e^{-C_1\theta^Ty} dPRM(\nu) \neq 0$ with probability $1$ by the definition of a PRM, by Lemma~\ref{dekhliyo} it follows that \begin{align*}
\mathbb P(R_{n,c_n\theta} \in A_n) = &\Phi_n\left(\sum_{i=1}^n e^{c_n\theta^T(x_{\theta}-X_i)}1_{A_n}, \sum_{i=1}^n e^{c_n\theta^T(x_{\theta}-X_i)}\right) \\ \overset{d}{\to} & \frac{\int e^{-C_1y}1_{y \in D} d PRM(\nu)}{\int e^{-C_1y} dPRM(\nu)}.
\end{align*}

Finally, note that $$
\mathbb P(R_{n,c_n\theta} \in A_n) = \mathbb P(c_n\theta(x_{\theta} - R_{n,c_n\theta}) \in c_n(x_{\theta} - A_n))
$$
Combining the two statements above, if $Z$ is a random vector such that
$$
\mathbb P(Z \in D) =  \frac{\int e^{-C_1y}1_{y \in D} d PRM(\nu)}{\int e^{-C_1y} dPRM(\nu)},
$$
then $$
\mathbb P(c_n(x_{\theta}-R_{n,c_n\theta}) \in D) \to \mathbb P(C_1Z\in D).
$$
Thus, the random variable $Z_{C_1,PRM} = C_1Z$ is the desired limit. Note that this random variable depends upon the Poisson random measure, while the scaling limit $Z$ from part(a) of the theorem does not. It follows that these two random variables are not the same, which concludes the proof.

We remark that, as in the one-dimensional case, the limiting random variable $Z_{C_1,PRM}$ is continuous, but do not prove this here.
\end{proof}

The only non-trivial part of the proof of Theorem~\ref{thm:emphd}(c) is finding the scaling limit for the sample maximizer of $\theta^T X_{(i)}$. Once this is done, our proof will follow exactly as in the one-dimensional case, and we can omit the proofs of these parts. The next lemma explains the origin of the random variable $V$ in Theorem~\ref{thm:emphd}(c) as ths aforementioned scaling limit.

\begin{lemma}\label{help4}
Let $X_1,X_2,...,X_n$ be iid $X$ and let $X_{(n)} = \argmax_i \theta^T X_{i}$ be the sample maximizer of $\theta^T X_{i}$. Then, $$
\frac{x-X_{(n)}}{U^{-1}(\frac 1n)} \to V
$$
for some random variable $V$ depending upon $\nu$.
\end{lemma}
\begin{proof}
Let $\mathcal{S} = \{w : \theta^T w = 1\}$, and define $T : \{y : \theta^Ty>0\} \to (0,\infty) \times \mathcal{S}$ by $$
T(y) = \left(\theta^Ty, \frac{y}{\theta^Ty}\right).
$$
Clearly $T$ is continuous on this space, and the support of $\frac{x-X}{t}$ is contained in the domain of $T$ for all $t>0$.

We claim, furthermore, that $T^{-1}(K)$ is compact for every $K \subset (0,\infty) \times \mathcal{S}$ which is compact. This is easy to see : indeed, if $K$ is compact then $K$ is contained in a set of the form $[a,b] \times J$ where $J \subset \mathcal{S}$ is compact. However, it's easy to see that $T^{-1}([a,b] \times J)$ is a closed and bounded rectangle, hence compact (this is particularly easy to see if $\theta$ is a multiple of one of the standard basis vectors, for instance). Hence, $T^{-1}(K)$ is a closed subset of a compact set, hence compact.

By \cite[Proposition 3.18]{Resnick1987}, we have that if $\frac 1{U(t)} \mathbb P(\frac{x-X}{t} \in \cdot) \to \nu(\cdot)$ vaguely, then $$
\frac 1{U(t)} \mathbb P\left(\left(\frac{\theta^T(x-X)}{t}, \frac{x-X}{\theta^T(x-X)}\right) \in \cdot\right) \to \nu \circ T^{-1}(\cdot)
$$
vaguely. But by Lemma~\ref{lem:nu}(b), $\nu \circ T^{-1}$ is a product measure, of the form $$
\nu \circ T^{-1} = dr \times ds.
$$

Now, by \cite[Proposition 3.21]{Resnick1987} and a restriction argument similar to the proofs in \cite[Section 3.3.2]{Resnick1987} we obtain $$
\sum_{i=1}^n \delta_{(\theta^T(x-X_{i}), \frac{x-X_{i}}{\theta^T(x-X_{i})})} \Rightarrow PRM(\nu \circ T^{-1}) = PRM(dr \times ds),
$$

 where we use $dr \times ds$ to denote the spherical and angular parts of $d\nu \circ T^{-1}$. By the continuous mapping theorem and L, it follows that $\frac{x-X_{(n)}}{U^{-1}(1/n)}$ converges to the atom $V$ of the PRM on which $\theta^Ty$ is minimized. 

It is not difficult to see the distribution of $V$. Note that the point $V$ is located at a point $(a,s)\in [0,\infty) \times \mathcal{S}$, if and only if the PRM has no point in the set $\{w : 0<\theta^T w < a\}$, which occurs with probability $e^{-\nu \circ T^{-1}((0,a] \times \mathcal{S})}$. Thus, it follows that 
$$
\mathbb P(V \in B) = \int_{B} e^{-\nu \circ T^{-1}((0,r] \times \mathcal{S})} PRM(dr \times ds).
$$
\end{proof}

We now complete the proof of Theorem~\ref{thm:emphd}(c).

\begin{proof}[Proof of Theorem~\ref{thm:emphd}(c)]

Suppose, without loss of generality that $X_1$ is the sample maximizer of $\theta^T x$. We claim that in probability, $$
\frac{\sum_{i=1}^n e^{c_n\theta^TX_{i}}}{e^{c_n\theta^TX_1}} \to 1,
$$
if $c_nU^{-1}(\frac 1n) \to +\infty$, which is implied by Theorem~\ref{m2tmt2hd} and $\frac{M_{\theta_n}}{n} \to \infty$. In particular, this implies that
\begin{equation}
\frac{\sum_{i=2}^n e^{c_n\theta^T X_{i}}}{e^{c_n\theta^TX_1}} \to 0.
\end{equation}

 To prove this claim, we have by Lemma~\ref{pppc} that $\sum_{i=1}^n \delta_{\frac{x_{\theta}-X_{i}}{U^{-1}(\frac 1n)}} \Rightarrow PRM( \nu)$ where $\nu$ is as in Assumption~\ref{ass:mvrv}. Now, note that the smallest point of the former point process is attained when $X_{i} = \max_j X_j$, which without loss of generality we assume is attained at $X_1$. Therefore, we get by Lemma~\ref{pppconv} that 
$$
\sum_{i=1}^n \delta_{\frac{X_1-X_{i}}{U^{-1}(\frac 1n)}} \Rightarrow PRM( \nu) - s_{PRM( \nu)}.
$$

Let $M>0$ be an arbitrary but fixed parameter. Let $f(x) = e^{-Mx}\1_{\|x\|<M}$. Note that $f$ is continuous and compactly supported. Thus, by the continuous mapping theorem (see \cite[Section 3.5]{Resnick1987}) we have \begin{equation}\label{eq:conv1}
\sum_{i=1}^n e^{-M\frac{\theta^T(X_1-X_{i})}{U^{-1}(\frac 1n)}}\1_{\|X_1-X_i\| < U^{-1}(\frac 1n)M} \xrightarrow{d} \int e^{-M\theta^T x}\1_{\|x\|<M}\left[PRM( \nu) - s_{PRM( \nu)}\right].
\end{equation}

Now, it is clear that for any $\epsilon>0$,
\begin{align*}
&\lim_{M \to \infty} \mathbb P\left[\sum_{i=1}^n e^{-c_n\theta^T\left(X_1-X_i\right)}\1_{\|X_1-X_i\| < MU^{-1}(\frac 1n)} > 1+\epsilon\right] \\ =& \mathbb P\left[\sum_{i=1}^n e^{-c_n\theta^T(X_i-X_1)} > 1+\epsilon\right].
\end{align*}

For any $\epsilon>0$ and $M>0$, employing \eqref{eq:conv1} gives \begin{align}
&\limsup_{n} \Prob\left[\sum_{i=1}^n e^{-c_nU^{-1}(\frac 1n)\theta^T \left(\frac{X_1-X_i}{U^{-1}(\frac 1n)}\right)}\1_{\|X_1-X_i\| < MU^{-1}(\frac 1n)} > 1+\epsilon\right] \nonumber \\
\leq & \limsup_{n \to \infty}\Prob\left[\sum_{i=1}^n e^{-M\left(\frac{\theta^T(X_1-X_i)}{U^{-1}(\frac 1n)}\right)}\1_{\|X_1-X_i\| < MU^{-1}(\frac 1n)} > 1+\epsilon\right] \nonumber\\ = & \Prob\left[\int e^{-M\theta^Tx}\1_{\|x\|<M} (PRM(\nu) - s_{PRM(\nu)}) > 1+\epsilon\right]. \label{seven}
\end{align}

Note that $e^{-M\theta^Tx}\1_{\|x\|<M}$ is dominated by $e^{-\theta^Tx}$ for $M>1$, since the support of $\nu$ is contained in $\{\theta^T y >0\}$. Furthermore, as $M \to \infty$, we have $e^{-M\theta^Tx}\1_{\|x\|<M} \to \1_{\theta^T y =0}$ pointwise. In order to use the dominated convergence theorem, we need to prove that $e^{-\theta^Tx}$ is integrable under $PRM(\nu) - s_{PRM(\nu)}$ a.s. : we will prove the stronger statement that $e^{-\theta^Tx}$ is integrable under $PRM(\nu)$ a.s.

In order to prove this, observe that $(1-e^{-\theta^Tx}) \leq \theta^Tx$ whenever $\theta^Tx>0$. Thus, 
$$\int_{\theta^Tx>0} (1-e^{-e^{-\theta^Tx}})d\nu(x) \leq \int_{\theta^Tx>0} e^{-\theta^Tx} d\nu(x)  <\infty,$$
by Lemma~\ref{lem:nu}(c). By \cite[Proposition 3.6(ii)]{Resnick1987}, the left hand side is equal to $-\ln \E[e^{-\int_{\theta^Ty>0} e^{-\theta^Tx} dPRM(\nu)}]$. It follows that $\int_{\theta^Ty>0} e^{-\theta^Tx} dPRM(\nu)$ is finite a.s.

By the dominated convergence theorem, the probability in \eqref{seven} converges to $0$. In particular, it follows that 

$$
\limsup_{n} \Prob\left[\sum_{i=1}^n e^{-c_nU^{-1}(\frac 1n)\theta^T \left(\frac{X_1-X_i}{U^{-1}(\frac 1n)}\right)} > 1+\epsilon\right] = 0,
$$
which instantly proves the claim, since the sum must always exceed $1$ (by taking the term $i=1$), therefore it was sufficient to check only exceedance by $\epsilon$ to confirm convergence in probability.

Having proved the claim, recall that
\begin{align*}
\mathbb P(R_{n,c_n\theta} \in A_n) &= \frac{\sum_{i=1}^n e^{c_n\theta^T X_{i}}\1_{X_{i} \in A_n}}{\sum_{i=1}^n e^{c_n\theta^T X_{i}}} \\
&= \frac{\sum_{i=1}^n e^{c_n\theta^T X_{i}}\1_{X_{i} \in A_n}}{e^{c_n\theta^T X_{1}}} \frac{e^{c_n\theta^T X_{1}}}{\sum_{i=1}^n e^{c_n\theta^T X_{i}}}
\end{align*}
The second term converges in probability to $1$ by our claim, while the first term can be written as 
$$
\frac{\sum_{i=1}^n e^{c_n\theta^T X_{i}}\1_{X_{i} \in A_n}}{e^{c_n\theta^T  X_{1}}} = \1_{c_n\theta^TX_{1} \in A_n} + \sum_{i=2}^n e^{c_n\theta^T (X_{i} - X_{1})}\1_{X_{i} \in A_n}.
$$
Again, by the claim, the second term converges to $0$ in probability. It now follows that $$
\mathbb P(R_{n,c_n\theta} \in A_n) - \mathbb P(X_1\in A_n) \to 0
$$
in probability. In particular, for any Borel set $B$ take $A_n = x_{\theta} - BU^{-1}(\frac 1n)$. Then, clearly $$
\mathbb P(X_{1} \in A_n) = \mathbb P\left(\frac{x_{\theta}- X_{1}}{U^{-1}(\frac 1n)} \in B\right) \to \mathbb P(V \in B).
$$
Therefore, we have $\mathbb P(R_{n , c_n\theta} \in A_n) \to \mathbb P(V \in B)$, which rearranges itself to $$
P\left(\frac{x_{\theta}-R_{n,c_n\theta}}{U^{-1}(\frac 1n)} \in B\right) \to \mathbb P(V \in B)
$$
as desired. This completes the proof.
\end{proof}

\section{Proofs for Section~\ref{unbdd}}\label{appendix:f}

In this section we will prove results from Section~\ref{unbdd}. Recall, in this section, that $X$ is a continuous random vector with full support on $\R^d$, and density $f$. Also recall that there are $\alpha>1, K,L>0$ such that \eqref{eq:growth} holds. Fix $\theta \in \R^d, \|\theta\|=1$ and $c_n>0$ which will tend to infinity.

For completeness, we restate Lemma~\ref{lem:phic}. However, we also include a part(c) in this statement which was omitted in the original, since it is important for the proof but not for the statements of the theorems in this section.

\begin{lemma}
For each $c>0$,  the function $\Phi_c(x)$\begin{enumerate}[label = (\alph*)]
\item attains its unique maximum at a point $m_{c} =  \left(\frac{c}{\alpha K}\right)^{\frac 1{\alpha-1}}\theta$, and $$
\Phi_c(m_c) =  (\alpha-1)\alpha^{-\alpha/(\alpha-1)}c^{\alpha/(\alpha-1)}K^{-1/(\alpha-1)}.
$$
\item is thrice continuously differentiable in $\mathbb R^d \setminus \{0\}$, and \begin{gather*}
\nabla^2 \Phi_c(m_c) =  -K\alpha \|m_c\|^{\alpha-2} (I+(\alpha-2)\theta \theta^T).
\end{gather*}
(Note : in the second term, the product of vectors is the outer product, not the inner product. Hence it leads to a rank-one matrix).
\item is uniformly second-order Taylor approximable for large $c$ i.e. for every $\epsilon>0$, there exists $c_0$ such that $c>c_0$ and $\|x-m_c\| \leq \frac{\|m_c\|}{2}$ implies 
\begin{equation*}
\left|\Phi_c(x) - \Phi_c(m_c) - \frac 12 (x-m_c)^T \nabla^2\Phi_c(m_c) (x-m_c)\right| \leq \frac{\epsilon}{2} (x-m_c)^T \nabla^2\Phi_c(m_c) (x-m_c).
\end{equation*}
\end{enumerate}
\end{lemma}

\begin{proof}
Note that as $\|x\| \to \infty$, since $\theta^T x \leq \|\theta\| \|x\|$ by the Cauchy Schwarz inequality, it follows that $\lim_{\|x\| \to \infty} \Phi_c(x) = -\infty$. Furthermore, $\Phi_c$ is differentiable everywhere except possibly at $0$, where $\Phi_c(0) = 0$. Therefore, the maximum of $\Phi_c$ must be attained at one of its critical points provided we show that the value at this point is positive. Write $\|x\|^{\alpha} = \left(\|x\|^2\right)^{\alpha/2}$, which is easier to differentiate since the derivative of $\|x\|^2$ equals $2x$.

Differentiating, \begin{equation}\label{der}
[\nabla \Phi_c(x)](h) = c \theta^Th - \frac{K\alpha}{2}\|x\|^{\alpha-2} \times 2x^Th.
\end{equation} 
Since we are searching for critical points, this must be set to $0$ for all $h$, whence it follows that $$
c \theta = K\alpha\|x\|^{\alpha-2}x.
$$
This tells us that $x$ is a multiple of $\theta$, say $x = \mathcal{C}\theta$. But then $$
c = K \alpha\mathcal{C}^{\alpha-1} \implies  \mathcal{C} = \left(\frac{c}{\alpha K}\right)^{\frac 1{\alpha-1}}.
$$
It follows that $\Phi_c$ possesses a unique critical point at $m_c = \left(\frac{c}{\alpha K}\right)^{\frac 1{\alpha-1}}\theta$. A quick computation now reveals that $$
\Phi_c(m_c) = c\mathcal{C} - K\mathcal{C}^{\alpha} = (\alpha-1)\alpha^{-\alpha/(\alpha-1)}c^{\alpha/(\alpha-1)}K^{-1/(\alpha-1)}
$$
which is a positive constant, completing the proof of part (a).

To compute the Hessian $\nabla^2 \Phi_c$ it suffices to differentiate \eqref{der} once more, leading to $$
\nabla^2 \Phi_c(x) = -K\alpha \left((\alpha-2)\|x\|^{\alpha-4}xx^{T} + \|x\|^{\alpha-2}I\right).
$$
Substituting $x = m_c$ gives \begin{equation}\label{hessian}
\nabla^2 \Phi_c(m_c) = - K\alpha \|m_c\|^{\alpha-2} (I+(\alpha-2)\theta \theta^T),
\end{equation}
completing the proof of part(b).

The third derivative of $\Phi_c$ can be obtained by differentiating the above expression once more. We do not include the calculations, but it suffices to see that for some constant $D>0$, $\lim_{\|x\| \to \infty} \frac{\nabla^3\Phi_c(x)}{D \|x\|^{\alpha-3}} = 1$. (This is rather easy to observe in the one-dimensional case, for instance). 

Therefore, by the Taylor Theorem with remainder applied to $\Phi_c(x)$ in the ball $B_c = \left\{x : \|x-m_c\| \leq \|m_c\|/2\right\}$, 
\begin{equation}\label{appr}
\Phi_c(x) = \Phi_c(m_c) + \frac 12 (x-m_c)^{T} \nabla^2 \Phi_c(m_c) (x-m_c)+ R_{3}(x),
\end{equation}
where $$\|R_3(x)\| \leq \sup_{\|x-m_c\| \leq \|m_c\|/2} \|\nabla^3\Phi_c(x)\| \leq M\|m_c\|^{\alpha-3}$$
for some constant $M$ and all $c$ large enough. In particular, this along with \eqref{hessian} implies that
$$
\lim_{c_n \to\infty} \sup_{x \in B_c}\frac{\|R_3(x)\|}{\frac{1}{2}\|(x-m_c)^T  \nabla^2 \Phi_c(m_c) (x-m_c)\|} = 0.
$$
For any $\epsilon>0$, we may use the above limit to see that for $c_n>0$ large enough, $$
\|R_3(x)\| \leq \frac{\epsilon}{2}\|(x-m_c)^T  \nabla^2 \Phi_c(m_c) (x-m_c)\|
$$
for all $x\in B_c$. Part (c) follows directly from \eqref{appr} and the triangle inequality.
\end{proof}

We will now prove Theorem~\ref{m2tmt2unbdd} and Theorem~\ref{slunbdd}.

\subsection{Proof of Theorem~\ref{m2tmt2unbdd} and Theorem~\ref{slunbdd}}

The following lemma, which we proceed to state, does all the heavy lifting and begets the scaling limit of $X_{c_n \theta}$ and the growth of $M_{c_n \theta}$ as corollaries.

\begin{lemma}\label{thm:unbddas}
Let $B \subset \mathbb R^d$ be any Borel set. Define 
$$
B_{c} = m_c + (-\nabla^2 \Phi_c(m_c))^{-\frac 12}B,
$$
where $m_c, \Phi_c$ are as in Lemma~\ref{lem:phic}. Then,
$$
\lim_{c_n \to\infty}\frac{\E[e^{c_n\theta^T X}1_{X \in B_{c_n}}]}{\frac{e^{\Phi_{c_n}(m_{c_n})}}{\sqrt{\det(-\nabla^2 \Phi_{c_n}(m_{c_n}))}}} = L\mathbb P(N(0,I) \in B).
$$
\end{lemma}

\begin{proof}
For notational convenience, let $B_{c_n} = \sqrt{\det(-\nabla^2 \Phi_{c_n}(m_{c_n}))}$. The proof will be furnished in two parts. We will prove that 
\begin{equation}\label{lbd}
\liminf_{c_n \to\infty} \frac{\E[e^{c_n\theta^T X}1_{X \in B_{c_n}}]}{\frac{e^{\Phi_{c_n}(m_{c_n})}}{B_{c_n}}}\geq L \mathbb P(N(0,I) \in B),
\end{equation}
and subsequently we will prove that 
\begin{equation}\label{ubd}
\limsup_{c_n \to\infty}\frac{\E[e^{c_n\theta^T X}1_{X \in B_{c_n}}]}{\frac{e^{\Phi_{c_n}(m_{c_n})}}{B_{c_n}}} \leq L \mathbb P(N(0,I) \in B),
\end{equation}
completing the proof.

To start the proof of \eqref{lbd} we let $\epsilon, \epsilon'>0$ be arbitrary. By \eqref{eq:growth} there exists $R_{\epsilon'}>0$ such that \begin{equation}\label{lbd1}
f(x) \geq (L-\epsilon')e^{-K\|x\|^{\alpha}}
\end{equation}
for all $\|x\| > R_{\epsilon'}$. By definition of the expectation and \eqref{lbd1} we have
\begin{align}
&\E[e^{c_n\theta^T X}1_{X \in B_{c_n}}] \nonumber\\=& \int_{B_{c_n}} e^{c_n \theta^T x} f(x) dx \nonumber \\
 \geq & (L-\epsilon')\int_{B_{c_n} \cap \{\|x\| \geq R_{\epsilon'}\}} e^{\Phi_{c_n}(x)} dx. \nonumber \\
 \geq & (L-\epsilon') \int_{B_{c_n} \cap \{\|x\| \geq R_{\epsilon'}\} \cap \{\|x-m_{c_n}\| \leq \frac{\|m_{c_n}\|}{2}\}} e^{\left(\Phi_{c_n}(m_{c_n}) - \frac {1-\epsilon}{2} (x-m_{c_n})^{T} \nabla^2 \Phi_{c_n}(m_{c_n}) (x-m_{c_n})\right)} dx,
\label{lbd2}
\end{align}
where in the last step we used Lemma~\ref{lem:phic}(c). We shall now perform a two-step substitution merely for the sake of clarity : this calculation will not be repeated in the lower bound. To begin with, we perform the substitution $y = x-m_{c_n}$. This gives \begin{multline}\label{lbd3}
 (L-\epsilon') \int_{B_{c_n} \cap \{\|x\| \geq R_{\epsilon'}\} \cap \{\|x-m_{c_n}\| \leq \frac{\|m_{c_n}\|}{2}\}} e^{\left(\Phi_{c_n}(m_{c_n}) - \frac{1-\epsilon}{2} (x-m_{c_n})^{T} \nabla^2 \Phi_{c_n}(m_{c_n}) (x-m_{c_n})\right)} dx\\
 =  (L-\epsilon') \int_{(B_{c_n}-m_{c_n}) \cap \{\|y+m_{c_n}\| \geq R_{\epsilon'}\} \cap \{\|y\| \leq \frac{\|m_{c_n}\|}{2}\}} e^{\left(\Phi_{c_n}(m_{c_n}) - \frac{1-\epsilon}{2} y^{T} \nabla^2 \Phi_{c_n}(m_{c_n})y\right)} dy.
\end{multline}
Subsequently, we let $z = (-\nabla^2 \Phi_{c_n}(m_{c_n}))^{-\frac 12}y$. Effecting this change, \begin{multline}\label{lbd4}
(L-\epsilon') \int_{(B_{c_n}-m_{c_n}) \cap \{\|y+m_{c_n}\| \geq R_{\epsilon'}\} \cap \{\|y\| \leq \frac{\|m_{c_n}\|}{2}\}} e^{\Phi_{c_n}(m_{c_n}) - \frac {1-\epsilon}{2} y^{T} \nabla^2 \Phi_{c_n}(m_{c_n})y} dy \\
= (L-\epsilon')\frac{e^{\Phi_{c_n}(m_{c_n})}}{B_{c_n}} \int_{J_{c_n}} e^{-\frac {1-\epsilon}{2} \|z\|^2} dz,
\end{multline}
where \begin{equation}\label{Jc}
J_{c_n} = B \cap \{\|(-\nabla^2 \Phi_{c_n}(m_{c_n}))^{\frac 12}z+m_{c_n}\| \geq R_{\epsilon'}\} \cap \left\{\|(-\nabla^2 \Phi_{c_n}(m_{c_n}))^{\frac 12}z\| \leq \frac{\|m_{c_n}\|}{2}\right\}.
\end{equation}
Combining \eqref{lbd1}, \eqref{lbd2}, \eqref{lbd3} and \eqref{lbd4}, \begin{equation}\label{lbd5}
\frac{\E[e^{c_n\theta^T X}1_{X \in B_{c_n}}]}{\frac{e^{\Phi_{c_n}(m_{c_n})}}{B_{c_n}}} \geq (L-\epsilon') \int_{J_{c_n}} e^{-\frac {1-\epsilon}{2} \|z\|^2} dz.
\end{equation}

The behavior of $J_{c_n}$ as $c_n \to\infty$ will now be analyzed. As a consequence of Lemma~\ref{lem:phic}(b), note that the matrix $M_{c_n} = (-\nabla^2 \Phi_{c_n}(m_{c_n}))^{\frac 12}$ exhibits three kinds of behaviours as $c_n \to\infty$ : if $\alpha>2$ then $M_{c_n} \to 0$. If $\alpha=2$ then $M_{c_n} \to K\alpha I$. Finally, if $\alpha<2$ then $\|M_{c_n}z\| \to \infty$ for every $z \neq 0$. Since $\|m_{c_n}\| \to \infty$ by Lemma~\ref{lem:phic}(a), it follows that as $c_n \to\infty$, $\{\|(-\nabla^2 \Phi_{c_n}(m_{c_n}))^{-\frac 12}z+m_{c_n}\| \geq R_{\epsilon'}\}$ eventually converges to either $\mathbb R^d$ or $\mathbb R^d \setminus \{0\}$.

Similarly, by Lemma~\ref{lem:phic}(b) it follows that the biggest eigenvalue of $\nabla^2 \Phi_{c_n}(m_{c_n})$ grows at the rate of $\|m_{c_n}\|^{2-\alpha/2} < \frac{\|m_{c_n}\|}{2}$ for any $\alpha>1$ and large enough $c_n$. Thus, $
 \left\{\|(\nabla^2 \Phi_{c_n}(m_{c_n}))^{-\frac 12}z\| \leq \frac{\|m_{c_n}\|}{2}\right\}$ converges to $\mathbb R^d$ as $c_n \to\infty$. It follows from these observations and the definition \eqref{Jc} of $J_{c_n}$, that either $J_{c_n} \to B$ or $J_{c_n} \to B \setminus \{0\}$ as $c_n \to\infty$.

Finally, letting $c_n \to\infty$ in \eqref{lbd5}, $$
\liminf_{c_n \to\infty} \frac{\E[e^{c_n\theta^T X}1_{X \in B_{c_n}}]}{\frac{e^{\Phi_{c_n}(m_{c_n})}}{B_{c_n}}} \geq (L-\epsilon') \int_{B} e^{-\frac {1-\epsilon}{2} \|z\|^2} dz.
$$
This equation holds for arbitrary $\epsilon', \epsilon>0$. Thus, letting these parameters tend to zero, the lower bound \eqref{lbd} follows.

We will now prove the upper bound \eqref{ubd}. Once again, let $\epsilon, \epsilon'>0$ be arbitrary constants. By \eqref{eq:growth} there exists $R_{\epsilon'}>0$ such that \begin{equation}\label{ubd1}
f(x) \leq (L+\epsilon')e^{-K\|x\|^{\alpha}}
\end{equation}
for all $\|x\| > R_{\epsilon'}$. Thus, \begin{align}
&\mathbb E[e^{c_n \theta^T X}1_{X \in B_{c_n}}] \nonumber\\
 =& \int_{B_{c_n}} e^{c_n \theta^T X} \nonumber \\
 \leq & \int_{B_{c_n} \cap \{\|x\| \leq R_{\epsilon'}\}} e^{c_n\theta^T X} f(x)dx + (L+\epsilon')\int_{B_{c_n} \cap \{\|x\| > R_{\epsilon'}\}} e^{\Phi_{c_n}(x)}dx. \label{ubd2}
\end{align}

By a discussion similar to that following \eqref{lbd5} it follows that $B_{c_n} \cap \{\|x\| \leq R_{\epsilon'}\}
$ converges to either the empty set or $\{0\}$ as $c_n \to\infty$, which is a set of measure zero. Thus, this term in \eqref{ubd2} doesn't contribute on the right hand side i.e. \begin{equation}\label{ubd3}
\lim_{c_n \to\infty} \frac{\int_{B_{c_n} \cap \{\|x\| \leq R_{\epsilon'}\}} e^{c_n\theta^T x} f(x)dx + (L+\epsilon')\int_{B_{c_n} \cap \{\|x\| > R_{\epsilon'}\}} e^{\Phi_{c_n}(x)}dx}{ (L+\epsilon')\int_{B_{c_n} \cap \{\|x\| > R_{\epsilon'}\}} e^{\Phi_{c_n}(x)}dx} = 1.
\end{equation}

Owing to this realization, we now switch our focus to the second term of \eqref{ubd2}, and for $c_n>0$ large enough apply Lemma~\ref{lem:phic}(b) to get 
\begin{multline}\label{ubd4}
(L+\epsilon')\int_{B_{c_n} \cap \{\|x\| > R_{\epsilon'}\}} e^{\Phi_{c_n}(x)}dx \\ \leq (L+\epsilon')\int_{B_{c_n} \cap \{\|x\| > R_{\epsilon'}\} \cap \{\|x - m_{c_n}\| \leq \frac{\|m_{c_n}\|}{2}\}} e^{\Phi_{c_n}(m_{c_n}) + \frac{1+\epsilon}{2}(x-m_{c_n})^T \nabla^2 \Phi_{c_n}(m_{c_n}) (x-m_{c_n})}dx.
\end{multline}

Just as in the proof of \eqref{lbd}, all we must do now is effect the change of variables and ensure consistency in the limit.  $x = m_{c_n} + (-\nabla^2 \Phi_{c_n}(m_{c_n}))^{\frac 12} z$. Then, 
\begin{multline}\label{ubd5}
(L+\epsilon')\int_{B_{c_n} \cap \{\|x\| > R_{\epsilon'}\} \cap \{\|x - m_{c_n}\| \leq \frac{\|m_{c_n}\|}{2}\}} e^{\Phi_{c_n}(m_{c_n}) + \frac{1+\epsilon}{2}(x-m_{c_n})^T \nabla^2 \Phi_{c_n}(m_{c_n}) (x-m_{c_n})}dx \\=(L+\epsilon')\frac{e^{\Phi_{c_n}(m_{c_n})}}{B_{c_n}} \int_{J_{c_n}} e^{-\frac {1+\epsilon}{2} \|z\|^2} dz,
\end{multline}
where $J_{c_n}$ is as in \eqref{Jc}. As in the discussion there, it follows that $J_{c_n} \to\mathbb R^d$ or $\mathbb R^d \setminus \{0\}$ as $c_n \to\infty$. From here, taking the limit superior on both sides of \eqref{ubd5} and noting the arbitrariness of $\epsilon, \epsilon'$, the proof is immediate. 
\end{proof}

The proofs of the two main theorems are now immediate.

\begin{proof}[Proof of Theorem~\ref{m2tmt2unbdd}]

We proceed exactly as in the proof of Theorem~\ref{thm:m2tmt21D}. Taking $B = \R^d$ in Lemma~\ref{thm:unbddas}, we obtain the following results. \begin{gather*}
\lim_{c_n \to\infty}\frac{\E[e^{2c_n\theta^T X}]}{\frac{e^{\Phi_{2c_n}(m_{2c_n})}}{\sqrt{\det(-\nabla^2 \Phi_{2c_n}(m_{2c_n}))}}} = L \\
\lim_{c_n \to\infty} \frac{\frac{e^{2\Phi_{c_n}(m_{c_n})}}{\det(-\nabla^2 \Phi_{c_n}(m_{c_n}))}}{\E[e^{c_n\theta^T X}]^2} = L^2.
\end{gather*}
Multiplying the two of them, $$
\lim_{c_n \to\infty} M_{c_n \theta} \times \frac{\frac{e^{2\Phi_{c_n}(m_{c_n})}}{\det(-\nabla^2 \Phi_{c_n}(m_{c_n}))}}{\frac{e^{\Phi_{2c}(m_{2c})}}{\sqrt{\det(-\nabla^2 \Phi_{2c}(m_{2c}))}}} = L.
$$
Observe that by Lemma~\ref{lem:phic}(b), we have  $$ 
\det(-\nabla^2 \Phi_{c_n}(m_{c_n})) = K^d\alpha^d m_{c_n}^{d(\alpha-2)} \det(I+(\alpha-2)\theta \theta^T) = K^d\alpha^d m_{c_n}^{d(\alpha-2)} (\alpha-1)
$$
since $I+(\alpha-2)\theta \theta^T$ has the simple eigenvalue $\alpha-1$ and the eigenvalue $1$ with multiplicity $d-1$. Now, using Lemma~\ref{lem:phic}(a),(b), following some algebra we have $$
\lim_{c_n \to\infty} \frac{M_{c_n \theta}}{\exp(p_{\alpha,K,d} c_n^{\alpha/(\alpha-1)}) c_n^{-d(\alpha-2)/(\alpha-1)}} = q_{\alpha,K,L,d}.
$$
where $$
q_{\alpha,K,L,d} = L(\alpha-1)^{1/2} 2^{-d(\alpha-2)/(2\alpha-2)} (\alpha K)^{d/(2\alpha - 2)}
$$
and
$$
p_{\alpha,K,d} = (\alpha-1)\alpha^{-\alpha/(\alpha-1)} K^{-1/(\alpha-1)} (2^{\alpha/(\alpha-1)} - 2).
$$
This completes the proof.
\end{proof}

\begin{proof}[Proof of Theorem~\ref{slunbdd}]
    Let $B \subset \mathbb R^d$ be an arbitrary Borel set. Then, if we apply Lemma~\ref{thm:unbddas} with $B$ and $\R^d$ separately, \begin{equation}\label{conindist}
\mathbb P(X_{c_n \theta} \in B_{c_n}) = \frac{\E[e^{c_n \theta^T X} \1_{X \in B_{c_n}}]}{\E[e^{c_n \theta^T X}]} \to \mathbb P(N(0,I) \in B)
    \end{equation}
    as $c_n \to\infty$, where $B_{c_n} = m_{c_n} + (-\nabla^2\Phi_{c_n}(m_{c_n}))^{-\frac 12} B$. However, $$
\mathbb P(X_{c_n \theta} \in B_{c_n}) = \mathbb P( (-\nabla^2\Phi_{c_n}(m_{c_n}))^{\frac 12}(X_{c_n \theta}-m_{c_n}) \in B).
    $$
    Therefore, \eqref{conindist} is equivalent to saying that $(-\nabla^2\Phi_{c_n}(m_{c_n}))^{\frac 12}(X_{c_n \theta}-m_{c_n}) \overset{d}{\to} N(0,I)$, as desired.
\end{proof}

Now, we will prove Theorem~\ref{empunbdd}.

\subsection{Proof of Theorem~\ref{empunbdd}}

We begin with the proof of Theorem~\ref{empunbdd}(a). Unsurprisingly, this is along the lines of the proof of Theorem~\ref{thm:slemp1D}(a) and Theorem~\ref{thm:emphd}(a).

\begin{proof}[Proof of Theorem~\ref{empunbdd}(a)]
We will use Proposition~\ref{posassist} again. Let $B$ be any Borel set of non-zero measure and let $B_{c_n} = m_{c_n} + (-\nabla^2\Phi_{c_n}(m_{c_n}))^{-\frac 12} B$.

By Lemma~\ref{thm:unbddas} it is clear that $\lim_{n \to \infty} \frac{\E[e^{c_n\theta^T X}]}{\E[e^{c_n\theta^T X}\1_{X \in B_{c_n}}]} < \infty$. Furthermore, $\frac{M_{c_n \theta}}{n} \to 0$ by assumption. Therefore, applying Proposition~\ref{posassist},
$$
\frac{\Prob(R_{n, c_n \theta} \in B_{c_n})}{\Prob(X_{c_n \theta} \in B_{c_n})} \to 1
$$
in probability.   Now, the denominator converges to a known quantity $\mathbb P(N(0,I) \in B)$ by Theorem~\ref{slunbdd}. Hence the numerator must also converge to the same quantity. The proof is immediate.
\end{proof}

 For the proof of Theorem~\ref{empunbdd}(b), we require the following vague convergence result. Following this, the PRM machinery makes everything go through as usual.

\begin{theorem}\label{thm:vague}
Let $a_{c_n} = m_{c_n} \theta$ and $b_{c_n} = \alpha^{-1} m_{c_n}^{1-\alpha}$. Then $$
L^{-1}e^{m_{c_n}^{\alpha}} \alpha^{d} m_{c_n}^{d(\alpha-1)}\mathbb P\left(\frac{X-a_{c_n}}{b_{c_n}} \in \cdot\right) \to \nu(\cdot)
$$
vaguely( where we recall vague convergence from Assumption~\ref{ass:mvrv}), where $\nu$ has density $e^{-\theta^Ty}$ over $\mathbb R^d$.
\begin{proof}
Let $A$ be any compact set, and consider $
\mathbb P\left(\frac{X-a_{c_n}}{b_{c_n}} \in A\right)$. Suppose that $\epsilon,\epsilon'>0$ are arbitrary. Let $R>0$ be large enough such that if $\|x\|>R$ then \begin{equation}\label{est}
(1-\epsilon)Le^{-K\|x\|^{\alpha}} \leq f(x) \leq (1+\epsilon)Le^{-K\|x\|^{\alpha}}.
\end{equation}

Since $A$ is compact and $\|m_{c_n}\| \to \infty$, we have that $b_{c_n} \to0$, and therefore $a_{c_n}+b_{c_n}A \subset \{\|x\| > R\}$ for large enough $c_n>0$. 

By the remainder version of Taylor's theorem applied to $f(y) = \|y\|^{\alpha}$ at the point $a_{c_n}$, for any vector $h \in A$ we have
\begin{equation}\label{tailor}
\|a_{c_n}+b_{c_n}h\|^{\alpha} =  \|a_{c_n}\|^{\alpha} + \alpha b_{c_n}\|a_{c_n}\|^{\alpha -2} a_{c_n}^Th + R_2
\end{equation}
where \begin{align}
R_2 = &\frac 12 h^T\left(b_{c_n}^2\alpha \|a_{c_n}+\lambda b_{c_n}h\|^{\alpha-2}I\right)h \nonumber\\ +& \frac 12 h^T\left(b_{c_n}^2\alpha(\alpha-2)\|a_{c_n} +\lambda b_{c_n}h\|^{\alpha-4}(a_{c_n}+\lambda b_{c_n} h)(a_{c_n}+\lambda b_{c_n}h)^T\right)h\label{rem} \end{align}
for some $\lambda \in (0,1)$. We will now show that this remainder term goes to $0$ uniformly over $h \in A$ and $\lambda \in [0,1]$ as $c_n \to\infty$. 

Since $A$ is compact, $\|\lambda h\|$ is bounded for all $\lambda \in [0,1], h \in A$. Since $\|a_{c_n}\| \to \infty$ and $b_{c_n} \to0$, it follows that $\frac{\|a_{c_n}+ \lambda b_{c_n}h\|}{\|a_{c_n}\|} \to 1$ uniformly over $h \in A$ and $\lambda \in [0,1]$. Furthermore, by the definition of $b_{c_n}$ and $a_{c_n}$ it is clear that $$b_{c_n}^2\|a_{c_n}\|^{\alpha-2} = \|a_{c_n}\|^{-\alpha} \to 0, \quad b_{c_n}^{4}\|a_{c_n}\|^{\alpha-4} = \|a_{c_n}\|^{-3\alpha} \to 0.$$ 
Therefore, both terms in \eqref{rem} go to zero as $c_n \to\infty$, implying that $R_2$ in \eqref{rem} goes to $0$.

So, by \eqref{tailor} we can pick $N$ large enough so that 
\begin{equation}\label{remtozero}
\left|\|a_{c_n}+b_{c_n}h\|^{\alpha} -  \|a_{c_n}\|^{\alpha} - \alpha b_{c_n}\|a_{c_n}\|^{\alpha -2} a_{c_n}^Th \right| \leq \epsilon'.
\end{equation}
for all $n>N$ and $h\in A$. We can now begin manipulating the original expression for $c_n$ large enough, first demonstrating a lower bound, followed by an upper one.

\begin{align}
\mathbb P\left(\frac{X-a_{c_n}}{b_{c_n}} \in A\right) = & \mathbb P\left(X \in a_{c_n}+b_{c_n} A\right)\nonumber \\
\overset{\eqref{est}}{\geq} & (1-\epsilon)L\int_{a_{c_n}+b_{c_n} A} e^{-K\|x\|^{\alpha}} dx \nonumber \\
= & b_{c_n}^{d}\int_{A} e^{-K\|a_{c_n}+b_{c_n}y\|^{\alpha}}dy\nonumber \\
\overset{\eqref{remtozero}}{\geq} & (1-\epsilon)e^{-\epsilon'} L e^{-\|a_{c_n}\|^{\alpha}} b_{c_n}^{d} \int_{A} e^{-\alpha b_{c_n} \|a_{c_n}\|^{\alpha-2}a_{c_n}^Ty} dy.\label{oneside}
\end{align}
Therefore, \begin{equation}\label{oneside1}
L^{-1}e^{\|a_{c_n}\|^{\alpha}} b_{c_n}^{-d}\mathbb P\left(\frac{X-a_{c_n}}{b_{c_n}} \in A\right) \geq  (1-\epsilon)e^{-\epsilon'}\int_{A} e^{-\alpha b_{c_n} \|a_{c_n}\|^{\alpha-2}a_{c_n}^Ty} dy.
\end{equation}
By taking $\liminf$ on both sides as $c_n\to \infty$, and seeing that $\alpha b_{c_n} \|a_{c_n}\|^{\alpha-2}a_{c_n}^Ty \to \theta^T y$ from the definitions of $b_{c_n}$ and $a_{c_n}$, we obtain $$
\liminf_{c_n \to\infty} L^{-1}e^{\|a_{c_n}\|^{\alpha}} b_{c_n}^{-d}\mathbb P\left(\frac{X-a_{c_n}}{b_{c_n}} \in A\right) \geq (1-\epsilon)e^{-\epsilon'}\int_{A} e^{-\theta^T y} dy.
$$
Since this is true for all $\epsilon,\epsilon'>0$ it follows that \begin{equation}\label{oneside2}
\liminf_{c_n \to\infty} L^{-1}e^{\|a_{c_n}\|^{\alpha}} b_{c_n}^{-d}\mathbb P\left(\frac{X-a_{c_n}}{b_{c_n}} \in A\right) \geq \int_{A} e^{-\theta^T y} dy.
\end{equation}
The other side of this equation i.e. \begin{equation}\label{otherside}
\limsup_{c_n \to\infty} L^{-1}e^{\|a_{c_n}\|^{\alpha}} b_{c_n}^{-d}\mathbb P\left(\frac{X-a_{c_n}}{b_{c_n}} \in A\right) \leq \int_{A} e^{-\theta^T y} dy.
\end{equation}
is derived exactly in the same way as \eqref{oneside2}, using the other sides of the inequalities \eqref{est} and \eqref{remtozero}. This is exactly the definition of vague convergence, whence the result follows.
\end{proof}
\end{theorem}

Once this is true, if $f(c_n) = L^{-1}e^{m_{c_n}^{\alpha}} \alpha^dm_{c_n}^{d(\alpha-1)}$, then $$
n \mathbb P\left(\frac{a_{f^{-1}(n)}-X}{b_{f^{-1}(n)}} \in \cdot\right) \to \nu(\cdot)
$$
where $\nu$ has density $-e^{-\theta^T y}$ with respect to the Lebesgue measure.

From here, the proof of Theorem~\ref{empunbdd}(b) follows exactly as the proof of the previous two theorems of this kind. 

\begin{proof}[Proof of Theorem~\ref{empunbdd}(b)]
Suppose that $\frac{M_{c_n \theta}}{n} \to C\in (0,\infty)$. Then, by Theorem~\ref{m2tmt2unbdd} we see that $\frac{c_n}{f^{-1}(n)} \to C_1$ for some $C_1>0$. Let $D$ be a Borel set, and $A_n = a_{f^{-1}(n)} - b_{f^{-1}(n)}D$.

Consider $\mathbb P(R_{n,c_n\theta} \in A_n)$, and divide the top and bottom by $e^{\theta^T a_{f^{-1}(n)}}$ to obtain 
\begin{align*}
\mathbb P(R_{n,c_n\theta} \in A_n) =& \frac{\sum_{i=1}^n e^{c_n\theta^T X_{i}}\1_{X_{i} \in A_n}}{\sum_{i=1}^n e^{c_n\theta^T X_{i}}}\\ =&  \frac{\sum_{i=1}^n e^{c_n\theta^T(X_{i}-a_{f^{-1}(n)})}\1_{X_{i} \in A_n}}{\sum_{i=1}^n e^{c_n\theta^T(X_{i}-a_{f^{-1}(n)})}} \\
=& \Phi_n\left(\sum_{i=1}^n e^{c_n\theta^T(a_{f^{-1}(n)}-X_i)}1_{A_n}, \sum_{i=1}^n e^{c_n\theta^T (a_{f^{-1}(n)}-X_i)}\right),
\end{align*}
where 
$$
\Phi(y,z) = \frac{y}{z}.
$$

At this point, we verify the hypotheses of Lemma~\ref{dekhliyo}. We take $\|\theta_n\| = c_n$ and $\theta$ as in the lemma itself. Let $a_n =a_{f^{-1}(n)}$, $C_1$ be as above, $\nu$ be as in \eqref{ass:mvrv} and $D$ as chosen. All the hypotheses of Lemma~\ref{dekhliyo} are easily verified.

Since $\Phi$ is a continuous mapping, and $\int e^{-C_1\theta^Ty} dPRM(\nu) \neq 0$ with probability $1$ by the definition of a PRM, by Lemma~\ref{dekhliyo} it follows that \begin{align*}
\mathbb P(R_{n,\theta} \in A_n) = & \Phi_n\left(\sum_{i=1}^n e^{c\theta^T(a_{f^{-1}(n)}-X_i)}1_{A_n}, \sum_{i=1}^n e^{c\theta^T(a_{f^{-1}(n)}-X_i)}\right) \\  \overset{d}{\to} & \frac{\int e^{-C_1y}1_{y \in D} d PRM(\nu)}{\int e^{-C_1y} dPRM(\nu)}.
\end{align*}

Finally, note that \begin{align*}
\mathbb P(R_{n,c_n\theta} \in A_n) = & \mathbb P(c_n\theta(a_{f^{-1}(n)} - R_{n,c_n\theta}) \in c_n(a_{f^{-1}(n)} - A_n)) \\ \approx & \mathbb P(c_n\theta(a_{f^{-1}(n)} - R_{n,c_n\theta}) \in C_1D)
\end{align*}
if $n$ is large enough. Combining the two statements above, if $Z$ is a random vector such that
$$
\mathbb P(Z \in D) =  \frac{\int e^{-C_1y}1_{y \in D} d PRM(\nu)}{\int e^{-C_1y} dPRM(\nu)},
$$
then $$
\mathbb P(c_n(a_{f^{-1}(n)}-R_{n,c_n\theta}) \in D) \to \mathbb P(C_1Z\in D).
$$
Thus, the random variable $Z_{C,PRM} = C_1Z$ is the desired limit. Note that this random variable depends upon the Poisson random measure, while the scaling limit $Z$ from part (a) of the theorem does not. It follows that these two random variables are not the same, which concludes the proof.

We remark that, as in the one-dimensional case, the limiting random variable $Z_{C,PRM}$ is continuous, but do not prove this here.
\end{proof}

Note that this part of the theorem, unlike the previous two parts, is slightly different in the following way. Let $X_{(n)}$ denote the sample maximizer of $\theta^T y$ as in the previous sections. It turns out that while $\theta^T X_{(n)}$ scales to a Gumbel random variable, the component $X_{(n)} - \theta^T X_{(n)}$ orthogonal to $\theta$ requires no scaling at all, to converge to a standard normal random vector. We omit the details since they are superfluous, but remark that the estimator $R_{n,\theta_n}$, therefore, doesn't admit a scaling limit in the traditional sense.

All we can show, is that $$R_{n,\theta_n} - X_{(n)} \overset{d}{\to} 0,$$ which follows exactly as in the proof of Theorem~\ref{thm:emphd}(c).
\end{appendices}

\end{document}